\documentclass[a4paper,fleqn]{cas-sc}

\usepackage[numbers]{natbib}

\newtheorem{theorem}{Theorem}
\newtheorem{definition}{Definition}

\newtheorem{lemma}{Lemma}
\newdefinition{remark}{Remark}
\newproof{proof}{Proof}

\begin{document}
\let\WriteBookmarks\relax
\def\floatpagepagefraction{1}
\def\textpagefraction{.001}

\shorttitle{Numerical Analysis of a Nonlinear Thermoelastic System}    

\shortauthors{Alcântara et al.}  

\title [mode = title]{Numerical Analysis for a Nonlinear Coupled Thermoelastic System with Acoustic Boundary Conditions}

\tnotemark[1] 

\tnotetext[1]{} 

%
\author[1]{Adriano A. Alc\^{a}ntara}[orcid=0000-0003-3422-7388]
\cormark[1]
\ead{adriano.alcantara@ppgi.ufrj.br}
\credit{
Formal analysis, 
Conceptualization, 
Methodology, 
Investigation, 
Writing -- original draft, 
Writing -- review and editing}

\author[1]{Bruno A. Carmo}[orcid=0000-0002-1278-2324]
\credit{
Software, 
Methodology, 
Investigation, 
Writing -- original draft, 
Writing -- review and editing}

\author[1]{Juan B. Límaco}[orcid=0000-0001-8695-8612]
\credit{
Supervision, 
Writing -- review and editing}

\author[1]{Carlos M. Guzmán}[orcid=0000-0002-7570-6775]
\credit{
Supervision, 
Writing -- review and editing}

\affiliation[1]{
  organization = {Institute of Mathematics and Statistics, Federal Fluminense University},
  city         = {Niterói},
  state        = {Rio de Janeiro},
  country      = {Brazil}}

\cortext[1]{Corresponding author}

\begin{abstract} 
This paper presents the first $L^2$-norm error analysis for a nonlinear thermoelastic system with acoustic boundary conditions.
The model involves three coupled variables: the acoustic potential $u$, the temperature $\theta$, and the normal displacement to the boundary $\delta$,  with nonlinearities $|u|^d u$ in the domain and $g(\cdot,u')$ on the boundary, in addition to a nonlocal term $\beta(\int_\Omega\theta dx)\Delta\theta$ in the heat equation. 
We establish error estimates for the fully discrete problem, using the Crank–Nicolson method for time discretization and finite element method with basis of degree $r \geqslant 1$ for space discretization. 
We demonstrate that the domain–boundary coupling induces a half-power loss in the spatial convergence order for the variable $u$, and this effect propagates to the remaining variables due to the thermoelastic and boundary couplings, resulting in the final error estimate  $\mathcal{O}(h^{r+\frac{1}{2}} + \tau^2)$. 
A bootstrap argument is employed to ensure the uniform boundedness of $\|U_h^n\|_{L^\infty}$, which is essential for handling the nonlinearity. 
Numerical results using Lagrange finite element basis of degrees $ r = 1, 2, 3 $
(linear, quadratic, and cubic) corroborate the theoretical estimates and illustrate the transition from suboptimal to optimal convergence orders when the exact solution possesses additional $L^\infty$ regularity.
\end{abstract}



\maketitle

\section{Introduction}
\label{sec:introducao}
\subsubsection{Acoustic boundary conditions: historical background and literature review}
The propagation of acoustic waves in fluid media is a classical problem in mathematical physics, traditionally described by the linear wave equation. However, in many real-world applications, the boundary enclosing the domain is not perfectly rigid but exhibits reactive properties that demand more sophisticated models. To capture this fluid–structure interaction, Beale and Rosencrans \cite{Artigo1974_Beale_Acoustic} introduced acoustic boundary conditions, in which the surface reacts locally to the acoustic pressure as a resistive harmonic oscillator, coupling the dynamics of the fluid with that of the boundary.

In this model, the normal displacement to the boundary, denoted by $\delta(x,t)$,  satisfies an ordinary differential equation on the boundary, while the acoustic field $u(x,t)$ obeys the wave equation in the domain. The impenetrability condition ensures compatibility between the fluid velocity and the boundary velocity, resulting in a hybrid system with hyperbolic dynamics in the domain and oscillatory dynamics on the boundary (for small amplitudes). Beale \cite{Artigo1976_Beale_Spectral,Artigo1977_Beale} analyzed the spectral properties of the operator associated with the linear problem with acoustic boundary condition, demonstrating that such a condition significantly alters the spectrum and the asymptotic behavior of the solutions.

Since then, a vast literature has consolidated the mathematical analysis of hyperbolic equations with acoustic boundary conditions. Notable contributions include the studies by Frota and Goldstein \cite{frota2000some} and Cousin, Frota and Larkin \cite{cousin2001system} on nonlinear problems; the works of Mugnolo \cite{mugnolo2006abstract}, Kobayashi and Tanaka \cite{kobayashi2008application}, Frigeri \cite{frigeri2010attractors}, and Limaco et al. \cite{limaco2011evolution}; and, more recently, results on stabilization and energy decay obtained by Graber \cite{graber2012uniform}, Vicente and Frota \cite{vicente2013mixed,vicente2013nonlinear}, Ha \cite{Artigo2016_Tae_Gab_Ha}, and Cavalcanti et al. \cite{Artigo2020_Cavalcanti}, as well as the recent works by Vitillaro \cite{Enzo2024,Enzo2025} and the references therein. In addition, acoustic boundary conditions also appear in various other models of interest, including the Klein–Gordon system studied in \cite{cousin2004system} and the quasilinear von Kármán equation with memory effects recently investigated in \cite{Artigo2021_LEE}.

\subsubsection{Linear and nonlinear models with an acoustic boundary: evolution of numerical analysis}
Numerical analysis for this class of problems is quite recent. To the best of our knowledge, the first paper to present a numerical simulation for a linear wave equation with variable coefficients and acoustic boundary conditions was published by Alcântara et al. \cite{Artigo2018_Adriano} in 2018. In the same year, Hipp et al. \cite{Artigo2018_HippKovacs2018_H1} published the pioneering numerical analysis of a linear wave equation with a constant coefficient and acoustic boundary conditions, establishing the first error estimates in the $H^1$-norm. Building on this work, Hipp and Kovács \cite{Artigo2020_HippKovacs_L2} extended the analysis to the $L^2$-norm in 2020, using linear Lagrange finite elements. Subsequently, nonlinearities were incorporated both in the domain and on the boundary by Alcântara et al. \cite{Artigo2021_Adriano,Artigo2022_Adriano}, who also investigated existence, uniqueness, exponential decay, and numerical simulations using linear and quadratic Lagrange basis. Leibold \cite{Artigo2022_Jan} extended the error estimation framework of \cite{Artigo2018_HippKovacs2018_H1} to include a boundary nonlinearity, while maintaining the $H^1$-norm estimate and the linear wave equation in the domain. Recently, in 2025, Alcântara et al. \cite{Artigo2025_Adriano} consolidated the $L^2$-norm numerical analysis for three wave models with acoustic boundary conditions, including cubic finite element basis, and highlighted a half-order loss in the spatial convergence rate of the acoustic variable due to the domain--boundary coupling. As far as we know, this was the first work in the literature to establish $L^2$-norm error estimates for wave equations with general nonlinearities both in the domain and on the boundary.

\subsubsection{The study model}
Although the works mentioned earlier represent significant advances in the numerical analysis of wave equations with acoustic boundary conditions, they are all restricted to purely hyperbolic models and do not take thermal effects into account. In this work, we focus on a more general and physically richer model: the thermoelastic system with acoustic boundary conditions proposed by Braz e Silva, Clark and Frota \cite{Artigo2017_SilvaClarkFrota}. This system couples the acoustic wave dynamics to a nonlinear and nonlocal heat equation, and it is the model for which we develop the numerical analysis in the forthcoming sections. More precisely, the model under study is
\begin{align}
\left\{
\begin{aligned}
&	u''(x,t)
	- \mu(t) \Delta u(x,t) 
	+ f(u(x,t))
	+ (\mathbf{a}\cdot\nabla) \theta(x,t)
	= 0\,\, \mbox{ in } \, \, Q,  
\\
&	\theta'(x,t) 
	- \beta\Big(\int_{\Omega}\theta(x,t)\,dx\Big)\Delta \theta(x,t)
	+(\mathbf{a}\cdot\nabla) u'(x,t)
	= 0 \,\mbox{ in } \, \,Q,
\\
&	u(x,t) = 0\; \mbox{ on }\; \Sigma_0, 
\\
&	u'(x,t) 
	+ f_1(x)\delta''(x,t) 
	+ f_2(x)\delta'(x,t) 
	+ f_3(x)\delta(x,t)=0 \,\mbox{ on } \,\Sigma_1, 
\\
&	\frac{\partial u }{\partial \nu}(x,t) 
	- \delta'(x,t) 
	+ g(x, u'(x,t))
	= 0 \; \mbox{ on }\,\Sigma_1, 
\\
&	\theta(x,t)=0\;\mbox{ on }\; \Sigma, 
\\
&	u(x,0) = u_{0}(x), \quad  u^{\prime}(x,0) = u_{1}(x) \quad
\mbox{ for } \,\, x \in  \Omega,
\\
&	\theta(x,0) = \theta_{0}(x) \quad  \mbox{ for } \,\, x \in  \Omega,
\\
&	\delta(x,0) = \delta_{0}(x), \quad \delta'(x,0)=\frac{\partial u_0 }{\partial
\nu}(x) - g(x, u_1(x)) \, \, \mbox{ for } \,\, x \in \Gamma_{1},
\end{aligned}
\right.\label{modelo1}
\end{align}
where  $f(u)=\lambda|u|^d u$, with $\lambda, d>0$ constants;
 $u_{0}, u_{1}, \theta_{0} : \Omega \rightarrow \mathbb{R}$ and
$\delta_{0}: \Gamma_{1} \rightarrow \mathbb{R}$ are given functions. The domain $\Omega$ is a bounded open subset of $\mathbb{R}^n$ with smooth boundary $\Gamma:=\Gamma_0\cup\Gamma_1$ such that $\Gamma_0\cap\Gamma_1=\emptyset$, $Q:=\Omega\times(0,T)$ with $T>0$, $\Sigma_0:=\Gamma_0\times(0,T)$, $\Sigma_1:=\Gamma_1\times(0,T)$, and $\Sigma=\Sigma_0\cup\Sigma_1$. Furthermore,
$\mu:[0,\infty) \rightarrow \mathbb{R}$ and $\beta:\mathbb{R}
\rightarrow \mathbb{R}$ are given functions; $\mathbf{a}= (a_{1}, \ldots, a_{n})$ 
is a constant vector in $\mathbb{R}^n$; $(\mathbf{a}\cdot\nabla)$ is the operator
$\sum\limits_{i=1}^{n} a_{i}\frac{\partial}{\partial x_{i}}$ and $\Delta
= \sum\limits_{i=1}^{n} \frac{\partial^{2}}{\partial x_{i}^{2}}$ is the Laplace operator;
$f_i$ para $i=1,2,3$, are given real-valued functions defined on $\Gamma_1$;  $\nu$  denotes the outward unit normal vector on $\Gamma_1$;
$g: \Gamma_1 \times \mathbb{R} \to \mathbb{R}$ is a given function. Here, ${^{\prime}}$ denotes the partial derivative with respect to time, and all derivatives are understood in the sense of distributions.

The physical interpretation of system \eqref{modelo1} is as follows. As in \cite{Artigo1974_Beale_Acoustic}, the domain $\Omega$ represents a region occupied by a fluid (gas or liquid) at rest, except for small-amplitude acoustic wave motions. 
The variable $u$ is the acoustic velocity potential in \eqref{modelo1}$_{1}$, $\theta$ is the thermal variation in \eqref{modelo1}$_{2}$, and $\delta$ is the normal displacement to the boundary in \eqref{modelo1}$_{4}$. 
The impenetrability condition in \eqref{modelo1}$_{5}$ couples the normal fluid velocity to the boundary velocity, with an additional dissipative term $g(u')$. 
The boundary is decomposed into two parts with distinct physical roles: $\Gamma_0$, the rigid portion of the boundary, where the classical Dirichlet condition $u=0$ is imposed; and $\Gamma_1$, the reactive portion, where each point behaves as a mass-spring-damper system. The fundamental novelty lies in the bidirectional thermoelastic coupling: on one hand, thermal gradients generate acoustic waves via the term $(\mathbf{a}\cdot\nabla)\theta$ in the equation for $u$; on the other, fluid motion alters the temperature field via the term $(\mathbf{a}\cdot\nabla)u'$ in the heat equation, characterizing the inverse thermoacoustic effect (see Corato et al., \cite{Artigo2025_Corato}). Furthermore, the nonlocal diffusion coefficient $\beta(\int_\Omega \theta\,dx)\Delta\theta$ generalizes the classical model by introducing a dependence on the average temperature within the thermal diffusion process.

\subsubsection{Objective, originality and contributions}
\label{subsec:objetivo}
Although system \eqref{modelo1} has been theoretically analyzed in \cite{Artigo2017_SilvaClarkFrota}, where existence, uniqueness, and exponential decay of solutions were established, there are, to the best of our knowledge, no numerical analysis results available for this model. This gap in the literature is particularly relevant because the coupling between the acoustic wave dynamics and the nonlocal heat equation introduces additional difficulties—such as the domain–boundary coupling via the trace operator and the presence of nonlinearities in both the domain and the boundary—that are not covered by existing numerical analyses for wave equations with acoustic boundary conditions.

The central objective of this work is therefore to develop, for the first time, a complete numerical analysis of the nonlinear thermoelastic system \eqref{modelo1} with acoustic boundary conditions. More precisely, we establish $L^2$-norm error estimates for the fully discrete problem obtained by combining the finite element method in space with the Crank–Nicolson scheme in time. The originality of our approach lies in extending the classical techniques for wave equations with acoustic boundary conditions to the thermoelastic setting, where the coupling with the heat equation and the nonlinearities require new arguments. The main contributions of this work are summarized below.

\begin{itemize}
    \item 
    \textit{First $L^2$-norm error estimates for the fully discrete thermoelastic system.}
    We derive error estimates in the $L^2$ norm for the fully discrete problem associated with \eqref{modelo1}, covering the three coupled variables $u$, $\theta$, and $\delta$, as well as their first-order time derivatives. To the best of our knowledge, this is the first numerical analysis of this model.

    \item 
    \textit{Identification of the half-order loss in the spatial convergence.}
    We show that the domain–boundary coupling induces a half-order loss in the spatial convergence of $u$, arising from the trace estimate
    \[
    \|\gamma_0(R_h^u u - u)\|_{L^2(\Gamma_1)} \leqslant C\, h^{r+\frac{1}{2}}\|u\|_{H^{r+1}(\Omega)}.
    \]
    This loss propagates to the final estimate and reduces the spatial convergence order to $\mathcal{O}(h^{r+\frac{1}{2}})$. This result extends the error estimates established in \cite{Artigo2015_Kashiwabara} and \cite{Artigo2025_Adriano} to the context of the nonlinear thermoelastic system with acoustic boundary conditions.

    \item 
    \textit{Treatment of the nonlinearity $|u|^d u$ via a bootstrap argument.}
    We handle the nonlinear term $f(u) = \lambda|u|^d u$ through a bootstrap  argument that guarantees constants independent of the mesh size $h$ and the time step $\tau$. This is essential because $f$ is not globally Lipschitz.

    \item 
    \textit{Generalization of the order-reduction strategy of Baker \cite{Artigo1976_Baker}.}
    We extend the order-reduction technique employed by Baker for linear wave equations to a substantially more general setting. While Baker transformed a linear hyperbolic problem into a first-order parabolic system to apply the Crank–Nicolson method, our system involves nonlinearities in the domain ($\lambda|u|^d u$) and on the boundary ($g(\cdot,u')$), as well as a thermoelastic coupling that introduces a nonlocal heat equation and a system of three coupled variables ($u$, $\theta$, and $\delta$). This requires a careful treatment of the interactions between the three variables and of the nonlocal term.

    \item 
    \textit{Numerical simulations in two dimensions with Lagrange finite elements.}
    We perform two-dimensional numerical simulations using Lagrange finite elements of degree $r \geqslant 1$ in space—including linear, quadratic, and cubic basis functions—combined with the Crank–Nicolson method in time. At each discrete time step, Newton's method is employed to solve the resulting nonlinear system. The numerical results corroborate the theoretical estimates and evidence the transition from suboptimal to optimal orders when the exact solution possesses additional $L^\infty$ regularity.

    \item 
    \textit{Reproducibility of the numerical results.}
    The numerical examples were implemented in the open-source programming language Julia. The complete source code, together with the instructions for reproducing all numerical results presented in this work, is available in the reproducibility repository \cite{ThermoelasticAcoustic}. This ensures full reproducibility and allows other researchers to easily extend the present numerical framework to related problems.
\end{itemize}

In summary, this work extends the $L^2$-error estimates of \cite{Artigo1976_Baker}, \cite{Artigo2015_Kashiwabara}, \cite{Artigo2020_HippKovacs_L2}, and \cite{Artigo2025_Adriano}, and the numerical methodology of \cite{Artigo2018_Adriano, Artigo2021_Adriano, Artigo2022_Adriano, Artigo2025_Adriano}, to a nonlinear thermoelastic system with acoustic boundary conditions. Compared with these references, the analysis and the simulation are more demanding: the thermoelastic coupling forces the simultaneous treatment of $u$, $\theta$, and $\delta$, and the acoustic boundary condition induces a half-order loss $\mathcal{O}(h^{r+\frac{1}{2}})$ that propagates to all variables. Numerically, the Crank--Nicolson scheme yields a nonlinear algebraic system solved by Newton's method, and reproducing the suboptimal-to-optimal transition requires exact solutions of controlled $L^\infty$ regularity and careful tuning of the Lagrange bases, well beyond the experiments of the cited works.

The paper is organized as follows. Section 2 presents the preliminary results: function spaces, notations, assumptions, and the existence and uniqueness theorem for the solution. Section 3 introduces the finite element spaces, the Ritz projections and their convergence and stability results, and 
develops the fully discrete error analysis in the $L^2(\Omega)$ and 
$L^2(\Gamma_1)$ norms, culminating in the main convergence theorem 
with the estimate $\mathcal{O}(h^{r+\frac{1}{2}} + \tau^2)$. 
Section 4 illustrates the theoretical convergence results numerically, highlighting the suboptimal and optimal cases. 
Finally, Section 5 presents the conclusions.
\section{Definitions and preliminary results}
\label{sec:var}
In addition to the classical function spaces (see \cite{Livro1969_Lions}), we establish the notation and the function spaces that will be used throughout the work. We define the following spaces:
\begin{itemize}
    \item $V = H_{\Gamma_0}^1(\Omega) = \{ u \in H^1(\Omega) : u|_{\Gamma_0} = 0 \}$, which is the space of functions satisfying the Dirichlet condition for $u$, with the norm $\|u\|_{V}^2 = \|\nabla u\|_{L^2(\Omega)}^2$.

    \item $H^1_0(\Omega) = \{ v \in H^1(\Omega) : v|_{\Gamma} = 0 \}$, which is the space of functions satisfying the Dirichlet condition for $\theta$; with norm $\|v\|_{H^1_0(\Omega)}^2 = \|\nabla v\|_{L^2(\Omega)}^2$.
    \item $H_\Delta(\Omega) = \{ w \in H^1(\Omega) : \Delta w \in L^2(\Omega) \}$, with the norm $
    \|w\|_{H_\Delta(\Omega)}^2 = \|w\|_{H^1(\Omega)}^2 + \|\Delta w\|_{L^2(\Omega)}^2.$
\end{itemize}
It is well known that the Poincaré inequality holds in both $V$ and $H^1_0(\Omega)$, and that the norms of $V$ and $H^1_0(\Omega)$ are equivalent to the $H^1(\Omega)$-norm. The trace operators are respectively defined as (zero and one trace)
  $
  \gamma_0: H^{1}(\Omega) \rightarrow H^{1/2}(\Gamma)\;\mbox{and}\;  \gamma_1: H_\Delta(\Omega)\rightarrow H^{-1/2}(\Gamma),
  $
where $\gamma_0(u)=u|_\Gamma$ and $\gamma_1(u)=\frac{\partial u}{\partial \nu}|_\Gamma$ for all $u \in \mathcal{D}(\overline{\Omega})$.  

For $T > 0$ and $X$ a Banach space, we denote by $L^p(0,T;X)$ the space of functions $\xi : (0,T) \to X$ such that $\|\xi\|_{L^p(0,T;X)} < \infty$, with
$\|\xi\|_{L^p(0,T;X)} = \big( \int_0^T \|\xi(t)\|_X^p \, dt \big)^{1/p}, \; 1\leqslant p < \infty,$
and $\|\xi\|_{L^\infty(0,T;X)} = \operatorname*{ess\,sup}_{0<t<T} \|\xi(t)\|_X$. 
We also define the Sobolev spaces with values in $X$,
\begin{align*}
H^m(0,T;X) = \left\{ \xi \in L^2(0,T;X) : \partial_t^i \xi \in L^2(0,T;X), \; i = 1, \ldots, m \right\},
\end{align*}
with norm $\|\xi\|_{H^m(0,T;X)}^2 = \sum_{i=0}^m \|\partial_t^i \xi\|_{L^2(0,T;X)}^2$.

From now on, we denote by $(\cdot , \cdot)$, $(\cdot , \cdot)_{\Gamma_1}$, $\|\cdot\|$, $\|\cdot\|_{\Gamma_1}$ the inner products and induced norms in the spaces $L^2(\Omega)$ and $L^2(\Gamma_1)$, respectively. When convenient, we will simplify the notation for the spaces $L^p(0,T;X)$,  writing simply $L^p(X)$. Moreover, throughout the text, $\|.\|_{\infty}$
denotes the maximum norm. 

The concept of global solution for model \eqref{modelo1} is given by the following definition.

\begin{definition}\label{D2.1}
A global solution to the nonlinear initial-boundary value problem \eqref{modelo1} is a triple $\{u, \theta, \delta\}$ in the class
\begin{equation}\label{2.1}
\begin{array}{l}
u \in L^\infty_{\rm loc}\big(0,\infty; \, V\cap H_\Delta(\Omega)\big),\,
u' \in L^\infty_{\rm loc}\big(0,\infty; \, V \big),\,
u'' \in L^\infty_{\rm loc}\big(0,\infty;\, L^2(\Omega)\big),
\\[2pt]
\theta,\; \theta' \in
L^{\infty}_{\rm loc}(0,\infty;L^2(\Omega))\cap
L^{2}_{\rm loc}(0,\infty;H^1_0(\Omega)),
\\[2pt]
\delta,\;\delta',\;\delta'' \in L^\infty_{\rm loc}\big(0,\infty; \, L^2(\Gamma_1)\big),
\end{array}
\end{equation}
such that for every arbitrarily fixed $T > 0$, the following equations are satisfied:
\begin{equation}\label{2.2}
\begin{array}{l}
u'' - \mu \Delta u+\lambda |u|^d u+(\mathbf{a} \cdot \nabla)\theta = 0 \quad \text{in } \quad L^\infty\big(0,T;L^2(\Omega)\big),
\\[4pt]
\dfrac{\partial u}{\partial \nu}-\delta'+g(\cdot,u')=0 \quad \text{in } \quad L^\infty\big(0,T;H^{-1/2}(\Gamma_1)\big),
\\[6pt]
\theta'-\beta\big(\int_{\Omega}\theta\,dx \big)\Delta \theta + (\mathbf{a} \cdot \nabla) u' = 0 \quad \text{in } \quad L^\infty\big(0,T;L^2(\Omega)\big),
\\[8pt]
u'+f_1\delta''+f_2\delta'+f_3\delta = 0 \quad \text{in } \quad L^\infty\big(0,T;L^{2}(\Gamma_1)\big),
\end{array}
\end{equation}
and the initial conditions \eqref{modelo1}$_{7,8,9}$ are also satisfied.
\end{definition}
\begin{theorem}(Existence and uniqueness).\label{T2.1}
Assume that:
\begin{itemize}
\item[a)] $\mathbf{a} = (a_1, \ldots, a_n) \in \mathbb{R}^n$ is a vector and $\lambda > 0$ is a real number;

\item[b)] $d \in \mathbb{R}$ with $0 < d < \infty$ if $n = 1$; $1/2 \leqslant d \leqslant \infty$ if $n = 2$; or $1/n \leqslant d \leqslant 2/(n-2)$ if $n \geqslant 3$;

\item[c)] $\mu \in L^{\infty}_{\rm loc}(\mathbb{R}^+)$, $\mu' \in L^1(\mathbb{R}^+)$ \,and\, $\mu(t) \geqslant \mu_0 > 0$ \; for almost every
\; $t \in \mathbb{R}^+$;

\item[d)] $\beta$ is differentiable with $\beta' \in L^{\infty}_{\rm loc}(\mathbb{R})$ and $\beta(\vartheta) \geqslant \beta_0 > 0$;

\item[e)] $f_i \in C^0(\Gamma_1)$, $i=1,2,3$, with $f_i(x) > 0$ for \,$i=1,3$\, and $f_2(x) \geqslant 0$\; for all\; $x \in \Gamma_1$;

\item[f)] The function $g:\Gamma_1\times\mathbb{R}\to\mathbb{R}$ is  a continuous function satisfying:
\begin{align}
|g(x,s)-g(x,r)| &\leqslant L_g|s-r|, \quad \forall x\in\Gamma_1,\; \forall s,r\in\mathbb{R}, \label{hip:eta_lip}\\
[g(x,s)-g(x,r)](s-r) &\geqslant g_0(s-r)^2, \quad g_0>0,\quad \text{ and } \quad g(x,0)=0, \; \forall x \in \Gamma_1.  \label{hip:eta_mon}
\end{align}
\end{itemize}
Thus, there exists a unique solution to the coupled nonlinear thermoelastic model \eqref{modelo1} in the sense of Definition \ref{D2.1}, taking the initial data in the spaces
\[
(u_0, u_1, \theta_0, \delta_0) \in (V\cap H_\Delta(\Omega)) \times V \times H_0^1(\Omega) \times L^2(\Gamma_1).
\]
\end{theorem}
\noindent $\mathbf{ Proof.}$ See \cite{Artigo2017_SilvaClarkFrota}.
\hfill$\square$

\section{$L^2$-norm error estimate}
\subsection{Spatial discretization: finite element spaces and Ritz projections}
\label{sec:estimativa}
As in \cite{Artigo2015_Kashiwabara}, for the error estimation in this work,
let $\Omega\subset\mathbb{R}^n$ ($n=2$ or $3$) and $\mathcal{T}_h$ be a quasi-uniform triangulation of the domain $\overline{\Omega}$ in finite elements satisfying the following properties: $\overline{\Omega} = \bigcup_{E \in \mathcal{T}_h} E$; for any distinct $E, F \in \mathcal{T}_h$, $int E \cap int F = \varnothing$; triangulation is \textit{regular} (or \textit{shape-regular}): there is a constant $\sigma > 0$ such that, for all $E \in \mathcal{T}_h$,
$\frac{h_E}{\rho_E} \le \sigma,$ where $h_E = \operatorname{diam}(E)$ and $\rho_E$ is the diameter of the largest ball inscribed in $E$. The global parameter $h = \max_{E\in\mathcal{T}_h} h_E$.
For the boundary $\Gamma_1$, the induced triangulation is defined
$
\mathcal{T}_h|_{\Gamma_1} = \big\{ B = \partial E \cap \Gamma_1 : E \in \mathcal{T}_h,\; \partial E \cap \Gamma_1 \neq \varnothing \big\}.$

For a fixed degree $j \geqslant 1$, we define the finite element spaces: $V_h = \big\{ v_h \in C(\overline{\Omega}) : v_h|_E \in \mathcal{P}_j(E), \;\forall E\in\mathcal{T}_h,\; v_h|_{\Gamma_0}=0 \big\} \subset V$; \, $W_h = \big\{ \theta_h \in C(\overline{\Omega}) : \theta_h|_E \in \mathcal{P}_j(E), \;\forall E\in\mathcal{T}_h,\; \theta_h|_{\Gamma}=0 \big\} \subset H^1_0(\Omega)$; and 
$Z_h = \big\{ \chi_h \in C(\overline{\Gamma}_1) : \chi_h|_B \in \mathcal{P}_j(B), \;\forall B \in \mathcal{T}_h|_{\Gamma_1} \big\} \subset L^2(\Gamma_1).$  Here, $\mathcal{P}_j(E)$ denotes the space of polynomials of total degree at most  $\leqslant j$ on the element $E$ and analogously on the element $B$.

\begin{definition}[Ritz projections]
\label{def:ritz}
For each $t \in [0, T]$, as in \cite{Artigo2015_Kashiwabara} and \cite{Artigo2025_Adriano}, we define the projection operators:
\begin{itemize}
    \item $R_h^u: V\to V_h$ such that $(\nabla R_h^u v(t),\nabla w_h)=(\nabla v(t),\nabla w_h)$, \; for all\; $w_h\in V_h$.
    \item $R_h^\theta: H^1_0(\Omega)\to W_h$ such that $(\nabla R_h^\theta\varphi(t),\nabla\psi_h)=(\nabla\varphi(t),\nabla\psi_h)$, \; for all\;
    $\psi_h\in W_h$.
    \item $R_h^\delta: L^2(\Gamma_1)\to Z_h$ such that $(R_h^\delta\zeta(t),\chi_h)_{\Gamma_1}=(\zeta(t),\chi_h)_{\Gamma_1}$, \; for all\; $\chi_h\in Z_h$.
\end{itemize}
\end{definition}
The existence and uniqueness of such operators follow from the application of the Lax-Milgram Theorem; see \cite{Livro2010_Evans}.

\begin{lemma}[Stability of Ritz projections]
\label{lem:estabilidade_ritz}
The Ritz projections in Definition \ref{def:ritz} satisfy, for each $t\in[0,T]$:
\begin{align}
\|\nabla R_h^u v(t)\| &\leqslant \|\nabla v(t)\|, \; \text{ for all } \; v(t)\in V, \label{eq:estab_u}\\
\|\nabla R_h^\theta \varphi(t)\| &\leqslant \|\nabla \varphi(t)\|, \; \text{ for all } \; \varphi(t)\in H^1_0(\Omega), \label{eq:estab_theta}\\
\|R_h^\delta \zeta(t)\|_{\Gamma_1} &\leqslant \|\zeta(t)\|_{\Gamma_1}, \; \text{ for all } \; \zeta(t)\in L^2(\Gamma_1). \label{eq:estab_delta}
\end{align}
\end{lemma}
\noindent $\mathbf{ Proof.}$
The proof is standard: one chooses the test function equal to the projection itself in Definition \ref{def:ritz}, and uses the Cauchy-Schwarz inequality.
\hfill$\square$

\begin{lemma}[Projection estimates]
\label{lem:proj}
Suppose that $u(t)\in H^{r+1}(\Omega)\cap V$, $\theta(t)\in H^{r+1}(\Omega)\cap H^1_0(\Omega)$, $\delta(t)\in H^{r+1}(\Gamma_1)$ for some $1\leqslant r\leqslant j$. Then there exist constants $C>0$ independent of $h$ such that:
\begin{align}
\|R_h^u u(t)-u(t)\| + h\|\nabla(R_h^u u(t)-u(t))\| &\leqslant C h^{r+1}\|u(t)\|_{H^{r+1}(\Omega)}, \label{est:proj_u}\\
\|R_h^\theta\theta(t)-\theta(t)\| + h\|\nabla(R_h^\theta\theta(t)-\theta(t))\| &\leqslant C h^{r+1}\|\theta(t)\|_{H^{r+1}(\Omega)}, \label{est:proj_theta}\\
\|R_h^\delta\delta(t)-\delta(t)\|_{\Gamma_1} &\leqslant C h^{r+1}\|\delta(t)\|_{H^{r+1}(\Gamma_1)}. \label{est:proj_delta}
\end{align}
\end{lemma}
\noindent $\mathbf{ Proof.}$ See \cite{Artigo2025_Adriano} and \cite{Artigo2015_Kashiwabara}.
An important consequence of applying \eqref{est:proj_u} relating boundary-domain estimation is
\begin{align}
\|\gamma_0(R_h^u u(t) - u(t))\|_{\Gamma_1} \leqslant C h^{r+\frac{1}{2}}\|u(t)\|_{H^{r+1}(\Omega)}. \label{est:trace}
\end{align}
Indeed, for $w \in V$, we have
$\|\gamma_0(w)\|_{\Gamma_1} \leqslant C\|\nabla w\|^{\frac{1}{2}}\|w\|^{\frac{1}{2}}.$
For all details, see \cite{Artigo2025_Adriano} and \cite{Artigo2020_HippKovacs_L2}.
\hfill$\square$

\subsection{Fully discrete problem}
\label{sec:total}
Let $T>0$ be the final time and $\tau\in(0,T)$ the time step.
We consider the uniform temporal grid $t_k = k\tau$, $k = 0,\dots,N$,
with $N$ the integer part of $T/\tau$,
and denote by 
$t_{k-\frac{1}{2}} = \frac{t_{k-1}+t_k}{2}$, $k=1,\dots,N$, the midpoint of each subinterval.
For a function $v$, we set:
\begin{align}
v^k = v(t_k),\quad \widetilde{v}^{\,k} = \frac{v^k+v^{k-1}}{2},\quad \widetilde{\partial}v^k = \frac{v^k-v^{k-1}}{\tau}.
\end{align}

To formulate the fully discrete problem, we take the variational formulation of problem \eqref{modelo1} and, following \cite{Artigo1976_Baker}, transform the second-order equations into a first-order system by introducing $p=u'$ and $q=\delta'$. As in \cite{Artigo1976_Baker}, we obtain the parabolic formulation:
\begin{align}
\left\{
\begin{aligned}
& (p'(t),\phi) + \mu(t)\big[(\nabla u(t),\nabla\phi) - (q(t)-g(p(t)),\phi)_{\Gamma_1}\big]
 + (f(u(t)),\phi) + ((\mathbf{a}\cdot\nabla)\theta(t),\phi) = 0, \; \phi\in V;
 \\
& (\theta'(t),\psi) + \beta\big(\hspace{-0.09cm}\int_\Omega\theta(t)dx\big)(\nabla\theta(t),\nabla\psi) + ((\mathbf{a}\cdot\nabla)p(t),\psi) = 0, \; \psi\in H^1_0(\Omega);
\\
& (p(t) + f_1q'(t) + f_2q(t) + f_3\delta(t),\chi)_{\Gamma_1}  = 0, \; \chi\in L^2(\Gamma_1);
\\
& p(t)=u'(t),\quad q(t)=\delta'(t). 
\end{aligned}
\right.\label{proj-1}
\end{align}
Evaluating equations \eqref{proj-1} at the midpoint $t=t_{k-\frac{1}{2}}$ and using Taylor expansions with integral remainder, as in \cite{Artigo2004_ErvinHeuer} and \citep{Artigo2025_Adriano}, and subsequently neglecting the truncation terms, we obtain the fully discrete problem: for each
$k=1,\dots,N$, find $U_h^k,P_h^k\in V_h$; \,$\Theta_h^k\in W_h$; \,$D_h^k,Q_h^k\in Z_h$ such that for all $\phi_h\in V_h$, $\psi_h\in W_h$, $\chi_h\in Z_h$:
\begin{align}
\left\{
\begin{aligned}
& (\widetilde{\partial}P_h^k,\phi_h) \hspace{-0.07cm}+\hspace{-0.07cm} \mu^{k-\frac{1}{2}}(\nabla\widetilde{U}_h^{\,k},\nabla\phi_h) \hspace{-0.07cm}+\hspace{-0.07cm} (f(\widetilde{U}_h^{\,k}),\phi_h) \hspace{-0.07cm}+\hspace{-0.07cm} ((\mathbf{a}\cdot\nabla)\widetilde{\Theta}_h^{\,k},\phi_h) 
=\hspace{-0.07cm} \mu^{k-\frac{1}{2}}(\widetilde{Q}_h^{\,k},\phi_h)_{\Gamma_1} 
\hspace{-0.07cm}-\hspace{-0.07cm} \mu^{k-\frac{1}{2}}(g(\widetilde{P}_h^{\,k}),\phi_h)_{\Gamma_1}, 
\\
& (\widetilde{\partial}\Theta_h^k,\psi_h) + \beta\big(\hspace{-0.09cm}\int_\Omega\widetilde{\Theta}_h^{\,k}dx\big)(\nabla\widetilde{\Theta}_h^{\,k},\nabla\psi_h) + ((\mathbf{a}\cdot\nabla)\widetilde{P}_h^{\,k},\psi_h) = 0, 
\\
& (\widetilde{P}_h^{\,k} 
+ f_1\widetilde{\partial}Q_h^k 
+ f_2\widetilde{Q}_h^{\,k} 
+ f_3\widetilde{D}_h^{\,k},\chi_h)_{\Gamma_1} 
= 0, 
\end{aligned}         
\right.\label{PTD}
\end{align}
with the consistency relations $\widetilde{\partial}U_h^k = \widetilde{P}_h^{\,k}$, $\widetilde{\partial}D_h^k = \widetilde{Q}_h^{\,k}$;
and the initial conditions $U_h^0$, $P_h^0$, $\Theta_h^0$, $D_h^0$, $Q_h^0$ are obtained by Ritz projections of the continuous initial data onto the respective finite element spaces.

\subsection{Error decomposition}
\label{sec:decomp}
Following \cite{Artigo1973_Wheeler, Artigo1976_Baker}, we define the \textbf{projection errors} as the difference between the projection of the exact solution onto each finite element space and the exact solution itself:
\begin{align}
\begin{aligned}
\omega_u^k &= R_h^u u^k - u^k; &
\omega_\theta^k &= R_h^\theta \theta^k - \theta^k; &
\omega_\delta^k &= R_h^\delta \delta^k - \delta^k; \\
\omega_u'^k &= (R_h^u u)'^k - u'^k; &
\omega_\delta'^k &= (R_h^\delta \delta)'^k - \delta'^k; &
\omega_\delta''^k &= (R_h^\delta \delta)''^k - \delta''^k;
\label{decomp-1}
\end{aligned}
\end{align}
where the superscript $k$ denotes evaluation at $t_k = k\tau$ and the prime indicates differentiation with respect to time.
Let $U_h^k$, $P_h^k$, $\Theta_h^k$, $D_h^k$, and $Q_h^k$ be the numerical approximations (Crank–Nicolson–Galerkin) for $u$, $u'$, $\theta$, $\delta$, and $\delta'$, respectively. We define the \textbf{discretization errors} as the difference between the numerical solution and the projection of the exact solution:
\begin{align}
\begin{aligned}
\rho_u^k = U_h^k - R_h^u u^k; \quad 
\rho_\theta^k = \Theta_h^k - R_h^\theta \theta^k; \quad 
\rho_\delta^k = D_h^k - R_h^\delta \delta^k; \quad
\rho_u'^k = P_h^k - (R_h^u u)'^k; \quad
\rho_\delta'^k = Q_h^k - (R_h^\delta \delta)'^k.
\label{decomp-2}
\end{aligned}
\end{align}
With these definitions, the total errors decompose as
$
U_h^k - u^k = \rho_u^k + \omega_u^k;\;
\Theta_h^k - \theta^k = \rho_\theta^k + \omega_\theta^k;\;
D_h^k - \delta^k = \rho_\delta^k + \omega_\delta^k,
$
and analogously for the respective time derivatives. This decomposition allows the error estimates to be obtained separately:
the projection errors $\omega$ are controlled by Ritz estimates (Lemma \ref{lem:proj}), which depend on the regularity of the exact solution, whereas the discretization errors $\rho$ constitute the central object of the error analysis in this work.

\subsubsection{Projection equation for $u$}

As in \cite{Livro2006_Thomee}, we start from the continuous variational equation at the midpoint; then, after adding and subtracting the projection terms, except for $f(u(t_{k-\frac{1}{2}})$, we obtain, for all $\phi_h\in V_h$,
\begin{align}
& (\widetilde{\partial}(R_h^u u)'^k,\phi_h) 
+ \mu^{k-\frac{1}{2}}(\nabla\widetilde{R_h^u u}^{\,k},\nabla\phi_h) 
+ \mu^{k-\frac{1}{2}}\big(g(\widetilde{R_h^u u}'^k),\phi_h\big)_{\Gamma_1}
+ ((\mathbf{a}\cdot\nabla)\widetilde{R_h^\theta\theta}^{\,k},\phi_h) + \mu^{k-\frac{1}{2}}(\widetilde{R_h^\delta\delta}'^k,\phi_h)_{\Gamma_1}\nonumber
\\
&\;
= - (f(u(t_{k-\frac{1}{2}})),\phi_h)  
+ (\widetilde{\partial}(R_h^u u)'^k - u''(t_{k-\frac{1}{2}}),\phi_h)
+ \mu^{k-\frac{1}{2}}(\nabla(\widetilde{R_h^u u}^{\,k} -  u(t_{k-\frac{1}{2}})),\nabla\phi_h) 
\nonumber
\\
&\quad
+ \mu^{k-\frac{1}{2}}\big(g(\widetilde{R_h^u u}'^k) \hspace{-0.07cm}-\hspace{-0.07cm} g(u'(t_{k-\frac{1}{2}})),\phi_h\big)_{\Gamma_1} 
\hspace{-0.07cm}+\hspace{-0.07cm} \big((\mathbf{a}\cdot\nabla)(\widetilde{R_h^\theta\theta}^{\,k} \hspace{-0.07cm}-\hspace{-0.07cm} \theta(t_{k-\frac{1}{2}})),\phi_h\big) 
\hspace{-0.07cm}+\hspace{-0.07cm} \mu^{k-\frac{1}{2}}\big(\widetilde{R_h^\delta\delta}'^k \hspace{-0.07cm}-\hspace{-0.07cm} \delta'(t_{k-\frac{1}{2}}),\phi_h\big)_{\Gamma_1}. 
\label{eq:proj_u1}
\end{align}
We can rewrite \eqref{eq:proj_u1} and highlight the truncation error part via Taylor expansion; to this end:
\begin{align}
& (\widetilde{\partial}(R_h^u u)'^{k} \pm \widetilde{\partial}u'^{k}  - u''(t_{k-\frac{1}{2}}),\phi_h) 
= (\widetilde{\partial}\omega_u'^k,\phi_h) 
+  (\widetilde{\partial}u'^{k} - u''(t_{k-\frac{1}{2}}),\phi_h);
\label{taylor1}
\\
&
  \mu^{k-\frac{1}{2}}(\nabla(\widetilde{R_h^u u}^{\,k} \pm \widetilde{u}^k-  u(t_{k-\frac{1}{2}})),\nabla\phi_h)
= \mu^{k-\frac{1}{2}}(\nabla\widetilde{\omega}_u^{\,k},\nabla\phi_h) + \mu^{k-\frac{1}{2}}(\nabla(\widetilde{u}^k -  u(t_{k-\frac{1}{2}})),\nabla\phi_h);
\label{taylor2}
\\
&
\big((\mathbf{a}\cdot\nabla)(\widetilde{R_h^\theta\theta}^{\,k} \pm \widetilde{\theta}^k -\theta(t_{k-\frac{1}{2}})),\phi_h\big)
= ((\mathbf{a}\cdot\nabla)\widetilde{\omega}_\theta^{\,k},\phi_h) 
+ \big((\mathbf{a}\cdot\nabla)(\widetilde{\theta}^k -\theta(t_{k-\frac{1}{2}})),\phi_h\big);
\\
&
\mu^{k-\frac{1}{2}}\big(\widetilde{R_h^\delta\delta}'^k \pm \widetilde{\delta}'^k - \delta'(t_{k-\frac{1}{2}}),\phi_h\big)_{\Gamma_1}=
\mu^{k-\frac{1}{2}}(\widetilde{\omega}_\delta'^{\,k},\phi_h)_{\Gamma_1}
+ \mu^{k-\frac{1}{2}}\big(\widetilde{\delta}'^k - \delta'(t_{k-\frac{1}{2}}),\phi_h\big)_{\Gamma_1}.
\label{taylor4}
\end{align}
Thus, using \eqref{taylor1}-\eqref{taylor4} in \eqref{eq:proj_u1}, it follows that
\begin{align}
& (\widetilde{\partial}(R_h^u u)'^k,\phi_h) 
+ \mu^{k-\frac{1}{2}}(\nabla\widetilde{R_h^u u}^{\,k},\nabla\phi_h) 
+ \mu^{k-\frac{1}{2}}\big(g(\widetilde{R_h^u u}'^k),\phi_h\big)_{\Gamma_1}
+ ((\mathbf{a}\cdot\nabla)\widetilde{R_h^\theta\theta}^{\,k},\phi_h) \nonumber
\\
&
\quad
+ \mu^{k-\frac{1}{2}}(\widetilde{R_h^\delta\delta}'^k,\phi_h)_{\Gamma_1}
= - (f(u(t_{k-\frac{1}{2}})),\phi_h)   
+ (\widetilde{\partial}\omega_u'^k,\phi_h)
+ (\widetilde{\partial}u'^{k} - u''(t_{k-\frac{1}{2}}),\phi_h)
\nonumber
\\
&\quad 
+ \mu^{k-\frac{1}{2}}(\nabla\widetilde{\omega}_u^{\,k},\nabla\phi_h)
+ \mu^{k-\frac{1}{2}}(\nabla(\widetilde{u}^k -  u(t_{k-\frac{1}{2}})),\nabla\phi_h)
+ \mu^{k-\frac{1}{2}}\big(g(\widetilde{R_h^u u}'^k) - g(u'(t_{k-\frac{1}{2}})),\phi_hz\big)_{\Gamma_1} 
\nonumber
\\
&\quad
+ ((\mathbf{a}\cdot\nabla)\widetilde{\omega}_\theta^{\,k},\phi_h)
+ \big((\mathbf{a}\cdot\nabla)(\widetilde{\theta}^k-\theta(t_{k-\frac{1}{2}})),\phi_h\big) 
+ \mu^{k-\frac{1}{2}}(\widetilde{\omega}_\delta'^{\,k},\phi_h)_{\Gamma_1}
+ \mu^{k-\frac{1}{2}}\big(\widetilde{\delta}'^k - \delta'(t_{k-\frac{1}{2}}),\phi_h\big)_{\Gamma_1}.
\label{eq:proj_u1_rewritten}
\end{align}
  
\subsubsection{Projection equation for $\theta$}

Analogously to what was done for $u$ in \eqref{eq:proj_u1}, now for the equation of $\theta$. For all $\psi_h\in W_h$,
\begin{align}
& (\widetilde{\partial}(R_h^\theta\theta)^k,\psi_h) 
+ \beta\big(\hspace{-0.09cm}\int_\Omega\widetilde{\Theta}_h^{\,k}dx\big)
(\nabla \widetilde{R_h^\theta\theta}^k,\nabla\psi_h)
+ ((\mathbf{a}\cdot\nabla)\widetilde{R_h^u u}'^k,\psi_h)
=
 \beta\big(\hspace{-0.09cm}\int_\Omega\widetilde{\Theta}_h^{\,k}dx\big)
(\nabla \widetilde{\theta}^{k},\nabla\psi_h)  \nonumber\\
&\quad 
- \beta\big(\hspace{-0.09cm}\int_\Omega\theta(t_{k-\frac{1}{2}})dx\big)
(\nabla \theta(t_{k-\frac{1}{2}}),\nabla\psi_h)
 + (\widetilde{\partial}(R_h^\theta\theta)^k - \theta'(t_{k-\frac{1}{2}}),\psi_h) 
+ \big((\mathbf{a}\cdot\nabla)\big(\widetilde{R_h^u u}'^k - u'(t_{k-\frac{1}{2}})\big),\psi_h\big),
\label{eq:proj_theta1}
\end{align}
where we used the fact that $(\nabla \widetilde{R_h^\theta\theta}^k,\nabla\psi_h)=(\nabla \widetilde{\theta}^{k},\nabla\psi_h)$ from the Ritz projection, see Definition \ref{def:ritz}.

We can rewrite \eqref{eq:proj_theta1} and highlight the truncation error part via Taylor expansion; to this end:
\begin{align}
& 
(\widetilde{\partial}(R_h^\theta\theta)^k \pm \widetilde{\partial}\theta^{k} - \theta'(t_{k-\frac{1}{2}}),\psi_h)  
= (\widetilde{\partial}\omega_\theta^k,\psi_h) 
+ (\widetilde{\partial}\theta^{k} - \theta'(t_{k-\frac{1}{2}}),\psi_h); \label{theta1-} 
\\
&
\big((\mathbf{a}\cdot\nabla)\big(\widetilde{R_h^u u}'^k \pm \widetilde{u}'^k - u'(t_{k-\frac{1}{2}})\big),\psi_h\big) 
= 
((\mathbf{a}\cdot\nabla)\widetilde{\omega}_u'^{\,k},\psi_h)
+  \big((\mathbf{a}\cdot\nabla)\big(\widetilde{u}'^k - u'(t_{k-\frac{1}{2}})\big),\psi_h\big). \label{theta2-}
\end{align}

Thus, returning to \eqref{eq:proj_theta1} and using \eqref{theta1-}-\eqref{theta2-}, it follows that:
\begin{align}
& (\widetilde{\partial}(R_h^\theta\theta)^k,\psi_h) 
+ \beta\big(\hspace{-0.09cm}\int_\Omega\widetilde{\Theta}_h^{\,k}dx\big)
(\nabla \widetilde{R_h^\theta\theta}^k,\nabla\psi_h)
+ ((\mathbf{a}\cdot\nabla)\widetilde{R_h^u u}'^k,\psi_h)
= \beta\big(\hspace{-0.09cm}\int_\Omega\widetilde{\Theta}_h^{\,k}dx\big)
(\nabla \widetilde{\theta}^{k},\nabla\psi_h) 
+ (\widetilde{\partial}\omega_\theta^k,\psi_h)
 \nonumber\\
&\;
\hspace{-0.09cm}-\hspace{-0.09cm} \beta\big(\hspace{-0.09cm}\int_\Omega\theta(t_{k-\frac{1}{2}})dx\big)
(\nabla \theta(t_{k-\frac{1}{2}}),\nabla\psi_h)
\hspace{-0.07cm}+\hspace{-0.07cm} (\widetilde{\partial}\theta^{k} \hspace{-0.09cm}-\hspace{-0.09cm} \theta'(t_{k-\frac{1}{2}}),\psi_h)
\hspace{-0.07cm}+\hspace{-0.07cm} ((\mathbf{a}\cdot\nabla)\widetilde{\omega}_u'^{\,k},\psi_h)
\hspace{-0.07cm}+\hspace{-0.07cm} \big((\mathbf{a}\cdot\nabla)\big(\widetilde{u}'^k \hspace{-0.09cm}-\hspace{-0.09cm} u'(t_{k-\frac{1}{2}})\big),\psi_h\big).
\label{eq:proj_theta1_rewritten}
\end{align}

\subsubsection{Projection equation for $\delta$}

For the equation of $\delta$, following analogously what was done in \eqref{eq:proj_u1}. For all $\chi_h\in Z_h$,
\begin{align}
& \big(\widetilde{R_h^u u}'^k, \chi_h\big)_{\Gamma_1}
 + \big(f_1\widetilde{\partial}(R_h^\delta\delta)'^k, \chi_h\big)_{\Gamma_1} 
 + \big(f_2\widetilde{R_h^\delta\delta}'^k , \chi_h\big)_{\Gamma_1}
 + \big(f_3\widetilde{R_h^\delta\delta}^k,\chi_h\big)_{\Gamma_1} 
= 
\big(\widetilde{R_h^u u}'^k - u'(t_{k-\frac{1}{2}}\big), \chi_h)_{\Gamma_1}
\nonumber\\
&\quad 
+ \big(f_1\big(\widetilde{\partial}(R_h^\delta\delta)'^k - \delta''(t_{k-\frac{1}{2}})\big) , \chi_h\big)_{\Gamma_1} 
+ \big(f_2\big(\widetilde{R_h^\delta\delta}'^k - \delta'(t_{k-\frac{1}{2}})\big) , \chi_h\big)_{\Gamma_1}
+ \big(f_3\big(\widetilde{R_h^\delta\delta}^{\,k} - \delta(t_{k-\frac{1}{2}})\big),\chi_h\big)_{\Gamma_1}.
\label{eq:proj_delta}
\end{align}
Analogously to $u$ and $\theta$, we can rewrite \eqref{eq:proj_delta} and highlight the truncation error part via Taylor expansion; to this end:
\begin{align}
& 
\big(\widetilde{R_h^u u}'^k \pm \widetilde{u}'^k- u'(t_{k-\frac{1}{2}}), \chi_h\big)_{\Gamma_1}
= 
(\widetilde{\omega}_u'^{\,k}, \chi_h)_{\Gamma_1}
+ (\widetilde{u}'^k- u'(t_{k-\frac{1}{2}}), \chi_h)_{\Gamma_1};
\label{delta1-}
\\
&
\big(f_1\big(\widetilde{\partial}(R_h^\delta\delta)'^k \pm \widetilde{\partial}\delta'^{k} - \delta''(t_{k-\frac{1}{2}})\big) , \chi_h\big)_{\Gamma_1} 
=
\big(f_1\widetilde{\partial}\omega_\delta'^k, \chi_h\big)_{\Gamma_1}  
+ \big(f_1(\widetilde{\partial}\delta'^{k} - \delta''(t_{k-\frac{1}{2}})) , \chi_h\big)_{\Gamma_1};
\\
& 
\big(f_2\big(\widetilde{R_h^\delta\delta}'^k \pm \widetilde{\delta}'^k - \delta'(t_{k-\frac{1}{2}})\big) , \chi_h\big)_{\Gamma_1}
=
\big(f_2\widetilde{\omega}_\delta'^{\,k}, \chi_h\big)_{\Gamma_1}  
+ \big(f_2(\widetilde{\delta}'^k - \delta'(t_{k-\frac{1}{2}})) , \chi_h\big)_{\Gamma_1};
\\ 
& 
\big(f_3\big(\widetilde{R_h^\delta\delta}^{\,k} \pm \widetilde{\delta}^k - \delta(t_{k-\frac{1}{2}})\big),\chi_h\big)_{\Gamma_1}
= 
\big(f_3\widetilde{\omega}_\delta^{\,k}, \chi_h\big)_{\Gamma_1}
+ 
\big(f_3(\widetilde{\delta}^k - \delta(t_{k-\frac{1}{2}})),\chi_h\big)_{\Gamma_1}.
\label{delta4-}
\end{align}
Thus, returning to \eqref{eq:proj_delta} and using \eqref{delta1-}-\eqref{delta4-}, it follows that:
\begin{align}
& \big(\widetilde{R_h^u u}'^k, \chi_h\big)_{\Gamma_1}
 + \big(f_1\widetilde{\partial}(R_h^\delta\delta)'^k, \chi_h\big)_{\Gamma_1} 
 + \big(f_2\widetilde{R_h^\delta\delta}'^k , \chi_h\big)_{\Gamma_1}
 + \big(f_3\widetilde{R_h^\delta\delta}^k,\chi_h\big)_{\Gamma_1} 
= 
(\widetilde{\omega}_u'^{\,k}, \chi_h)_{\Gamma_1}
\nonumber\\
&\quad 
+ (\widetilde{u}'^k- u'(t_{k-\frac{1}{2}}), \chi_h)_{\Gamma_1}
+ \big(f_1\widetilde{\partial}\omega_\delta'^k, \chi_h\big)_{\Gamma_1}  
+ \big(f_1(\widetilde{\partial}\delta'^{k} - \delta''(t_{k-\frac{1}{2}})) , \chi_h\big)_{\Gamma_1}
+ \big(f_2\widetilde{\omega}_\delta'^{\,k}, \chi_h\big)_{\Gamma_1}
\nonumber\\
&\quad
+ \big(f_2(\widetilde{\delta}'^k - \delta'(t_{k-\frac{1}{2}})) , \chi_h\big)_{\Gamma_1} 
+ \big(f_3\widetilde{\omega}_\delta^{\,k}, \chi_h\big)_{\Gamma_1}
+ \big(f_3(\widetilde{\delta}^k - \delta(t_{k-\frac{1}{2}})),\chi_h\big)_{\Gamma_1}.\label{eq:proj_delta_rewritten}
\end{align}
\begin{remark}
\label{remark}
For the error analysis, in addition to the hypotheses of Theorem \ref{T2.1}, we assume the Lipschitz conditions for $\beta$, the boundedness of $\mu$ and $f_i$ (see \eqref{hipoteses}), as well as the regularity conditions from Theorem \ref{teo:convergencia} given in \eqref{u-regularidade}-\eqref{delta-regularidade}. These are necessary for the Ritz projection estimates (Lemma \ref{lem:proj}), the bootstrap argument, and the control of the truncation terms.
\begin{align}
|\beta(r_1)\hspace{-0.07cm}-\hspace{-0.07cm}\beta(r_2)|\hspace{-0.07cm}\leqslant \hspace{-0.07cm}L_1 |r_1\hspace{-0.07cm}-\hspace{-0.07cm}r_2|, \; \forall \hspace{-0.02cm}r_1,r_2 \in \mathbb{R};
\; 0<\mu_0 \leqslant \mu(t) \leqslant \mu_1, \forall \hspace{-0.02cm}t\hspace{-0.02cm}\in[0,T];
\; 0< \hspace{-0.07cm}\underline{f}_i\leqslant f_i(x) \leqslant \hspace{-0.07cm}\overline{f}_i,  \forall x\in\Gamma_1;
\label{hipoteses}
\end{align}
with $\mu\in C^1([0,T])$.
\end{remark}
\subsection{Equations for the discretization errors $\rho$}
Subtracting the fully discrete numerical scheme \eqref{PTD} from the projection equations \eqref{eq:proj_u1_rewritten}, \eqref{eq:proj_theta1_rewritten}, and \eqref{eq:proj_delta_rewritten}, using the error decompositions defined in \eqref{decomp-1}–\eqref{decomp-2}, and rearranging the terms on the left-hand side in a manner convenient for estimation, we obtain:
\subsubsection{Equation for $\rho_u$}
For all $\phi_h\in V_h$, and using $\mu^{k-\frac{1}{2}}(\nabla\widetilde{\omega}_u^{\,k},\nabla\phi_h)=0$ (Ritz projection), see Definition \ref{def:ritz}, the result is
\begin{align}
& (\widetilde{\partial}\rho_u'^k,\phi_h) 
+ \mu^{k-\frac{1}{2}}(\nabla\widetilde{\rho}_u^{\,k},\nabla\phi_h) 
+ \mu^{k-\frac{1}{2}}\big(g(\widetilde{P}_h^{\,k}) - (g(\widetilde{R_h^u u}'^k),\phi_h\big)_{\Gamma_1} 
= - \big(f(\widetilde{U}_h^{\,k}) - f(u(t_{k-\frac{1}{2}})),\phi_h\big)
\nonumber
\\
&\quad 
- ((\mathbf{a}\cdot\nabla)\widetilde{\rho}_\theta^{\,k},\phi_h) 
+ \mu^{k-\frac{1}{2}}(\widetilde{\rho}_\delta{'}^{k},\phi_h)_{\Gamma_1} 
- (\widetilde{\partial}\omega_u'^k,\phi_h) - ((\mathbf{a}\cdot\nabla)\widetilde{\omega}_\theta^{\,k},\phi_h) 
+ \mu^{k-\frac{1}{2}}(\widetilde{\omega}_\delta'^{\,k},\phi_h)_{\Gamma_1} 
+ \mathcal{R}_u^k(\phi_h), \label{eq:erro_u}
\end{align}
where 
\begin{align}
 \mathcal{R}_u^k(\phi_h)&= -(\widetilde{\partial}u'^{k} - u''(t_{k-\frac{1}{2}}),\phi_h)
- \mu^{k-\frac{1}{2}}(\nabla(\widetilde{u}^k -  u(t_{k-\frac{1}{2}})),\nabla\phi_h)
\nonumber\\
&\quad 
-\big((\mathbf{a}\cdot\nabla)(\widetilde{\theta}^k -\theta(t_{k-\frac{1}{2}})),\phi_h\big)
-\mu^{k-\frac{1}{2}}\big(\widetilde{\delta}'^k - \delta'(t_{k-\frac{1}{2}}),\phi_h\big)_{\Gamma_1}.\label{resto_u}
\end{align}
\subsubsection{Equation for $\rho_\theta$}
For all $\psi_h\in W_h$, 
\begin{align}
& (\widetilde{\partial}\rho_\theta^k,\psi_h)
+\beta\big(\hspace{-0.09cm}\int_\Omega\widetilde{\Theta}_h^{\,k}dx\big)
(\nabla \widetilde{\rho}_\theta^{\,k},\nabla\psi_h)  
= -((\mathbf{a}\cdot\nabla)\widetilde{\rho}_u'^{\,k},\psi_h) - (\widetilde{\partial}\omega_\theta^k,\psi_h) - ((\mathbf{a}\cdot\nabla)\widetilde{\omega}_u'^{\,k},\psi_h) 
\nonumber
\\
&\qquad 
	- \Big[\beta\big(\hspace{-0.09cm}\int_\Omega \widetilde{\Theta}_h^{\,k}dx\big)(\nabla \widetilde{\theta}^{k},\nabla\psi_h)
	- \beta\big(\hspace{-0.09cm}\int_\Omega\theta(t_{k-\frac{1}{2}})dx\big)(\nabla\theta(t_{k-\frac{1}{2}}),\nabla\psi_h)\Big] 
	+ \mathcal{R}_\theta^k(\psi_h), \label{eq:erro_theta}
\end{align}
where 
\begin{align}
 \mathcal{R}_\theta^k(\psi_h)=-(\widetilde{\partial}\theta^{k} - \theta'(t_{k-\frac{1}{2}}),\psi_h)
- \big((\mathbf{a}\cdot\nabla)(\widetilde{u}'^k - u'(t_{k-\frac{1}{2}})),\psi_h\big).\label{resto_theta}
\end{align}
\subsubsection{Equation for $\rho_\delta$}
For all $\chi_h\in Z_h$,
\begin{align}
& (f_1\widetilde{\partial}\rho_\delta'^k,\chi_h)_{\Gamma_1} + (f_2\widetilde{\rho}_\delta'^{\,k},\chi_h)_{\Gamma_1} + (f_3\widetilde{\rho}_\delta^{\,k},\chi_h)_{\Gamma_1} = -(\widetilde{\rho}_u'^{\,k},\chi_h)_{\Gamma_1} 
\nonumber
\\
&\quad -(\widetilde{\omega}_u'^{\,k},\chi_h)_{\Gamma_1} - (f_1\widetilde{\partial}\omega_\delta'^k,\chi_h)_{\Gamma_1} - (f_2\widetilde{\omega}_\delta'^{\,k},\chi_h)_{\Gamma_1} - (f_3\widetilde{\omega}_\delta^{\,k},\chi_h)_{\Gamma_1} + \mathcal{R}_\delta^k(\chi_h), \label{eq:erro_delta}
\end{align}
where                
\begin{align}
\mathcal{R}_\delta^k(\chi_h)&=
- (\widetilde{u}'^k- u'(t_{k-\frac{1}{2}}), \chi_h)_{\Gamma_1}
- \big(f_1(\widetilde{\partial}\delta'^{k} - \delta''(t_{k-\frac{1}{2}})) , \chi_h\big)_{\Gamma_1}\nonumber
\\
&\quad 
- \big(f_2(\widetilde{\delta}'^k - \delta'(t_{k-\frac{1}{2}})) , \chi_h\big)_{\Gamma_1}
- \big(f_3(\widetilde{\delta}^k - \delta(t_{k-\frac{1}{2}})),\chi_h\big)_{\Gamma_1}.\label{resto_delta}
\end{align}

Observing that, regarding the terms $\omega^k_u$, $\omega^k_\theta$ and $\omega^k_\delta$, we already have the error estimates from Lemma \ref{lem:proj}, what remains to be estimated are the terms involving $\rho^k_u$, $\rho^k_\theta$ and $\rho^k_\delta$.

\subsubsection{Estimate for $\rho_u$}
We choose $\phi_h = \widetilde{\rho}_u'^{\,k}$ in \eqref{eq:erro_u} and \eqref{resto_u}. We will work on each term.

%
\noindent\textbf{Term}: $(\widetilde{\partial}\rho_u'^k,\widetilde{\rho}_u'^{k}).$
 Using the definition of $\widetilde{\partial}\rho_u'^k$ and $\widetilde{\rho}_u'^{k}$, it follows
\begin{align}
(\widetilde{\partial}\rho_u'^k,\widetilde{\rho}_u'^{k}) = 
\Big(\frac{\rho_{u}'^{k}-\rho_{u}'^{k-1}}{\tau},  \frac{\rho_{u}'^{k}+\rho_{u}'^{k-1}}{2}\Big)
  =\frac{1}{2\tau}\big(\|\rho_u'^k\|^2 - \|\rho_u'^{k-1}\|^2\big)
  =\frac{1}{2}\widetilde{\partial}\|\rho_u'^k\|^2.\label{termo_u_1}
\end{align}

\noindent\textbf{Term}: $\mu^{k-\frac{1}{2}} \big(\nabla \widetilde{\rho}_{u}^k , \nabla \widetilde{\rho}_u'^{k}\big)$. Note that, using the linearity of the operators $\widetilde{\partial}$ and $R_h^u$, it follows by definition,
\begin{align}
\widetilde{\rho}_u'^{k}=\widetilde{P}_h^k - \widetilde{R_h^u u}'^k
=\widetilde{\partial}U_h^k - \widetilde{R_h^u u}'^k 
&= \big(\widetilde{\partial}U_h^k - \widetilde{\partial}(R_h^u u)^k\big) + \big(\widetilde{\partial}(R_h^u u)^k - \widetilde{R_h^u u}'^k\big)
= \widetilde{\partial}\rho_{u}^k + R_h^u(\widetilde{\partial}u^k -\widetilde{u}'^k).\label{def1}
\end{align}
Thus, using \eqref{def1},
\begin{align}
\mu^{k-\frac{1}{2}} \big(\nabla \widetilde{\rho}_{u}^k , \nabla \widetilde{\rho}_u'^{k}\big)&=\mu^{k-\frac{1}{2}}\big(\nabla \widetilde{\rho}_{u}^k, \nabla\widetilde{\partial}\rho_{u}^k\big) + \mu^{k-\frac{1}{2}}\big(\nabla \widetilde{\rho}_{u}^k, \nabla R_h^u(\widetilde{\partial}u^k -\widetilde{u}'^k) \big)\
\nonumber
\\ &=\frac{\mu^{k-\frac{1}{2}}}{2}\widetilde{\partial}\|\nabla \rho_u^k\|^2
+ \mu^{k-\frac{1}{2}}\big(\nabla \widetilde{\rho}_{u}^k, \nabla R_h^u(\widetilde{\partial}u^k -\widetilde{u}'^k) \big). \label{pho1}
\end{align}
Now, the second term of the equality above is moved to the right-hand side of \eqref{eq:erro_u}, and using Cauchy-Schwarz and $\epsilon$-Young, along with the fact that $\|\nabla R_h^u \varphi\|\leqslant \|\nabla \varphi\|$ for all $\varphi \in V$ (see \eqref{eq:estab_u}), it follows that
\begin{align}
\hspace{-0.07cm}-\mu^{k-\frac{1}{2}}\big(\nabla \widetilde{\rho}_{u}^k, \nabla R_h^u(\widetilde{\partial}u^k \hspace{-0.07cm}-\hspace{-0.07cm}\widetilde{u}'^k) \big) 
&\hspace{-0.07cm}\leqslant\hspace{-0.07cm}
\epsilon_0||\nabla \widetilde{\rho}_{u}^k||^2 
     \hspace{-0.07cm}+\hspace{-0.07cm} \frac{\mu_1^{2}}{4\epsilon_0}||\nabla R_h^u(\widetilde{\partial}u^k -\widetilde{u}'^k)||^2  
     \hspace{-0.07cm}\leqslant\hspace{-0.07cm}
   \epsilon_0||\nabla \widetilde{\rho}_{u}^k||^2 
   \hspace{-0.07cm} +\hspace{-0.07cm} \frac{\mu_1^{2}}{4\epsilon_0}||\nabla(\widetilde{\partial}u^k \hspace{-0.07cm}-\hspace{-0.07cm} \widetilde{u}'^k)||^2. \label{pho2}
\end{align}

\noindent\textbf{Term}: $\mu^{k-\frac{1}{2}}\big(g(\widetilde{P}_h^{\,k})-g(\widetilde{R_h^u u}'^k),\widetilde{\rho}_u'^{\,k}\big)_{\Gamma_1}.$ 
Using the monotonicity property of the nonlinear function $g$ in \eqref{hip:eta_mon} and the hypothesis regarding $\mu$ in \eqref{hipoteses},
and recalling that $\widetilde{\rho}_u'^{\,k}=\widetilde{P}_h^{\,k} - \widetilde{R_h^u u}'^k$, it follows that
\begin{align}
\mu^{k-\frac{1}{2}}\big(g(\widetilde{P}_h^{\,k})-g(\widetilde{R_h^u u}'^k),\widetilde{\rho}_u'^{\,k}\big)_{\Gamma_1}
\geqslant \mu_0 g_0\|\widetilde{\rho}_u'^{\,k}\|_{\Gamma_1}^2
\geqslant \frac{\mu_0 g_0}{2}\|\widetilde{\rho}_u'^{\,k}\|_{\Gamma_1}^2.\label{pho3}
\end{align}

The terms $- ((\mathbf{a}\cdot\nabla)\widetilde{\rho}_\theta^{\,k}, \widetilde{\rho}_u'^{\,k})$ and $\mu^{k-\frac{1}{2}}(\widetilde{\rho}_\delta'^{\,k},\widetilde{\rho}_u'^{\,k})_{\Gamma_1}$ will not be estimated, as they will subsequently be cancelled: the first with a term in Section ref{sec:tho-theta} and the second with the term in Section \ref{sec:rho-delta}.

\noindent\textbf{Term}: $(\widetilde{\partial}\omega_u'^k,\widetilde{\rho}_u'^{\,k}).$ By the Cauchy–Schwarz and $\epsilon$-Young inequalities,
\begin{align}
|(\widetilde{\partial}\omega_u'^k,\widetilde{\rho}_u'^{\,k})| 
&\leqslant \frac{1}{4\epsilon_1}\|\widetilde{\partial}\omega_u'^k\|^2 
+ \epsilon_1\|\widetilde{\rho}_u'^{\,k}\|^2.\label{rho4}
\end{align}

\noindent\textbf{Term}: $-((\mathbf{a}\cdot\nabla)\widetilde{\omega}_\theta^{\,k},\widetilde{\rho}_u'^{\,k})$. As $\widetilde{\omega}_\theta^{\,k}\in H^1_0(\Omega)$, by Green's Theorem, 
$-((\mathbf{a}\cdot\nabla)\widetilde{\omega}_\theta^{\,k},\widetilde{\rho}_u'^{\,k}) = (\widetilde{\omega}_\theta^{\,k},\, (\mathbf{a}\cdot\nabla)\widetilde{\rho}_u'^{\,k}).$
Note that we removed the gradient of $\widetilde{\omega}^k_\theta$, but it cannot accompany $\widetilde{\rho}_{u'}^k$, since there is no $\|\nabla\widetilde{\rho}_{u'}^k\|^2$ on the left-hand side of the inequality in \eqref{resumo_u}.
We will circumvent this situation by resorting to the relation between $\rho_{u}'^{k}$ and $\rho_u^k$  established in \eqref{decomp-2}, namely,
$
  \widetilde{\rho}{_{u}'^{k}}
= \widetilde{\partial}\rho_u^k + R^u_h(\widetilde{\partial}u^k - \widetilde{u}'^k).
$
Thus,
\begin{align}
  \big((\mathbf{a}\cdot \nabla)\widetilde{\rho}{_{u}'^{k}},\widetilde{\omega}^k_\theta\big)
=&
  \big((\mathbf{a}\cdot \nabla)\widetilde{\partial}\rho_u^k,\widetilde{\omega}^k_\theta\big)
+ \big((\mathbf{a}\cdot \nabla)R^u_h(\widetilde{\partial}u^k - \widetilde{u'}^k),\widetilde{\omega}^k_\theta\big).\label{indent0}
\end{align}
However, we have a new problem in
$\big((\mathbf{a}\cdot \nabla)\widetilde{\partial}\rho_u^k,\widetilde{\omega}^k_\theta\big)$; we need to remove the operator $\widetilde{\partial}$ from $\rho_u^k$.
To this end, we will use the identity
\begin{align}
\big((\mathbf{a}\cdot \nabla)\widetilde{\partial}\rho_u^k,\widetilde{\omega}^k_\theta\big)
=&
  \widetilde{\partial}\big((\mathbf{a}\cdot \nabla)\rho_u^k,\widetilde{\omega}^k_\theta\big)
- \big((\mathbf{a}\cdot \nabla)\rho_u^{k-1},\widetilde{\partial}\widetilde{\omega}^k_\theta\big),\label{indent1}
\end{align}
where, 
$
  \widetilde{\partial}\big((\mathbf{a}\cdot \nabla)\rho_u^k,\widetilde{\omega}^k_\theta\big)
=
  \big[\big((\mathbf{a}\cdot \nabla)\rho_u^k,\widetilde{\omega}^k_\theta\big)
- \big((\mathbf{a}\cdot \nabla)\rho_u^{k-1},\widetilde{\omega}^{k-1}_\theta\big)\big]/ \tau
\;\mbox{ and }\;
  \widetilde{\partial}\widetilde{\omega}^k_\theta=\big[\widetilde{\omega}^k_\theta-\widetilde{\omega}^{k-1}_\theta\big]/ \tau.
$
Therefore, from \eqref{indent1} into \eqref{indent0}, using the Cauchy-Schwarz and Young inequalities and the stability of the projection in \eqref{eq:estab_u},
\begin{align}\label{cap3:I_1^n:parcela1}
\big((\mathbf{a}\cdot \nabla)\widetilde{\rho}_{u'}^k,\widetilde{\omega}^k_\theta\big)
&\leqslant 
  \widetilde{\partial}\big((\mathbf{a}\cdot \nabla)\rho_u^k,\widetilde{\omega}^k_\theta\big)
  + \|\nabla\rho_u^{k-1}\|^2
  + \frac{\|\mathbf{a}\|_{\infty}^2}{2}\|\widetilde{\partial}\widetilde{\omega}^k_\theta\|^2 
  + \frac{\|\mathbf{a}\|_{\infty}^2}{4}\|\nabla(\widetilde{\partial}u^k - \widetilde{u}'^k)\|^2 
  + \|\widetilde{\omega}^k_\theta\|^2.
\end{align}
Since $\widetilde{\partial}\big((\mathbf{a}\cdot \nabla)\rho_u^k,\widetilde{\omega}^k_\theta\big)$ is a telescoping term, after summing the estimates, the two remaining terms will be handled via Cauchy-Schwarz and Young inequalities.

\noindent\textbf{Term}: $\mu^{k-\frac{1}{2}}(\widetilde{\omega}_\delta'^{\,k},\widetilde{\rho}_u'^{\,k})_{\Gamma_1}$. By Cauchy--Schwarz and $\epsilon$-Young,
\begin{align}
|\mu^{k-\frac{1}{2}}(\widetilde{\omega}_\delta'^{\,k},\widetilde{\rho}_u'^{\,k})_{\Gamma_1}| 
&\leqslant \frac{\mu_1^{2}}{4\epsilon_2}\|\widetilde{\omega}_\delta'^{\,k}\|_{\Gamma_1}^2 
+ \epsilon_2\|\widetilde{\rho}_u'^{\,k}\|_{\Gamma_1}^2.\label{rho6}
\end{align}

\noindent\textbf{Nonlinear term with $g$}: 
Using the hypothesis that $g$ is Lipschitz (see \eqref{hip:eta_lip}), and \( 0 < \mu_0 \leqslant \mu(t) \leqslant \mu_1 \) (see \eqref{hipoteses}), together with the Cauchy-Schwarz and $\epsilon$-Young inequalities,
\begin{align}
&  \mu^{k-\frac{1}{2}}\big(g(\widetilde{R_h^u u}'^k) \hspace{-0.05cm}-\hspace{-0.05cm}
g(u'(t_{k-\frac{1}{2}})), \widetilde{\rho}_u'^{\,k}\big)_{\Gamma_1} 
\hspace{-0.05cm}\leqslant \hspace{-0.05cm} 
  L_g\mu_1 \|\widetilde{R_h^u u}'^k \hspace{-0.07cm}\pm \widetilde{u}'^k  \hspace{-0.05cm}-\hspace{-0.05cm} u'(t_{k-\frac{1}{2}})\|_{\Gamma_1}\|\widetilde{\rho}_u'^{\,k}\|_{\Gamma_1} 
\nonumber 
\\
& \quad
\leqslant     
\hspace{-0.07cm}L_g\mu_1 \Big(\|\widetilde{\omega}_u'^k\|_{\Gamma_1} 
\hspace{-0.07cm}+\hspace{-0.07cm} \|\widetilde{u}'^k \hspace{-0.07cm}-\hspace{-0.07cm} u'(t_{k-\frac{1}{2}})\|_{\Gamma_1}\Big)\|\widetilde{\rho}_u'^{\,k}\|_{\Gamma_1}
\hspace{-0.07cm}\leqslant\hspace{-0.07cm}
\frac{L^2_g\mu_1^{2}}{2\epsilon_2}
\|\widetilde{\omega}_u'^k\|^2_{\Gamma_1} 
\hspace{-0.07cm}+\hspace{-0.07cm} \frac{L^2_g\mu_1^{2}}{2\epsilon_2}\|\widetilde{u}'^k \hspace{-0.07cm}-\hspace{-0.07cm} u'(t_{k-\frac{1}{2}})\|^2_{\Gamma_1}
\hspace{-0.07cm}+\hspace{-0.07cm} \epsilon_2\|\widetilde{\rho\,}_u'^{\,k}\|^2_{\Gamma_1}.
\end{align}

\noindent\textbf{Nonlinear term with $f$}: $\mathcal{N}_f =- \big(f(\widetilde{U}^k_h) - f(u(t_{k-\frac{1}{2}})),\widetilde{\rho}_u'^{\,k}\big).$
To begin the induction argument, we assume that there exists a constant $M>0$ (to be determined) such that, for all indices $i < k$ (with $k\geqslant 1$),
$
\|\widetilde{U}^i_h\|_{L^\infty(\Omega)} \leqslant M.
$

\noindent\textit{Induction base case ($i=0$).}
Since, by \eqref{u-regularidade}, $u_0\in H^{r+1}(\Omega)\cap V$ with $r+1 > n/2$ (a condition satisfied for $n=2,3$ and $r\geqslant 1$), the continuous Sobolev embedding $H^{r+1}(\Omega)\hookrightarrow L^\infty(\Omega)$ guarantees the existence of a constant $C_0>0$ such that $\|u_0\|_{L^\infty(\Omega)}\leqslant C_0\|u_0\|_{H^{r+1}(\Omega)}$. Furthermore, the Ritz projection $R_h^u$ satisfies the uniform error estimate (Lemma 4 of \cite{Artigo2025_Alfaro}),
$
\|R_h^u u_0 - u_0\|_{L^\infty(\Omega)}\leqslant C_1 h^{r+1-\frac{n}{2}}\|u_0\|_{H^{r+1}(\Omega)},
$
with $C_1$ independent of $h$. Since the exponent $r+1-\frac{n}{2}>0$, there exists $h_0^{(0)}>0$ such that for all $h\leqslant h_0^{(0)}$,
$
C_1 h^{r+1-\frac{n}{2}}\|u_0\|_{H^{r+1}(\Omega)}\hspace{-0.07cm} \leqslant \hspace{-0.07cm}1.
$
From this and the triangle inequality, and by choosing $h\leqslant h_0^{(0)}$, we obtain
$
\|U_h^0\|_{L^\infty(\Omega)} \leqslant \|u_0\|_{L^\infty(\Omega)}+1 =: M_0,$ since $U_h^0=R_h^u u_0.$
This establishes the base case of the induction. Note that $\widetilde{U}^0_h$ is not defined; the induction begins with $i=0$ for $U_h^0$, and for $i \geqslant 1$, we use $\widetilde{U}^i_h$.\\
\noindent\textit{Inductive step.}
For the $k$-th step ($k\geqslant 1$), we write the decomposition $ \widetilde{U}^k_h = \widetilde{R_h^u u}^{\,k} + \widetilde{\rho}_u^{\,k}, \) where \( \widetilde{R_h^u u}^{\,k} $ is the Crank--Nicolson average of the Ritz projection. To obtain an $ L^\infty $ bound for $ \widetilde{R_h^u u}^{\,k} $, we use the uniform convergence of the projection. By finite element theory (cf. \cite{Livro2002_Ciarlet}), for each $ t \in [0,T]$,
$$
\|R_h^u u(t) - u(t)\|_{L^\infty(\Omega)} \leqslant C_2 h^{r+1-\frac{n}{2}}\|u(t)\|_{H^{r+1}(\Omega)},
$$
with $C_2$ independent of $h$ and $t$. This estimate is obtained by combining the error estimate in $H^1$ (or, in our case, in $V=H^1_{\Gamma_0}(\Omega)$), the inverse inequality in $V_h$, and the Sobolev embedding $H^{r+1}(\Omega)\hookrightarrow L^\infty(\Omega)$, , which is valid since $r+1 > n/2$ ($n=2,3)$. The Sobolev embedding yields
$\operatorname{ess\,sup}_{t\in[0,T]}\|u(t)\|_{L^\infty(\Omega)} \hspace{0.07cm}\leqslant \hspace{0.07cm}C_3 \operatorname{ess\,sup}_{t\in[0,T]}\|u(t)\|_{H^{r+1}(\Omega)} < \infty,$ since $u\in L^\infty(0,T;H^{r+1}(\Omega))$, see \eqref{u-regularidade}.
We denote $\|u\|_{L^\infty([0,T]\times\Omega)} := \operatorname{ess\,sup}_{t\in[0,T]}\|u(t)\|_{L^\infty(\Omega)}$, which is a finite constant. We choose $h_0^{(1)}>0$ such that, for all $h\leqslant h_0^{(1)}$,
$
C_2 h^{r+1-\frac{n}{2}}\|u\|_{L^\infty(0,T;H^{r+1})} \leqslant 1.
$
So, for all $t\in[0,T]$ (in the sense q.s),
$
\|R_h^u u(t)\|_{L^\infty(\Omega)} \leqslant \|u(t)\|_{L^\infty(\Omega)}+1.
$
Since the Crank–Nicolson average does not amplify the $L^\infty$ norm (i.e., $\|\widetilde{R_h^u u}^{\,k}\|_{L^\infty}\leqslant \max\{\|R_h^u u(t_{k-1})\|_{L^\infty},\|R_h^u u(t_k)\|_{L^\infty}\}$), we obtain $\|\widetilde{R_h^u u}^{\,k}\|_{L^\infty(\Omega)} \leqslant \|u\|_{L^\infty([0,T]\times\Omega)}+1 =: C_* < \infty.$
The error $\widetilde{\rho}_u^{\,k}$ satisfies the  estimate (valid for all $k \geqslant 1$ under the induction hypothesis):
$
\|\widetilde{\rho}_u^{\,k}\|_{L^2} \leqslant C_4(M)\big(h^{r+\frac{1}{2}}+\tau^2\big), 
$
where $C_4(M)$ depends on $M$ but is independent of $h$ and $\tau$. Since $\widetilde{\rho}_u^{\,k}\in V_h$, the inverse inequality (cf. \cite{Livro2002_Ciarlet, Livro2008_Brenner&Scott}) yields
\begin{align}
\|\widetilde{\rho}_u^{\,k}\|_{L^\infty} \leqslant C_{\text{inv}} h^{-\frac{n}{2}}\|\widetilde{\rho}_u^{\,k}\|_{L^2}
\leqslant C_{\text{inv}}C_4(M)\big(h^{r+\frac{1}{2}-\frac{n}{2}} + h^{-n/2}\tau^2\big). \label{pho}
\end{align} 
For $n=2,3$, we have $r+\frac{1}{2}-\frac{n}{2}\geqslant 0$, since $r\geqslant 1$. The first term in \eqref{pho} tends to zero as $h \rightarrow 0$. The second term, $h^{-n/2}\tau^2$, requires a stability condition. We assume the condition $\tau = \mathcal{O}(h^{n/4})$, that is, $\tau^2/h^{n/2}\to 0$ as $h\rightarrow 0$. Therefore, for any $\epsilon>0$, there exists $h_0^{(2)}(\epsilon)>0$ such that $\tau^2 \leqslant \epsilon h^{n/2}$ for all $h\leqslant h_0^{(2)}$. We choose $\epsilon = \frac{1}{2C_{\text{inv}}C_4(M)}$. 
Furthermore, we choose $h_0^{(3)}>0$ such that $C_{\text{inv}}C_4(M) h^{r+\frac{1}{2}-\frac{n}{2}} \leqslant \frac{1}{2}$ for all $0<h\leqslant h_0^{(3)}$. Taking
\begin{align}
h \leqslant h_0=\min \{1, h_0^{(0)}, h_0^{(1)}, h_0^{(2)}, h_0^{(3)}\}\label{h}
\end{align}
and the corresponding $\tau$, we obtain
$
C_{\text{inv}}C_4(M) h^{r+\frac{1}{2}-\frac{n}{2}} \leqslant \frac{1}{2},\;
C_{\text{inv}}C_4(M) h^{-n/2}\tau^2 \leqslant \frac{1}{2}.
$
Therefore, $\|\widetilde{\rho}_u^{\,k}\|_{L^\infty} \leqslant 1$. Consequently,
$
\|\widetilde{U}^k_h\|_{L^\infty} \leqslant \|\widetilde{R_h^u u}^{\,k}\|_{L^\infty} + \|\widetilde{\rho}_u^{\,k}\|_{L^\infty} \leqslant \big(\|u\|_{L^\infty([0,T]\times\Omega)}+1\big) + 1 = \|u\|_{L^\infty([0,T]\times\Omega)}+2.
$
Now, we define $M = \max\{M_0,\; \|u\|_{L^\infty([0,T]\times\Omega)}+2\}$. The induction hypothesis extends to $k$, completing the bootstrap argument.

\noindent{\it Pointwise estimate of  $f$.}
Since $f(s)=\lambda|s|^d s$ is of class $C^1(\mathbb{R})$ with $f'(s)=\lambda(d+1)|s|^d$, and $\|\widetilde{U}^k_h\|_{L^\infty}, \|u\|_{L^\infty}\leqslant M$ (for sufficiently small $h$), the Mean Value Theorem yields
$
|f(\widetilde{U}^k_h(x))-f(u(t_{k-\frac{1}{2}},x))|
\leqslant c_f\,|\widetilde{U}^k_h(x)-u(t_{k-\frac{1}{2}},x)|,\; \forall x\in\Omega,$
where $c_f = \lambda(d+1)\big(M+\|u\|_{L^\infty([0,T]\times\Omega)}\big)^d$.
Of this, and by Cauchy–Schwarz inequality in $\mathcal{N}_f$,
\begin{align}
|\mathcal{N}_f| \leqslant c_f \|\widetilde{U}^k_h - u(t_{k-\frac{1}{2}})\|\|\widetilde{\rho}_u'^{\,k}\|.\label{Nf1}
\end{align}
But, $\widetilde{U}^k_h \hspace{-0.07cm}-\hspace{-0.07cm} u(t_{k-\frac{1}{2}}) = \widetilde{\rho}_u^{\,k} \hspace{-0.07cm}+\hspace{-0.07cm} \widetilde{\omega}_u^{\,k} \hspace{-0.07cm}+\hspace{-0.07cm} (\widetilde{u}^k \hspace{-0.07cm}-\hspace{-0.07cm} u(t_{k-\frac{1}{2}})),$ see \eqref{decomp-1}-\eqref{decomp-2}, applying the triangle inequality and the \(\epsilon\)-Young inequality to \eqref{Nf1}, followed by the Poincaré inequality, we obtain
\begin{align}
|\mathcal{N}_f| 
\hspace{-0.09cm}\leqslant 
\frac{3c_f^2c_P}{4\epsilon_1}\|\nabla\widetilde{\rho}_u^{\,k}\|^2 
+ \frac{3c_f^2}{4\epsilon_1}\|\widetilde{\omega}_u^{\,k}\|^2
+ \frac{3c_f^2}{4\epsilon_1}\|\widetilde{u}^k - u(t_{k-\frac{1}{2}})\|^2
+ \epsilon_1\|\widetilde{\rho}_u'^{\,k}\|^2,\label{Nf2}
\end{align}
where $c_P$ is the Poincaré constant.

\noindent\textbf{Term}: $\mathcal{R}_u^k(\widetilde{\rho}_u'^{\,k})$. Taking $\phi_h=\widetilde{\rho}_u'^{\,k}$ in \eqref{resto_u}, it follows
\begin{align}
 \mathcal{R}_u^k(\phi_h)&=
- (\widetilde{\partial}u'^{k} - u''(t_{k-\frac{1}{2}}),\widetilde{\rho}_u'^{\,k})
- \mu^{k-\frac{1}{2}}(\nabla(\widetilde{u}^k - u(t_{k-\frac{1}{2}})),\nabla\widetilde{\rho}_u'^{\,k})
\nonumber\\
&\quad 
-\big((\mathbf{a}\cdot\nabla)(\widetilde{\theta}^k -\theta(t_{k-\frac{1}{2}})),\widetilde{\rho}_u'^{\,k}\big)
-\mu^{k-\frac{1}{2}}\big(\widetilde{\delta}'^k - \delta'(t_{k-\frac{1}{2}}),\widetilde{\rho}_u'^{\,k}\big)_{\Gamma_1}.\label{resto_u1}
\end{align}
Using the Cauchy-Schwartz and $\epsilon$-Young inequalities,
\begin{align}
&  -(\widetilde{\partial}u'^{k} - u''(t_{k-\frac{1}{2}}),\widetilde{\rho}_u'^{\,k}) 
 \leqslant 
 \frac{1}{4\epsilon_1}\|\widetilde{\partial}u'^{k} - u''(t_{k-\frac{1}{2}})\|^2 
 +\epsilon_1\|\widetilde{\rho}_u'^{\,k}\|^2;
\\
& -\big((\mathbf{a}\cdot\nabla)(\widetilde{\theta}^k -\theta(t_{k-\frac{1}{2}})),\widetilde{\rho}_u'^{\,k}\big)
 \leqslant 
 \frac{\|a\|^2_{\infty}}{4\epsilon_1}\|\nabla(\widetilde{\theta}^k -\theta(t_{k-\frac{1}{2}}))\|^2 
 +\epsilon_1\|\widetilde{\rho}_u'^{\,k}\|^2;
\\
& -\mu^{k-\frac{1}{2}}\big(\widetilde{\delta}'^k - \delta'(t_{k-\frac{1}{2}}),\widetilde{\rho}_u'^{\,k}\big)_{\Gamma_1} 
 \leqslant 
 \frac{\mu_1^2}{4\epsilon_2}\|\widetilde{\delta}'^k - \delta'(t_{k-\frac{1}{2}})\|_{\Gamma_1}^2 
 +\epsilon_2\|\widetilde{\rho}_u'^{\,k}\|_{\Gamma_1}^2.
\end{align}
Since we do not control the term $\nabla\widetilde{\rho}_u'^{\,k}$ in the estimates, we will apply Green's theorem to the term $- \mu^{k-\frac{1}{2}}(\nabla(\widetilde{u}^k - u(t_{k-\frac{1}{2}})),\nabla\widetilde{\rho}_u'^{\,k})$ to circumvent, hence,
\begin{align*}
 - \mu^{k-\frac{1}{2}}(\nabla(\widetilde{u}^k \hspace{-0.09cm}-\hspace{-0.09cm} u(t_{k-\frac{1}{2}})),\nabla\widetilde{\rho}_u'^{\,k})
          = \mu^{k-\frac{1}{2}}\Big[\big( \Delta (\widetilde{u}^k \hspace{-0.09cm}-\hspace{-0.09cm} u(t_{k-\frac{1}{2}})), \widetilde{\rho}_u'^{\,k}\big) 
            \hspace{-0.09cm}-\hspace{-0.09cm} \big(\frac{\partial}{\partial \nu}(\widetilde{u}^k \hspace{-0.09cm}-\hspace{-0.09cm} u(t_{k-\frac{1}{2}})), \widetilde{\rho}_u'^{\,k} \big)_{\Gamma_1}\Big].
\end{align*}
By Cauchy-Schwarz and $\epsilon$-Young, and the hypothesis about $\mu$, see  \eqref{hipoteses}, it follows
\begin{align}
\hspace{-0.04cm}-\hspace{-0.07cm} \mu^{k-\frac{1}{2}}(\nabla(\widetilde{u}^k \hspace{-0.07cm}-\hspace{-0.07cm} u(t_{k-\frac{1}{2}})),\nabla\widetilde{\rho}_u'^{\,k})
\hspace{-0.07cm}
\leqslant  
\hspace{-0.07cm}\frac{\mu_1^2}{4\epsilon_1}\|\Delta (\widetilde{u}^k \hspace{-0.07cm}- \hspace{-0.07cm}u(t_{k-\frac{1}{2}}))\|^2 
		  \hspace{-0.07cm}+\hspace{-0.07cm}\frac{\mu_1^2}{4\epsilon_2}\|\frac{\partial}{\partial \nu}(\widetilde{u}^k \hspace{-0.07cm}-\hspace{-0.07cm} u(t_{k-\frac{1}{2}}))\|^2_{\Gamma_1}
    \hspace{-0.07cm}+\hspace{-0.07cm}\epsilon_1\|\widetilde{\rho}_u'^{\,k}\|^2	  	
	\hspace{-0.07cm}+\hspace{-0.07cm}\epsilon_2\|\widetilde{\rho}_u'^{\,k}\|^2_{\Gamma_1}.	\label{ultima-par}	              
\end{align}	

Returning to  \eqref{eq:erro_u}  and  \eqref{resto_u}, using \eqref{termo_u_1}-\eqref{ultima-par}, and rearranging the terms, we have
\begin{align}
& \frac{1}{2}\widetilde{\partial}\|\rho_u'^k\|^2 
+ \frac{\mu^{k-\frac{1}{2}}}{2}\widetilde{\partial}\|\nabla\rho_u^k\|^2 
+ \frac{\mu_0g_0}{2}\|\widetilde{\rho}_u'^{\,k}\|_{\Gamma_1}^2 
\leqslant  
\widetilde{\partial}\big((\mathbf{a}\cdot \nabla)\rho_u^k,\widetilde{\omega}^k_\theta\big) 
- ((\mathbf{a}\cdot\nabla)\widetilde{\rho}_\theta^{\,k}, \widetilde{\rho}_u'^{\,k})
+ \mu^{k-\frac{1}{2}}(\widetilde{\rho}_\delta'^{\,k},\widetilde{\rho}_u'^{\,k})_{\Gamma_1}\nonumber
\\
& \quad
+\Big(\epsilon_0 + \frac{3c_f^2c_P}{4\epsilon_1}\Big)\|\nabla \widetilde{\rho}_{u}^k\|^2 
+ \|\nabla\rho_u^{k-1}\|^2
+ 4\epsilon_1\|\widetilde{\rho}_u'^{\,k}\|^2
+ 4\epsilon_2\|\widetilde{\rho}_u'^{\,k}\|_{\Gamma_1}^2
+\frac{1}{4\epsilon_1}\|\widetilde{\partial}\omega_u'^k\|^2
+\frac{L^2_g\mu_1^{2}}{2\epsilon_2}\|\widetilde{\omega}_u'^k\|^2_{\Gamma_1}\nonumber
\\
& \quad 
+ \frac{3c_f^2}{4\epsilon_1}\|\widetilde{\omega}_u^{\,k}\|^2
+\frac{\|\mathbf{a}\|_{\infty}^2}{4}\|\widetilde{\partial}\widetilde{\omega}^k_\theta\|^2 
+ \|\widetilde{\omega}^k_\theta\|^2
+ \frac{\mu_1^{2}}{4\epsilon_2}\|\widetilde{\omega}_\delta'^{\,k}\|_{\Gamma_1}^2 
+ \frac{3c_f^2}{4\epsilon_1}\|\widetilde{u}^k - u(t_{k-\frac{1}{2}})\|^2
+ \frac{1}{4\epsilon_1}\|\widetilde{\partial}u'^{k} - u''(t_{k-\frac{1}{2}})\|^2\nonumber
  \\
& \quad 
+ \Big(\frac{\mu_1^{2}}{4\epsilon_0}+\frac{\|\mathbf{a}\|_{\infty}^2}{4}\Big)\|\nabla(\widetilde{\partial}u^k -\widetilde{u}'^k)\|^2
+ \frac{\|a\|^2_{\infty}}{4\epsilon_1}\|\nabla(\widetilde{\theta}^k -\theta(t_{k-\frac{1}{2}}))\|^2 
+ \frac{\mu_1^2}{4\epsilon_2}\|\widetilde{\delta}'^k - \delta'(t_{k-\frac{1}{2}})\|_{\Gamma_1}^2
\nonumber
  \\
& \quad  
+ \frac{L^2_g\mu_1^{2}}{2\epsilon_2}\|\widetilde{u}'^k - u'(t_{k-\frac{1}{2}})\|^2_{\Gamma_1} 
+\frac{\mu_1^2}{4\epsilon_1}\|\Delta (\widetilde{u}^k-u(t_{k-\frac{1}{2}}))\|^2 
+\frac{\mu_1^2}{4\epsilon_2}\|\frac{\partial}{\partial \nu}(\widetilde{u}^k-u(t_{k-\frac{1}{2}}))\|^2_{\Gamma_1}.\label{resumo_u}
\end{align}	

\subsubsection{Estimate for $\rho_\theta$}\label{sec:tho-theta}
We choose $\psi_h = \widetilde{\rho}_\theta^{\,k}$ in \eqref{eq:erro_theta} and \eqref{resto_theta}. 

\noindent\textbf{Term}: $(\widetilde{\partial}\rho_\theta^k,\widetilde{\rho}_\theta^{\,k})$. It follows analogously to the calculations in \eqref{termo_u_1}:
\begin{align}
&(\widetilde{\partial}\rho_\theta^k,\widetilde{\rho}_\theta^{\,k}) = \frac{1}{2}\widetilde{\partial}\|\rho_\theta^k\|^2.\label{theta1}
\end{align}
\noindent\textbf{Term}: $\beta\big(\int_\Omega\widetilde{\Theta}_h^{\,k}dx\big)(\nabla \widetilde{\rho}_\theta^{\,k},\nabla \widetilde{\rho}_\theta^{\,k})$. Using the hypothesis on \( \beta \) (see Theorem \ref{T2.1}),
\begin{align}
\frac{\beta_0}{2}\|\nabla \widetilde{\rho}_\theta^{\,k}\|^2
\leqslant
\beta_0\|\nabla \widetilde{\rho}_\theta^{\,k}\|^2
\leqslant
\beta\big(\hspace{-0.07cm}\int_\Omega\widetilde{\Theta}_h^{\,k}dx\big)(\nabla \widetilde{\rho}_\theta^{\,k},\nabla \widetilde{\rho}_\theta^{\,k}). \label{theta2}
\end{align} 
{\it Cancellation of coupling terms:}
The terms $-((\mathbf{a}\cdot\nabla)\widetilde{\rho}_\theta^{\,k},\widetilde{\rho}_u'^{\,k})$ and $-((\mathbf{a}\cdot\nabla)\widetilde{\rho}_u'^{\,k},\widetilde{\rho}_\theta^{\,k})$ cancel out when summing the equations, since by Green's Theorem, since $\widetilde{\rho}_\theta^{\,k}\in H^0_1(\Omega)$, we have: 
$-((\mathbf{a}\cdot\nabla)\widetilde{\rho}_u'^{\,k},\widetilde{\rho}_\theta^{\,k}) = (\widetilde{\rho}_u'^{\,k},(\mathbf{a}\cdot\nabla)\widetilde{\rho}_\theta^{\,k}).$

\noindent\textbf{Term}: $-(\widetilde{\partial}\omega_\theta^k,\widetilde{\rho}_\theta^{\,k})$. Using the Cauchy–Schwarz and $\epsilon$-Young inequalities,
\begin{align}
-(\widetilde{\partial}\omega_\theta^k,\widetilde{\rho}_\theta^{\,k})
\leqslant
|(\widetilde{\partial}\omega_\theta^k,\widetilde{\rho}_\theta^{\,k})|
\leqslant 
\frac{1}{4\epsilon_3}\|\widetilde{\partial}\omega_\theta^k\|^2 
+ \epsilon_3\|\widetilde{\rho}_\theta^{\,k}\|^2.\label{theta3}
\end{align}
\noindent\textbf{Term}: $-((\mathbf{a}\cdot\nabla)\widetilde{\omega}_u'^{\,k},\widetilde{\rho}_\theta^{\,k})$.
Using Green's theorem to transfer the gradient so as not to lose spatial convergence order, we obtain
$
-((\mathbf{a}\cdot\nabla)\widetilde{\omega}_u'^{\,k},\widetilde{\rho}_\theta^{\,k}) 
= (\widetilde{\omega}_u'^{\,k}, (\mathbf{a}\cdot\nabla)\widetilde{\rho}_\theta^{\,k}),
$
since $\widetilde{\rho}_\theta^{\,k} \in H^1_0(\Omega)$. Of this, and by Cauchy–Schwarz and $\epsilon$-Young,
\begin{align}
-((\mathbf{a}\cdot\nabla)\widetilde{\omega}_u'^{\,k},\widetilde{\rho}_\theta^{\,k})
=(\widetilde{\omega}_u'^{\,k}, (\mathbf{a}\cdot\nabla)\widetilde{\rho}_\theta^{\,k})
& \leqslant  
\|\mathbf{a}\|_\infty \|\widetilde{\omega}_u'^{\,k}\| \|\nabla\widetilde{\rho}_\theta^{\,k}\|
\leqslant
\frac{\|\mathbf{a}\|_\infty^2}{\epsilon_{4}}\|\widetilde{\omega}_u'^{\,k}\|^2 
+ \epsilon_{4}\|\nabla\widetilde{\rho}_\theta^{\,k}\|^2.\label{theta4}
\end{align}

\noindent\textbf{Nonlinear term with $\beta$}:
$-\Big[\beta\big(\hspace{-0.09cm}\int_\Omega \widetilde{\Theta}_h^{\,k}dx\big)(\nabla \widetilde{\theta}^{k},\nabla\widetilde{\rho}_{\theta}^k)
	- \beta\big(\hspace{-0.09cm}\int_\Omega\theta(t_{k-\frac{1}{2}})dx\big)(\nabla\theta(t_{k-\frac{1}{2}}),\nabla\widetilde{\rho}_{\theta}^k)\Big].$
Adding and subtracting  $\beta\big(\hspace{-0.09cm}\int_\Omega\theta(t_{k-\frac{1}{2}})dx\big)(\nabla\widetilde{\theta}^k,\nabla\widetilde{\rho}_\theta^{\,k})$ in the above term, the term to be estimated is:
\begin{align}
\mathcal{N}_\beta :=&-\Big[\beta\big(\hspace{-0.09cm}\int_\Omega\widetilde{\Theta}_h^{\,k}dx\big)-\beta\big(\hspace{-0.09cm}\int_\Omega\theta(t_{k-\frac{1}{2}})dx\big)\Big](\nabla\widetilde{\theta}^k,\nabla\widetilde{\rho}_{\theta}^k)
  - \beta\big(\hspace{-0.09cm}\int_\Omega\theta(t_{k-\frac{1}{2}})dx\big)(\nabla\widetilde{\theta}^k-\nabla\theta(t_{k-\frac{1}{2}}),\nabla\widetilde{\rho}_{\theta}^k). \label{beta0}
\end{align}
$\mathcal{N}_{\beta,1}:=-\Big[\beta\big(\int_\Omega\widetilde{\Theta}_h^{\,k}dx\big)-\beta\big(\int_\Omega\theta(t_{k-\frac{1}{2}})dx\big)\Big](\nabla\widetilde{\theta}^k,\nabla\widetilde{\rho}_{\theta}^k)$. 
By the hypothesis on the function $\beta$ in \eqref{hipoteses},
together with the Cauchy-Schwarz inequality,
\begin{align}
\Big|\beta\big(\hspace{-0.09cm}\int_\Omega\widetilde{\Theta}_h^{\,k}dx\big)-\beta\big(\hspace{-0.09cm}\int_\Omega\theta(t_{k-\frac{1}{2}})dx\big)\Big| 
\leqslant 
L_1 \Big|\int_\Omega\big(\widetilde{\Theta}_h^{\,k} - \theta(t_{k-\frac{1}{2}})\big)dx\Big|
\leqslant 
L_1|\Omega|^{1/2} \|\widetilde{\Theta}_h^{\,k}-\theta(t_{k-\frac{1}{2}})\|.\label{beta1}
\end{align}
The difference $\widetilde{\Theta}_h^{\,k} - \theta(t_{k-\frac{1}{2}})$ é is decomposed as
$
\widetilde{\Theta}_h^{\,k} - \theta(t_{k-\frac{1}{2}}) = \widetilde{\rho}_\theta^{\,k} + \widetilde{\omega}_\theta^{\,k} + (\widetilde{\theta}^k - \theta(t_{k-\frac{1}{2}})),
$    
where $\widetilde{\rho}_\theta^{\,k} = \widetilde{\Theta}_h^{\,k} - \widetilde{R_h^\theta\theta}^k$ and $\widetilde{\omega}_\theta^{\,k} = \widetilde{R_h^\theta\theta}^k - \widetilde{\theta}^k$.
By the Cauchy-Schwarz inequality,
\begin{align}
    |(\nabla\widetilde{\theta}^k,\nabla\widetilde{\rho}_{\theta}^k)| \leqslant \|\nabla\widetilde{\theta}^k\| \|\nabla\widetilde{\rho}_{\theta}^k\| \leqslant \|\theta\|_{L^\infty(0,T;H_0^1(\Omega))} \|\nabla\widetilde{\rho}_{\theta}^k\|.\label{beta2}
\end{align}
Thus, using \eqref{beta1} and \eqref{beta2}, the above decomposition, and the triangle inequality,
\begin{align}
&\mathcal{N}_{\beta,1} 
\leqslant
\big|\beta\big(\hspace{-0.09cm}\int_\Omega\widetilde{\Theta}_h^{\,k}dx\big)-\beta\big(\hspace{-0.09cm}\int_\Omega\theta(t_{k-\frac{1}{2}})dx\big)\big|
|(\nabla\widetilde{\theta}^k,\nabla\widetilde{\rho}_{\theta}^k)| 
\leqslant 
	L_1 |\Omega|^{1/2} \|\widetilde{\Theta}_h^{\,k}-\theta(t_{k-\frac{1}{2}})\| \|\nabla\widetilde{\theta}^k\| \|\nabla\widetilde{\rho}_{\theta}^k\| \nonumber\\
&\qquad 
	\leqslant 
	c_\beta \|\widetilde{\Theta}_h^{\,k}-\theta(t_{k-\frac{1}{2}})\|\|\nabla\widetilde{\rho}_{\theta}^k\| 
	\leqslant c_\beta \big( \|\widetilde{\rho}_\theta^{\,k}\| + \|\widetilde{\omega}_\theta^{\,k}\| + \|\widetilde{\theta}^k-\theta(t_{k-\frac{1}{2}})\| \big) \|\nabla\widetilde{\rho}_{\theta}^k\|,\label{eq:cont22}
\end{align}
where $c_\beta:=L_1 |\Omega|^{1/2} \|\theta\|_{L^\infty(0,T;H_0^1(\Omega))}<+\infty$, since, by \eqref{theta-regularidade}, $\theta \in H^2(0,T; H^1_0(\Omega)) \hookrightarrow L^\infty(0,T;H_0^1(\Omega))$ (Sobolev embedding).
In \eqref{eq:cont22}, applying the $\epsilon$-Young inequality
and subsequently the bound $(x+y+z)^2 \leqslant 3(x^2+y^2+z^2)$,
\begin{align}
&\mathcal{N}_{\beta,1} 
\leqslant  
\frac{3c_\beta^2}{4\epsilon_4}
\|\widetilde{\rho}_\theta^{\,k}\|^2 
+ \frac{3c_\beta^2}{4\epsilon_4}\|\widetilde{\omega}_\theta^{\,k}\|^2 
+ \frac{3c_\beta^2}{4\epsilon_4}\|\widetilde{\theta}^k-\theta(t_{k-\frac{1}{2}})\|^2 
+ \epsilon_4 \|\nabla\widetilde{\rho}_{\theta}^k\|^2. \label{beta4_1}
\end{align}

\noindent$
\mathcal{N}_{\beta,2}:=- \beta\big(\hspace{-0.09cm}\int_\Omega\theta(t_{k-\frac{1}{2}})dx\big)(\nabla\widetilde{\theta}^k-\nabla\theta(t_{k-\frac{1}{2}}),\nabla\widetilde{\rho}_{\theta}^k).
$
Again, by the hypothesis regarding the function $\beta$ in \eqref{hipoteses},
$
|\beta\big(\hspace{-0.09cm}\int_\Omega\theta(t_{k-\frac{1}{2}})dx\big)|-|\beta(0)|
\leqslant |\beta\big(\hspace{-0.09cm}\int_\Omega\theta(t_{k-\frac{1}{2}})dx\big)-\beta(0)|
\leqslant \widehat{L}_1 |\int_\Omega\theta(t_{k-\frac{1}{2}})dx|.$
Then, 
$
|\beta\big(\hspace{-0.09cm}\int_\Omega\theta(t_{k-\frac{1}{2}})dx\big)|
\leqslant |\beta(0)| + \widehat{L}_1|\int_\Omega\theta(t_{k-\frac{1}{2}})dx|.
$
By the Cauchy-Schwarz inequality,
$|\int_\Omega\theta(t_{k-\frac{1}{2}})dx| \leqslant |\Omega|^{1/2}\|\theta(t_{k-\frac{1}{2}})\|
\leqslant |\Omega|^{1/2}\|\theta\|_{C(0,T; L^2(\Omega))}.$

Now, defining the constant
$
\overline{c}_\beta := \widehat{L}_1 |\Omega|^{1/2}\|\theta\|_{C(0,T; L^2(\Omega))} + |\beta(0)|<\infty,
$ see \eqref{theta-regularidade}. So,
$|\beta\big(\hspace{-0.09cm}\int_\Omega\theta(t_{k-\frac{1}{2}})dx\big)|
\leqslant \overline{c}_\beta \hspace{-0.04cm}<\hspace{-0.04cm}\infty.$
Of this, and by Cauchy--Schwarz and $\epsilon$-Young inequalities,
\begin{align}
\mathcal{N}_{\beta,2} 
\leqslant 
|\mathcal{N}_{\beta,2}|
\leqslant
 \frac{\overline{c}_\beta^2}{4\epsilon_4} \|\nabla(\widetilde{\theta}^k - \theta(t_{k-\frac{1}{2}}))\|^2 
+ \epsilon_4 \|\nabla\widetilde{\rho}_\theta^{\,k}\|^2.\label{beta7}
\end{align}
Thus, returning to \eqref{beta0}, and using \eqref{beta4_1} and \eqref{beta7}, we conclude that
\begin{align}
\mathcal{N}_\beta \hspace{-0.07cm}=\hspace{-0.07cm} \mathcal{N}_{\beta,1} \hspace{-0.07cm}+\hspace{-0.07cm} \mathcal{N}_{\beta,2} 
\leqslant
& 
\frac{3c_\beta^2}{4\epsilon_4} \|\widetilde{\rho}_\theta^{\,k}\|^2 
\hspace{-0.07cm}+\hspace{-0.07cm} 2\epsilon_4 \|\nabla\widetilde{\rho}_{\theta}^k\|^2 
\hspace{-0.07cm}+\hspace{-0.07cm} \frac{3c_\beta^2}{4\epsilon_4}\|\widetilde{\omega}_\theta^{\,k}\|^2  
\hspace{-0.07cm}+\hspace{-0.07cm} \frac{3c_\beta^2}{4\epsilon_4}\|\widetilde{\theta}^k-\theta(t_{k-\frac{1}{2}})\|^2 
\hspace{-0.07cm}+\hspace{-0.07cm} \frac{\overline{c}_\beta^2}{4\epsilon_4} \|\nabla(\widetilde{\theta}^k - \theta(t_{k-\frac{1}{2}}))\|^2 .\label{theta_NL}
\end{align}

\noindent\textbf{Term}: $\mathcal{R}_\theta^k(\widetilde{\rho}_\theta^{\,k})$. Taking $\psi_h \hspace{-0.07cm}=\hspace{-0.07cm} \widetilde{\rho}_\theta^{\,k}$ in \eqref{resto_theta}, $
 \mathcal{R}_\theta^k(\widetilde{\rho}_\theta^{\,k})=
 -(\widetilde{\partial}\theta^{k} \hspace{-0.07cm}-\hspace{-0.07cm} \theta'(t_{k-\frac{1}{2}}),\widetilde{\rho}_\theta^{\,k})
\hspace{-0.07cm}-\hspace{-0.07cm} \big((\mathbf{a}\cdot\nabla)\big(\widetilde{u}'^k - u'(t_{k-\frac{1}{2}})\big),\widetilde{\rho}_\theta^{\,k}\big).$ 
Using the Cauchy–Schwarz and $\epsilon$-Young inequalities,
\begin{align}
& 
 -(\widetilde{\partial}\theta^{k} - \theta'(t_{k-\frac{1}{2}}),\widetilde{\rho}_\theta^{\,k})
\leqslant 
\frac{1}{4\epsilon_3}\|\widetilde{\partial}\theta^{k} - \theta'(t_{k-\frac{1}{2}})\|^2 
+\epsilon_3\|\widetilde{\rho}_\theta^{\,k}\|^2;
\\
& 
- \big((\mathbf{a}\cdot\nabla)\big(\widetilde{u}'^k - u'(t_{k-\frac{1}{2}})\big),\widetilde{\rho}_\theta^{\,k}\big)
\leqslant 
\frac{\|a\|^2_{\infty}}{4\epsilon_3}\|\nabla(\widetilde{u}'^k - u'(t_{k-\frac{1}{2}}))\|^2 
+\epsilon_3\|\widetilde{\rho}_\theta^{\,k}\|^2.\label{est-theta}
\end{align}
Therefore, in \eqref{eq:erro_theta}–\eqref{resto_delta}, using \eqref{theta1}–\eqref{est-theta} and rearranging the terms, we obtain
\begin{align}
& \frac{1}{2}\widetilde{\partial}\|\rho_\theta^k\|^2
\hspace{-0.09cm}+\hspace{-0.09cm} \frac{\beta_0}{2}\|\nabla \widetilde{\rho}_\theta^{\,k}\|^2
  \hspace{-0.07cm}\leqslant\hspace{-0.07cm} 
\big((\mathbf{a}\cdot\nabla)\widetilde{\rho}_\theta^{\,k}, \widetilde{\rho}_u'^{\,k}\big) 
\hspace{-0.09cm} +\hspace{-0.09cm} \frac{1}{4\epsilon_3}\|\widetilde{\partial}\omega_\theta^k\|^2 
\hspace{-0.09cm}+\hspace{-0.09cm} \Big(\frac{3c_\beta^2}{4\epsilon_4} \hspace{-0.09cm}+\hspace{-0.09cm} 3\epsilon_3\Big)\|\widetilde{\rho}_\theta^{\,k}\|^2
\hspace{-0.09cm}+\hspace{-0.09cm} 3\epsilon_{4}\|\nabla\widetilde{\rho}_\theta^{\,k}\|^2 
\hspace{-0.09cm}+\hspace{-0.09cm} \frac{3c_\beta^2}{4\epsilon_4}\|\widetilde{\theta}^k-\theta(t_{k-\frac{1}{2}})\|^2 
\nonumber
\\
&
\;\;
\hspace{-0.09cm}+\hspace{-0.09cm} \frac{\|\mathbf{a}\|_\infty^2}{\epsilon_{4}}\|\widetilde{\omega}_u'^{\,k}\|^2  
\hspace{-0.09cm}+\hspace{-0.09cm} \frac{3c_\beta^2}{4\epsilon_4}\|\widetilde{\omega}_\theta^{\,k}\|^2  
\hspace{-0.09cm}+\hspace{-0.09cm} \frac{\overline{c}_\beta^2}{4\epsilon_4} \|\nabla(\widetilde{\theta}^k \hspace{-0.09cm}-\hspace{-0.09cm} \theta(t_{k-\frac{1}{2}}))\|^2 
\hspace{-0.09cm}+\hspace{-0.09cm} \frac{1}{4\epsilon_3}\|\widetilde{\partial}\theta^{k} \hspace{-0.09cm}-\hspace{-0.09cm} \theta'(t_{k-\frac{1}{2}})\|^2
\hspace{-0.09cm}+\hspace{-0.09cm} \frac{\|a\|^2_{\infty}}{4\epsilon_3}\|\nabla(\widetilde{u}'^k \hspace{-0.09cm}-\hspace{-0.09cm} u'(t_{k-\frac{1}{2}}))\|^2 .
\label{resumo_theta}
\end{align}            

\subsubsection{Estimate for $\rho_\delta$}\label{sec:rho-delta}
We choose $\chi_h = \widetilde{\rho}_\delta'^{\,n}$ in \eqref{eq:erro_delta} and \eqref{resto_delta}.
In order to cancel the term $\mu^{k-\frac{1}{2}}(\widetilde{\rho}_u'^{\,k},\widetilde{\rho}_\delta'^{\,k})_{\Gamma_1}$ in \ref{sec:soma-estima} with the same term coming from the estimate for $\rho_u$, we consider the system of equations \eqref{eq:erro_delta}-\eqref{resto_delta} multiplied by $\mu^{k-\frac{1}{2}}$, that is,
\begin{align}
&  \mu^{k-\frac{1}{2}}(f_1\widetilde{\partial}\rho_\delta'^k,\widetilde{\rho}_\delta'^{\,k})_{\Gamma_1} 
+ \mu^{k-\frac{1}{2}}(f_2\widetilde{\rho}_\delta'^{\,k},\widetilde{\rho}_\delta'^{\,k})_{\Gamma_1} 
+ \mu^{k-\frac{1}{2}}(f_3\widetilde{\rho}_\delta^{\,k},\widetilde{\rho}_\delta'^{\,k})_{\Gamma_1} 
=
-\mu^{k-\frac{1}{2}}(\widetilde{\rho}_u'^{\,k},\widetilde{\rho}_\delta'^{\,k})_{\Gamma_1}
-\mu^{k-\frac{1}{2}}(\widetilde{\omega}_u'^{\,k},\widetilde{\rho}_\delta'^{\,k})_{\Gamma_1}  
\nonumber\\
&\quad 
- \mu^{k-\frac{1}{2}}(f_1\widetilde{\partial}\omega_\delta'^k,\widetilde{\rho}_\delta'^{\,k})_{\Gamma_1} 
- \mu^{k-\frac{1}{2}}(f_2\widetilde{\omega}_\delta'^{\,k},\widetilde{\rho}_\delta'^{\,k})_{\Gamma_1} 
- \mu^{k-\frac{1}{2}}(f_3\widetilde{\omega}_\delta^{\,k},\widetilde{\rho}_\delta'^{\,k})_{\Gamma_1} 
- \mu^{k-\frac{1}{2}}\mathcal{R}_\delta^k(\widetilde{\rho}_\delta'^{\,k}). \label{eq:erro_delta1}
\end{align}

\noindent\textbf{Term}: $\mu^{k-\frac{1}{2}}(f_1\widetilde{\partial}\rho_\delta'^k,\widetilde{\rho}_\delta'^{\,k})_{\Gamma_1}$. Using the definition of $\widetilde{\partial}\rho_\delta'^k = (\rho_\delta'^k - \rho_\delta'^{k-1})/\tau$ and $\widetilde{\rho}_\delta'^{\,k} = (\rho_\delta'^k + \rho_\delta'^{k-1})/2$, it follows that
\begin{align}
    \mu^{k-\frac{1}{2}}(f_1\widetilde{\partial}\rho_\delta'^k,\widetilde{\rho}_\delta'^{\,k})_{\Gamma_1} &= \frac{\mu^{k-\frac{1}{2}}}{2\tau}\big( \|f_1^{1/2}\rho_\delta'^k\|_{\Gamma_1}^2 - \|f_1^{1/2}\rho_\delta'^{k-1}\|_{\Gamma_1}^2 \big)= \mu^{k-\frac{1}{2}}\frac{1}{2}\widetilde{\partial}\|f_1^{1/2}\rho_\delta'^k\|_{\Gamma_1}^2.\label{delta1}
\end{align}

\noindent\textbf{Term}: $\mu^{k-\frac{1}{2}}(f_2\widetilde{\rho}_\delta'^{\,k},\widetilde{\rho}_\delta'^{\,k})_{\Gamma_1}$. Since, by \eqref{hipoteses}, $f_2(x) \geqslant \underline{f}_2 \geqslant 0$ for all $x\in\Gamma_1$ and $\mu(t) \geqslant \mu_0>0,$ we have:
\begin{align}
\mu^{k-\frac{1}{2}}(f_2\widetilde{\rho}_\delta'^{\,k},\widetilde{\rho}_\delta'^{\,k})_{\Gamma_1} \geqslant  \mu_0\underline{f}_2 \|\widetilde{\rho}_\delta'^{\,k}\|_{\Gamma_1}^2.\label{delta2}
\end{align}

\noindent\textbf{Term}: $\mu^{k-\frac{1}{2}}(f_3\widetilde{\rho}_\delta^{\,k},\widetilde{\rho}_\delta'^{\,k})_{\Gamma_1}$. To rewrite it, we decompose, $\widetilde{\rho}_\delta'^{\,k}$. From the definitions in \eqref{decomp-2}, and the linearity of the projection, we have:
$\widetilde{\rho}_\delta'^{\,k} = \widetilde{\partial}\rho_\delta^k + R_h^\delta(\widetilde{\partial}\delta^k - \widetilde{\delta'}^k)$. Therefore,
    \begin{align}
   \mu^{k-\frac{1}{2}} (f_3\widetilde{\rho}_\delta^{\,k},\widetilde{\rho}_\delta'^{\,k})_{\Gamma_1} 
   &= \mu^{k-\frac{1}{2}}(f_3\widetilde{\rho}_\delta^{\,k},\widetilde{\partial}\rho_\delta^k)_{\Gamma_1} 
   + \mu^{k-\frac{1}{2}}(f_3\widetilde{\rho}_\delta^{\,k}, R_h^\delta(\widetilde{\partial}\delta^k - \widetilde{\delta'}^k))_{\Gamma_1}.\label{delta3}
    \end{align}
The first term in \eqref{delta3} is a discrete derivative:
    $   \mu^{k-\frac{1}{2}} (f_3\widetilde{\rho}_\delta^{\,k},\widetilde{\partial}\rho_\delta^k)_{\Gamma_1} = \mu^{k-\frac{1}{2}}\frac{1}{2}\widetilde{\partial}\|f_3^{1/2}\rho_\delta^k\|_{\Gamma_1}^2.$
    The second term in \eqref{delta3} will be moved to the right-hand side. Using the Cauchy-Schwarz and $\epsilon$-Young inequalities,
\begin{align}
\big|\mu^{k-\frac{1}{2}}(f_3\widetilde{\rho}_\delta^{\,k}, R_h^\delta(\widetilde{\partial}\delta^k - \widetilde{\delta'}^k))_{\Gamma_1}\big| \leqslant 
\frac{\mu_1^2\overline{f}_3}{4\epsilon_5} \|\widetilde{\partial}\delta^k - \widetilde{\delta'}^k\|_{\Gamma_1}^2
+ \epsilon_5  \|\widetilde{\rho}_\delta^{\,k}\|_{\Gamma_1}^2,
\label{delta4}
\end{align}
where $\|R_h^\delta X\|_{\Gamma_1}\hspace{-0.07cm}\leqslant\hspace{-0.07cm} \|X\|_{\Gamma_1},\forall X \in L^2(\Gamma_1)$, see \eqref{eq:estab_delta}, and see the hypothesis regarding $f_3$ and  $\mu$,  in \eqref{hipoteses}.

\noindent\textbf{Term}: $-\mu^{k-\frac{1}{2}}(\widetilde{\rho}_u'^{\,k},\widetilde{\rho}_\delta'^{\,k})_{\Gamma_1} $. 
This term represents the coupling between $\rho_u'$ and $\rho_\delta'$; it will not be estimated, as it will be canceled later with an analogous term coming from the equation for $\rho_u$.

We estimate the following terms using the Cauchy-Schwarz and $\epsilon$-Young inequalities, always with the aim of isolating $ \|\widetilde{\rho}_\delta'^{\,k}\|_{\Gamma_1}^2$ and the projection errors.
    
\noindent\textbf{Terms}: $-\mu^{k-\frac{1}{2}}(\widetilde{\omega}_u'^{\,k},\widetilde{\rho}_\delta'^{\,k})_{\Gamma_1}$; $-\mu^{k-\frac{1}{2}}(f_1\widetilde{\partial}\omega_\delta'^k,\widetilde{\rho}_\delta'^{\,k})_{\Gamma_1}$; $-\mu^{k-\frac{1}{2}}(f_1\widetilde{\partial}\omega_\delta'^k,\widetilde{\rho}_\delta'^{\,k})_{\Gamma_1}$; $-\mu^{k-\frac{1}{2}}(f_2\widetilde{\omega}_\delta'^{\,k},\widetilde{\rho}_\delta'^{\,k})_{\Gamma_1}$. Hence,
\begin{align}
\begin{aligned}
& |\mu^{k-\frac{1}{2}}(\widetilde{\omega}_u'^{\,k},\widetilde{\rho}_\delta'^{\,k})_{\Gamma_1}| 
     \leqslant 
\frac{\mu_1^2}{4\epsilon_6}\|\widetilde{\omega}_u'^{\,k}\|_{\Gamma_1}^2 
\hspace{-0.07cm}+\hspace{-0.07cm} \epsilon_6\|\widetilde{\rho}_\delta'^{\,k}\|_{\Gamma_1}^2; 
\quad
|\mu^{k-\frac{1}{2}}(f_1\widetilde{\partial}\omega_\delta'^k,\widetilde{\rho}_\delta'^{\,k})_{\Gamma_1}| 
    \leqslant 
\frac{\mu_1^2\overline{f}_1}{4\epsilon_6}\|\widetilde{\partial}\omega_\delta'^k\|_{\Gamma_1}^2 
    \hspace{-0.07cm}+\hspace{-0.07cm} \epsilon_6\|\widetilde{\rho}_\delta'^{\,k}\|_{\Gamma_1}^2;
    \\
    & |\mu^{k-\frac{1}{2}}(f_2\widetilde{\omega}_\delta'^{\,k},\widetilde{\rho}_\delta'^{\,k})_{\Gamma_1}| 
    \leqslant 
    \frac{\mu_1^2\overline{f}_2}{4\epsilon_6}\|\widetilde{\omega}_\delta'^{\,k}\|_{\Gamma_1}^2 
    \hspace{-0.07cm}+\hspace{-0.07cm} \epsilon_6\|\widetilde{\rho}_\delta'^{\,k}\|_{\Gamma_1}^2;
\quad
    |\mu^{k-\frac{1}{2}}(f_3\widetilde{\omega}_\delta^{\,k},\widetilde{\rho}_\delta'^{\,k})_{\Gamma_1}| 
    \leqslant 
    \frac{\mu_1^2\overline{f}_3}{4\epsilon_6}\|\widetilde{\omega}_\delta^{\,k}\|_{\Gamma_1}^2 
    \hspace{-0.07cm}+\hspace{-0.07cm} \epsilon_6\|\widetilde{\rho}_\delta'^{\,k}\|_{\Gamma_1}^2.
\end{aligned}
\label{delta5+}
\end{align}

\noindent\textbf{Term}: $\mu^{k-\frac{1}{2}}\mathcal{R}_{\delta}^k(\widetilde{\rho}_{\delta}'^{\,k})$. From \eqref{resto_delta},
$\mu^{k-\frac{1}{2}}\mathcal{R}_\delta^k(\widetilde{\rho}_{\delta}'^{\,k})
=
- \mu^{k-\frac{1}{2}}(\widetilde{u}'^k - u'(t_{k-\frac{1}{2}}), \widetilde{\rho}_{\delta}'^{\,k})_{\Gamma_1}
- \mu^{k-\frac{1}{2}}\big(f_1(\widetilde{\partial}\delta'^{k} - \delta''(t_{k-\frac{1}{2}})) , \widetilde{\rho}_{\delta}'^{\,k}\big)_{\Gamma_1}
- \mu^{k-\frac{1}{2}}\big(f_2(\widetilde{\delta}'^k - \delta'(t_{k-\frac{1}{2}})) , \widetilde{\rho}_{\delta}'^{\,k}\big)_{\Gamma_1}.$
Using hypotheses for $\mu$ and $f_i$ (see \eqref{hipoteses}) and by Cauchy–Schwarz and $\epsilon$-Young,
\begin{align}
& 
 - \mu^{k-\frac{1}{2}}(\widetilde{u}'^k - u'(t_{k-\frac{1}{2}}), \widetilde{\rho}_{\delta}'^{\,k})_{\Gamma_1}
 \leqslant 
 \frac{\mu_1^2}{4\epsilon_6}\|\widetilde{u}'^k - u'(t_{k-\frac{1}{2}})\|_{\Gamma_1}^2 
 +\epsilon_6\|\widetilde{\rho}_{\delta}'^{\,k}\|_{\Gamma_1}^2;
\\
&
- \big(f_1(\widetilde{\partial}\delta'^{k} - \delta''(t_{k-\frac{1}{2}})) , \widetilde{\rho}_{\delta}'^{\,k}\big)_{\Gamma_1}
 \leqslant 
 \frac{\overline{f}_1\mu_1^2}{4\epsilon_6}\|\widetilde{\partial}\delta'^{k} - \delta''(t_{k-\frac{1}{2}})\|_{\Gamma_1}^2 
 +\epsilon_6\|\widetilde{\rho}_{\delta}'^{\,k}\|_{\Gamma_1}^2;
\\
&
- \mu^{k-\frac{1}{2}}\big(f_2(\widetilde{\delta}'^k - \delta'(t_{k-\frac{1}{2}})) , \widetilde{\rho}_{\delta}'^{\,k}\big)_{\Gamma_1}
 \leqslant 
 \frac{\overline{f}_2\mu_1^2}{4\epsilon_6}\|\widetilde{\delta}'^k - \delta'(t_{k-\frac{1}{2}})\|_{\Gamma_1}^2 
 +\epsilon_6\|\widetilde{\rho}_{\delta}'^{\,k}\|_{\Gamma_1}^2;
 \\
&
- \mu^{k-\frac{1}{2}}\big(f_3(\widetilde{\delta}^k - \delta(t_{k-\frac{1}{2}})),\widetilde{\rho}_{\delta}'^{\,k}\big)_{\Gamma_1}
 \leqslant 
 \frac{\overline{f}_3\mu_1^2}{4\epsilon_6}\|\widetilde{\delta}^k - \delta(t_{k-\frac{1}{2}})\|_{\Gamma_1}^2 
 +\epsilon_6\|\widetilde{\rho}_{\delta}'^{\,k}\|_{\Gamma_1}^2.\label{taylor_delta}
\end{align}

After all the estimations, returning to \eqref{eq:erro_delta1}, which is linked to \eqref{eq:erro_delta}-\eqref{resto_delta}, and using \eqref{delta1}-\eqref{taylor_delta}, noting that $\mu_0 \underline{f}_2 \|\widetilde{\rho}_\delta'^{\,k}\|_{\Gamma_1}^2 \geqslant 0$, we obtain the estimate:
\begin{align}
&\mu^{k-\frac{1}{2}}\frac{1}{2}\widetilde{\partial}\|f_1^{1/2}\rho_\delta'^k\|_{\Gamma_1}^2 
+  \mu^{k-\frac{1}{2}}\frac{1}{2}\widetilde{\partial}\|f_3^{1/2}\rho_\delta^k\|_{\Gamma_1}^2
\leqslant
-\mu^{k-\frac{1}{2}}(\widetilde{\omega}_u'^{\,k},\widetilde{\rho}_\delta'^{\,k})_{\Gamma_1}
+ 8\epsilon_6\|\widetilde{\rho}_\delta'^{\,k}\|_{\Gamma_1}^2 
+ \epsilon_5  \|\widetilde{\rho}_\delta^{\,k}\|_{\Gamma_1}^2 
+ \frac{\mu_1^2\overline{f}_1}{4\epsilon_6}\|\widetilde{\partial}\omega_\delta'^k\|_{\Gamma_1}^2\nonumber
\\
&\quad  
+ \frac{\mu_1^2\overline{f}_2}{4\epsilon_6}\|\widetilde{\omega}_\delta'^{\,k}\|_{\Gamma_1}^2
+ \frac{\mu_1^2\overline{f}_3}{4\epsilon_6}\|\widetilde{\omega}_\delta^{\,k}\|_{\Gamma_1}^2
+\frac{\mu_1^2}{4\epsilon_6}\|\widetilde{\omega}_u'^{\,k}\|_{\Gamma_1}^2
+ \frac{\mu_1^2}{4\epsilon_6}\|\widetilde{u}'^k - u'(t_{k-\frac{1}{2}})\|_{\Gamma_1}^2 
+ \frac{\mu_1^2\overline{f}_3}{4\epsilon_5} \|\widetilde{\partial}\delta^k - \widetilde{\delta}'^k\|_{\Gamma_1}^2
\nonumber 
\\ 
&\quad 
+ \frac{\overline{f}_1\mu_1^2}{4\epsilon_6}\|\widetilde{\partial}\delta'^{k} - \delta''(t_{k-\frac{1}{2}})\|_{\Gamma_1}^2
+ \frac{\overline{f}_2\mu_1^2}{4\epsilon_6}\|\widetilde{\delta}'^k - \delta'(t_{k-\frac{1}{2}})\|_{\Gamma_1}^2 
+ \frac{\overline{f}_3\mu_1^2}{4\epsilon_6}\|\widetilde{\delta}^k - \delta(t_{k-\frac{1}{2}})\|_{\Gamma_1}^2.
\label{resumo_delta}
\end{align}

\subsubsection{Summation of the estimates and application of the discrete Gronwall lemma}\label{sec:soma-estima}
We sum the three inequalities \eqref{resumo_u}, \eqref{resumo_theta}, and \eqref{resumo_delta} and use the cancellation of the coupling terms. Choosing $\epsilon_2=\frac{\mu_0g_0}{8}>0$ and $\epsilon_4=\frac{\beta_0}{6}>0$, we obtain the recursion:
\begin{align}
&\frac{1}{2}\widetilde{\partial}\Big[\|\rho_u'^k\|^2 
+\|\rho_\theta^k\|^2 \Big]
+ \mu^{k-\frac{1}{2}}\frac{1}{2}\widetilde{\partial}
\Big[
\|\nabla\rho_u^k\|^2 
+\|f_1^{1/2}\rho_\delta'^k\|_{\Gamma_1}^2 
+\|f_3^{1/2}\rho_\delta^k\|_{\Gamma_1}^2\Big]
\leqslant 
\widetilde{\partial}\big((\mathbf{a}\cdot \nabla)\rho_u^k,\widetilde{\omega}^k_\theta\big) 
+ c_1\|\nabla \widetilde{\rho}_{u}^k\|^2\nonumber 
\\
& \quad
+ c_1\|\nabla\rho_u^{k-1}\|^2
+ c_1\|\widetilde{\rho}_u'^{\,k}\|^2
+ c_1\|\widetilde{\rho}_\theta^{\,k}\|^2
+ c_1\|\widetilde{\rho}_\delta'^{\,k}\|_{\Gamma_1}^2 
+ c_1  \|\widetilde{\rho}_\delta^{\,k}\|_{\Gamma_1}^2  
+ \mathcal{F}_u^k 
+ \mathcal{F}_\theta^k 
+ \mathcal{F}_\delta^k
+ \mathcal{\overline{R}}_u^k 
+ \mathcal{\overline{R}}_{\theta}^k
+ \mathcal{\overline{R}}_{\delta}^k,\label{des_geral1}
\end{align}
with $c_1=\max\big\{\epsilon_0 + (3c_f^2c_P4)/\epsilon_1, 1, 4\epsilon_1, (3c_\beta^2)/4\epsilon_4 + 3\epsilon_3, 8\epsilon_6, \epsilon_5 \big\}>0,$ \,and
\begin{align}
&\mathcal{F}_u^k := 
\frac{1}{4\epsilon_1}\|\widetilde{\partial}\omega_u'^k\|^2 
+ \frac{3c_f^2}{4\epsilon_1}\|\widetilde{\omega}_u^{\,k}\|^2
+ \frac{\|\mathbf{a}\|_\infty^2}{\epsilon_{4}}\|\widetilde{\omega}_u'^{\,k}\|^2 
+\Big(\frac{L^2_g\mu_1^{2}}{2\epsilon_2} + \frac{\mu_1^2}{4\epsilon_6}\Big)\|\widetilde{\omega}_u'^{\,k}\|_{\Gamma_1}^2;\label{f_u}
\\
&\mathcal{F}_\theta^k := 
\frac{\|\mathbf{a}\|_{\infty}^2}{4}\|\widetilde{\partial}\widetilde{\omega}^k_\theta\|^2 
+ \frac{1}{4\epsilon_3}\|\widetilde{\partial}\omega_\theta^k\|^2 
+ \Big(\frac{3c_\beta^2}{4\epsilon_4}+1 \Big)\|\widetilde{\omega}_\theta^{\,k}\|^2;\label{f_theta}
\\
&\mathcal{F}_\delta^k := 
\frac{\mu_1^2\overline{f}_1}{4\epsilon_6}\|\widetilde{\partial}\omega_\delta'^k\|_{\Gamma_1}^2 
+ \Big(\frac{\mu_1^{2}}{4\epsilon_2}
+ \frac{\mu_1^2\overline{f}_2}{4\epsilon_6}\Big)\|\widetilde{\omega}_\delta'^{\,k}\|_{\Gamma_1}^2
+ \frac{\mu_1^2\overline{f}_3}{4\epsilon_6}\|\widetilde{\omega}_\delta^{\,k}\|_{\Gamma_1}^2;\label{f_delta}
\end{align}
where $\mathcal{F}_u^k, \mathcal{F}_\theta^k, \mathcal{F}_\delta^k $ gathers the projection error terms, and $\mathcal{\overline{R}}_u^k, \mathcal{\overline{R}}_{\theta}^k, \mathcal{\overline{R}}_{\delta}^k $  the Taylor truncation errors, as follows:
\begin{align}
&\mathcal{\overline{R}}_u^k  
\hspace{-0.07cm} := \hspace{-0.07cm}\frac{3c_f^2}{4\epsilon_1}\|\widetilde{u}^k \hspace{-0.07cm}-\hspace{-0.07cm} u(t_{k-\frac{1}{2}})\|^2
\hspace{-0.07cm}+\hspace{-0.07cm} \frac{1}{4\epsilon_1}\|\widetilde{\partial}u'^{k} \hspace{-0.05cm}-\hspace{-0.05cm} u''(t_{k-\frac{1}{2}})\|^2
\hspace{-0.07cm}+\hspace{-0.07cm} \Big(\frac{\mu_1^{2}}{4\epsilon_0}\hspace{-0.07cm}+\hspace{-0.07cm}\frac{\|\mathbf{a}\|_{\infty}^2}{4}\Big)\|\nabla(\widetilde{\partial}u^k \hspace{-0.05cm}-\hspace{-0.05cm}\widetilde{u}'^k)\|^2
\hspace{-0.07cm}+\hspace{-0.07cm} \frac{\mu_1^2}{4\epsilon_1}\|\Delta (\widetilde{u}^k \hspace{-0.07cm}-\hspace{-0.07cm}u(t_{k-\frac{1}{2}}))\|^2\nonumber 
\\
& \qquad  
\hspace{-0.05cm}+\hspace{-0.05cm} \frac{\mu_1^2}{4\epsilon_2}\|\frac{\partial}{\partial \nu}(\widetilde{u}^k \hspace{-0.05cm}-\hspace{-0.05cm} u(t_{k-\frac{1}{2}}))\|^2_{\Gamma_1}
\hspace{-0.05cm}+\hspace{-0.05cm} \frac{\|a\|^2_{\infty}}{4\epsilon_3}\|\nabla(\widetilde{u}'^k \hspace{-0.05cm}-\hspace{-0.05cm} u'(t_{k-\frac{1}{2}}))\|^2 +\Big(\frac{\mu_1^2}{4\epsilon_6} 
\hspace{-0.05cm}+\hspace{-0.05cm} \frac{L^2_g\mu_1^{2}}{2\epsilon_2}\Big)\|\widetilde{u}'^k \hspace{-0.05cm}-\hspace{-0.05cm} u'(t_{k-\frac{1}{2}})\|_{\Gamma_1}^2;
\label{resto_u-}
\\
&  \mathcal{\overline{R}}_{\theta}^k
\hspace{-0.07cm}:=\hspace{-0.07cm} 
\frac{1}{4\epsilon_3}\|\widetilde{\partial}\theta^{k} \hspace{-0.07cm}-\hspace{-0.07cm} \theta'(t_{k-\frac{1}{2}})\|^2
\hspace{-0.07cm}+\hspace{-0.07cm} \frac{3c_\beta^2}{4\epsilon_4}\|\widetilde{\theta}^k \hspace{-0.07cm}-\hspace{-0.07cm}\theta(t_{k-\frac{1}{2}})\|^2 
\hspace{-0.07cm}+\hspace{-0.07cm} \Big(\frac{\|a\|^2_{\infty}}{4\epsilon_1} \hspace{-0.07cm}+\hspace{-0.07cm} \frac{\overline{c}_\beta^2}{4\epsilon_4}\Big)\|\nabla(\widetilde{\theta}^k \hspace{-0.07cm}- \hspace{-0.07cm}\theta(t_{k-\frac{1}{2}}))\|^2;\label{resto_theta-}
\\
& \mathcal{\overline{R}}_{\delta}^k
\hspace{-0.07cm}:=\hspace{-0.07cm} \frac{\overline{f}_3\mu_1^2}{4\epsilon_6}\|\widetilde{\delta}^k \hspace{-0.07cm}-\hspace{-0.07cm} \delta(t_{k-\frac{1}{2}})\|_{\Gamma_1}^2
\hspace{-0.09cm}+\hspace{-0.09cm} \frac{\mu_1^2\overline{f}_3}{4\epsilon_5} \|\widetilde{\partial}\delta^k \hspace{-0.07cm}-\hspace{-0.07cm} \widetilde{\delta}'^k\|_{\Gamma_1}^2
\hspace{-0.09cm}+\hspace{-0.09cm} \frac{\overline{f}_1\mu_1^2}{4\epsilon_6}\|\widetilde{\partial}\delta'^{k} \hspace{-0.07cm}-\hspace{-0.07cm} \delta''(t_{k-\frac{1}{2}})\|_{\Gamma_1}^2\nonumber
\\
&\qquad   \hspace{-0.09cm}+\hspace{-0.09cm} \Big(\hspace{-0.09cm}\frac{\overline{f}_2\mu_1^2}{4\epsilon_6}\hspace{-0.09cm}+\hspace{-0.09cm}\frac{\mu_1^2}{4\epsilon_2}\hspace{-0.09cm}\Big)\|\widetilde{\delta}'^k \hspace{-0.07cm}-\hspace{-0.07cm} \delta'(t_{k-\frac{1}{2}})\|_{\Gamma_1}^2.\label{resto_delta-}
\end{align}

\noindent\textit{Preparation for the Gronwall lemma}.
In order to simplify the calculations later on, consider
\begin{align}
\mathcal{E}^k := \|\rho_u'^k\|^2 + \|\rho_\theta^k\|^2 + \|\nabla\rho_u^k\|^2 + \|\rho_\delta'^k\|_{\Gamma_1}^2 + \|\rho_\delta^k\|_{\Gamma_1}^2.\label{energia}
\end{align}

We multiply \eqref{des_geral1} by $2\tau$ and sum from $k=1$ to $m$, with $m\in\{0,1, \ldots, N\}$,
\begin{align}
&\sum_{k=1}^m \hspace{-0.09cm}\tau\,\overline{\partial}\big(\|\rho_u'^k\|^2 \hspace{-0.09cm}+\hspace{-0.09cm}\|\rho_\theta^k\|^2\big)
\hspace{-0.09cm}+\hspace{-0.09cm} \sum_{k=1}^m \hspace{-0.09cm}\mu^{k-\frac{1}{2}}\tau\,\overline{\partial}\Upsilon^k
\hspace{-0.09cm}\leqslant \hspace{-0.09cm}
\sum_{k=1}^m \hspace{-0.09cm}2\tau\,\overline{\partial}\big((\mathbf{a}\cdot\nabla)\rho_u^k,\widetilde{\omega}_\theta^k\big) 
\hspace{-0.07cm}+\hspace{-0.07cm} c_1\sum_{k=1}^m \hspace{-0.07cm}2\tau \|\nabla \widetilde{\rho}_{u}^k\|^2 
\hspace{-0.07cm}+\hspace{-0.07cm} c_1\sum_{k=1}^m \hspace{-0.07cm}2\tau \|\nabla\rho_u^{k-1}\|^2
\nonumber
\\
&
\;\;
\hspace{-0.09cm}+\hspace{-0.09cm} c_1\sum_{k=1}^m \hspace{-0.09cm}2\tau \|\widetilde{\rho}_u'^{\,k}\|^2
\hspace{-0.09cm}+\hspace{-0.09cm} c_1 \sum_{k=1}^m \hspace{-0.09cm}2\tau\|\widetilde{\rho}_\theta^{\,k}\|^2
\hspace{-0.09cm}+\hspace{-0.09cm} c_1\sum_{k=1}^m \hspace{-0.09cm}2\tau\|\widetilde{\rho}_\delta'^{\,k}\|_{\Gamma_1}^2
\hspace{-0.09cm}+\hspace{-0.09cm} c_1\sum_{k=1}^m \hspace{-0.09cm}2\tau\|\widetilde{\rho}_\delta^{\,k}\|_{\Gamma_1}^2 
\hspace{-0.09cm}+\hspace{-0.09cm} \sum_{k=1}^m \hspace{-0.09cm}2\tau 
\big[
 \mathcal{F}_u^k 
\hspace{-0.09cm}+\hspace{-0.09cm} \mathcal{F}_\theta^k 
\hspace{-0.09cm}+\hspace{-0.09cm} \mathcal{F}_\delta^k
\hspace{-0.09cm}+\hspace{-0.09cm} \mathcal{\overline{R}}_u^k 
\hspace{-0.09cm}+\hspace{-0.09cm} \mathcal{\overline{R}}_{\theta}^k
\hspace{-0.09cm}+\hspace{-0.09cm} \mathcal{\overline{R}}_{\delta}^k
\big],\label{soma_1}
\end{align}
where $\Upsilon^k = \|\nabla\rho_u^k\|^2 + \|f_1^{1/2}\rho_\delta'^k\|_{\Gamma_1}^2 + \|f_3^{1/2}\rho_\delta^k\|_{\Gamma_1}^2$.

\noindent{\it Left-hand side of \eqref{soma_1} (telescoping sum)}: Since $\tau\,\overline{\partial}c^k = c^k - c^{k-1}$, we obtain
\begin{align}
\sum_{k=1}^m \tau\,\widetilde{\partial}\big(\|\rho_u'^k\|^2+\|\rho_\theta^k\|^2\big) = \big(\|\rho_u'^m\|^2+\|\rho_\theta^m\|^2\big) - \big(\|\rho_u'^0\|^2+\|\rho_\theta^0\|^2\big).\label{tele1}
\end{align}
We apply the summation-by-parts identity to:
\begin{align}
\sum_{k=1}^m \mu^{k-\frac{1}{2}}\tau\,\overline{\partial}\Upsilon^k = \sum_{k=1}^m \mu^{k-\frac{1}{2}}(\Upsilon^k \hspace{-0.09cm}-\hspace{-0.09cm} \Upsilon^{k-1}) &= \mu^{m-\frac{1}{2}}\Upsilon^m \hspace{-0.09cm}-\hspace{-0.09cm} \mu^{\frac{1}{2}}\Upsilon^0
 \hspace{-0.09cm}-\hspace{-0.09cm} \sum_{k=1}^{m-1}(\mu^{k+\frac{1}{2}}\hspace{-0.09cm}-\hspace{-0.09cm}\mu^{k-\frac{1}{2}})\Upsilon^k.
\label{tele2}
\end{align}
The term $\mu^{m-\frac{1}{2}} \Upsilon^m$ is positive and will be kept on the left-hand side of \eqref{soma_1}; the terms $-\mu^{\frac{1}{2}} \Upsilon^0$ (constant) and $-\sum_{k=1}^{m-1}(\mu^{k+\frac{1}{2}}-\mu^{k-\frac{1}{2}})\Upsilon^k$ will be moved to the right-hand side of \eqref{soma_1}. Then, since $\mu\in C^1([0,T])$ (see \eqref{hipoteses}), there exists $L_\mu = \|\mu'\|_\infty$ such that $|\mu^{k+\frac{1}{2}}-\mu^{k-\frac{1}{2}}| \leqslant L_\mu\tau$ by the Mean Value Theorem. Furthermore, given the hypotheses on $f_i$ (see \eqref{hipoteses}), then
\begin{align}
\Big|\sum_{k=1}^{m-1}(\mu^{k+\frac{1}{2}}-\mu^{k-\frac{1}{2}})\Upsilon^k\Big| 
&\leqslant 
\sum_{k=1}^{m-1}|\mu^{k+\frac{1}{2}}-\mu^{k-\frac{1}{2}}|\Upsilon^k 
\leqslant  
L_\mu\tau\sum_{k=1}^{m-1}\Upsilon^k
\leqslant  
L_\mu\tau\sum_{k=0}^{m-1}\Upsilon^k
\nonumber
\\
&  
\leqslant
L_\mu\tau\sum_{k=0}^{m-1}\big(\|\nabla\rho_u^k\|^2 
+ \overline{f}_1\|\rho_\delta'^k\|_{\Gamma_1}^2 
+ \overline{f}_3\|\rho_\delta^k\|_{\Gamma_1}^2 \big) .\label{tele3}
\end{align}

\noindent{\it Right-hand side of \eqref{soma_1}}: The term $2\tau\,\widetilde{\partial}\big((\mathbf{a}\cdot\nabla)\rho_u^k,\widetilde{\omega}_\theta^k\big)$ is telescopic:
\begin{align}
2\sum_{k=1}^m \tau\,\widetilde{\partial}\big((\mathbf{a}\cdot\nabla)\rho_u^k,\widetilde{\omega}_\theta^k\big) = 2\Big[\big((\mathbf{a}\cdot\nabla)\rho_u^m,\widetilde{\omega}_\theta^m\big) 
- \big((\mathbf{a}\cdot\nabla)\rho_u^0,\widetilde{\omega}_\theta^0\big)\Big].
\end{align}
Using the Cauchy–Schwarz and $\epsilon$-Young inequalities, and noting that $\mu_0>0$ and $\mu_1>0$ (see \eqref{hipoteses}),
\begin{align}
&2| \big((\mathbf{a}\cdot\nabla)\rho_u^m,\widetilde{\omega}_\theta^m\big)| 
\leqslant 
2\|\mathbf{a}\|_\infty \|\nabla\rho_u^m\|\|\widetilde{\omega}_\theta^m\| 
\leqslant 
\frac{\mu_0}{2}\|\nabla\rho_u^m\|^2 + \frac{2\|\mathbf{a}\|^2_\infty}{\mu_0}\|\widetilde{\omega}_\theta^m\|^2,
\\
& 2|\big((\mathbf{a}\cdot\nabla)\rho_u^0,\widetilde{\omega}_\theta^0\big)| 
\leqslant
2\|\mathbf{a}\|_\infty \|\nabla\rho_u^0\|\|\widetilde{\omega}_\theta^0\|  
\leqslant 
\mu_1\|\nabla\rho_u^0\|^2 + \frac{\|\mathbf{a}\|^2_\infty}{\mu_1}\|\widetilde{\omega}_\theta^0\|^2,
\end{align}
where the term $\frac{\mu_0}{2}\|\nabla\rho_u^m\|^2$ will be moved to the left-hand side of \eqref{soma_1}; $\mu_1\|\nabla\rho_u^0\|^2$ o the right-hand side of \eqref{soma_1}; and $\frac{2\|\mathbf{a}\|^2_\infty}{\mu_0}\|\widetilde{\omega}_\theta^m\|^2$  and $\frac{\|\mathbf{a}\|^2_\infty}{\mu_1}\|\widetilde{\omega}_\theta^0\|^2$ will be incorporated into \eqref{f_theta}, which is redefined as:
\begin{align}
\mathcal{\overline{F}}_\theta^k := 
\frac{\|\mathbf{a}\|_{\infty}^2}{4}\|\widetilde{\partial}\widetilde{\omega}^k_\theta\|^2 
\hspace{-0.07cm}+\hspace{-0.07cm} \frac{1}{4\epsilon_3}\|\widetilde{\partial}\omega_\theta^k\|^2 
\hspace{-0.07cm}+\hspace{-0.07cm} \Big(\frac{3c_\beta^2}{4\epsilon_4}+1 \Big)\|\widetilde{\omega}_\theta^{\,k}\|^2
\hspace{-0.07cm}+\hspace{-0.07cm} \frac{2\|\mathbf{a}\|^2_\infty}{\mu_0}\|\widetilde{\omega}_\theta^m\|^2
\hspace{-0.07cm}+\hspace{-0.07cm} \frac{\|\mathbf{a}\|^2_\infty}{\mu_1}\|\widetilde{\omega}_\theta^0\|^2.
\label{f_theta1}
\end{align}

\noindent{\it Bounding of the averages:} We will use the following bound for the averages on the right-hand side of \eqref{soma_1}: for any function $Y$, note that
\begin{align}
\|\widetilde{Y}^{\,k}\|^2=\Big\|\frac{Y^{\,k}+Y^{\,k-1}}{2}\Big\|^2 \leqslant \frac{1}{2}\big(\|Y^k\|^2+\|Y^{k-1}\|^2\big).\label{padrão1}
\end{align}
Using \eqref{padrão1},
\begin{align}
&c_1\sum_{k=1}^m 2\tau \|\nabla \widetilde{\rho}_{u}^k\|^2 
\leqslant 
c_1\sum_{k=1}^m \tau \big(\|\nabla\rho_u^k\|^2+\|\nabla\rho_u^{k-1}\|^2\big)\nonumber
\\
&\qquad
\leqslant
c_1
\tau\Big(\|\nabla\rho_u^0\|^2  \hspace{-0.09cm}+ \hspace{-0.09cm} \sum_{k=1}^{m-1} \|\nabla\rho_u^k\|^2 
\hspace{-0.09cm} + \hspace{-0.09cm} \|\nabla\rho_u^m\|^2
\hspace{-0.09cm}+ \hspace{-0.09cm} \sum_{k=0}^{m-1} \|\nabla\rho_u^k\|^2\Big)
=2c_1\tau \sum_{k=0}^{m-1} \|\nabla\rho_u^k\|^2 
+ c_1\tau \|\nabla\rho_u^m\|^2.\label{media1}
\end{align}
Analogously to \eqref{media1}, the remaining terms follow, 
\begin{align}
\begin{aligned}
& c_1\sum_{k=1}^m 2\tau \|\nabla\rho_u^{k-1}\|^2= 2c_1\tau \sum_{k=0}^{m-1}  \|\nabla\rho_u^{k}\|^2; 
\qquad
 c_1\sum_{k=1}^m 2\tau \|\widetilde{\rho}_u'^{\,k}\|^2;
\leqslant 
2c_1 \tau \sum_{k=0}^{m-1} \|\rho_u'^k\|^2 
+ c_1\tau\|\rho_u'^m\|^2;
\\
&
c_1\sum_{k=1}^m 2\|\widetilde{\rho}_\theta^{\,k}\|^2
\leqslant 
2c_1\tau\sum_{k=0}^{m-1} \|\rho_\theta^k\|^2
+ c_1\tau \|\rho_\theta^m\|^2;
\qquad
 c_1\sum_{k=1}^m 2\tau\|\widetilde{\rho}_\delta'^{\,k}\|_{\Gamma_1}^2 
\leqslant 
2c_1 \tau\sum_{k=0}^{m-1}\|\rho_\delta'^k\|_{\Gamma_1}^2
+ c_1 \tau \|\rho_\delta'^m\|_{\Gamma_1}^2;
\\
&
c_1  \sum_{k=1}^m 2\tau\|\widetilde{\rho}_\delta^{\,k}\|_{\Gamma_1}^2
\leqslant 
2c_1\tau\sum_{k=0}^{m-1}\|\rho_\delta^k\|_{\Gamma_1}^2
+ c_1\tau \|\rho_\delta^m\|_{\Gamma_1}^2.
\end{aligned}\label{media_geral}
\end{align}
Since $0< \mu_0 \leqslant \mu^{k-\frac{1}{2}} \leqslant \mu_1 $) and $ 0<\underline{f}_i \leqslant f_i(x) \leqslant \overline{f}_i $ (see \eqref{hipoteses}), using \eqref{tele1}--\eqref{media_geral} in \eqref{soma_1} and choosing \( \tau \) such that
\begin{align}
0<\tau < \tau_0=\min \Big\{\frac{1}{c_1}, \frac{\mu_0}{2c_1}, \frac{\mu_0\underline{f}_1}{c_1}, \frac{\mu_0\underline{f}_3}{c_1}\Big\},\label{tau}
\end{align}
so,
\begin{align}
& (1-c_1\tau)\|\rho_u'^m\|^2 
\hspace{-0.03cm}+\hspace{-0.03cm} (1\hspace{-0.03cm}-\hspace{-0.03cm}c_1\tau)\|\rho_\theta^m\|^2 
\hspace{-0.03cm}+\hspace{-0.03cm} \Big(\frac{\mu_0}{2}\hspace{-0.03cm}- \hspace{-0.03cm}\tau c_1\Big)\|\nabla\rho_u^m\|^2
\hspace{-0.03cm}+\hspace{-0.03cm} (\mu_0\underline{f}_1\hspace{-0.03cm}-\hspace{-0.03cm}c_1\tau)\|\rho_\delta'^m\|_{\Gamma_1}^2\nonumber
\\     
& 
\quad + (\mu_0\underline{f}_3\hspace{-0.07cm}-\hspace{-0.07cm}c_1\tau)\|\rho_\delta^m\|_{\Gamma_1}^2
\leqslant 
\|\rho_u'^0\|^2+\|\rho_\theta^0\|^2
+\mu_1\big(2\|\nabla\rho_u^0\|^2 
+ \overline{f}_1\|\rho_\delta'^0\|_{\Gamma_1}^2 
+ \overline{f}_3\|\rho_\delta^0\|_{\Gamma_1}^2\big)\nonumber
\\
& 
\quad 
+2c_1\tau\sum_{k=0}^{m-1} \|\rho_\theta^k\|^2
\hspace{-0.07cm}+\hspace{-0.07cm} (2c_1\hspace{-0.07cm}+\hspace{-0.07cm}L_\mu\overline{f}_1) \tau\sum_{k=0}^{m-1}\|\rho_\delta'^k\|_{\Gamma_1}^2
\hspace{-0.07cm}+\hspace{-0.07cm} (2c_1\hspace{-0.07cm}+\hspace{-0.07cm}L_\mu\overline{f}_3)\tau\sum_{k=0}^{m-1}\|\rho_\delta^k\|_{\Gamma_1}^2
\hspace{-0.07cm}+\hspace{-0.07cm}2c_1\tau \sum_{k=0}^{m-1} \|\rho_u'^k\|^2
\nonumber
\\
& 
\quad 
+(4c_1\hspace{-0.07cm}+\hspace{-0.07cm}L_\mu)\tau \sum_{k=0}^{m-1} \|\nabla\rho_u^k\|^2 
+ \sum_{k=1}^m 2\tau \big[ \mathcal{F}_u^k 
+ \mathcal{\overline{F}}_\theta^k
+ \mathcal{F}_\delta^k\big]
+ \sum_{k=1}^m 2\tau \big[\mathcal{\overline{R}}_u^k 
+ \mathcal{\overline{R}}_{\theta}^k
+ \mathcal{\overline{R}}_{\delta}^k
\big].\label{soma2}
\end{align}
Now, in \eqref{soma2} defining the constants $c_2\hspace{-0.07cm}=\hspace{-0.07cm}\min\big\{1\hspace{-0.05cm} - \hspace{-0.05cm}c_1\tau, \frac{\mu_0}{2}\hspace{-0.05cm}- \hspace{-0.05cm}\tau c_1, \mu_0\underline{f}_1\hspace{-0.05cm}-\hspace{-0.05cm}c_1\tau, \mu_0\underline{f}_3\hspace{-0.05cm}-\hspace{-0.05cm}c_1\tau \big\}>0$ and $c_3=\max\big\{1, 2\mu_1, \mu_1\overline{f}_1, \mu_1\overline{f}_3, 2c_1, 2c_1\hspace{-0.07cm}+\hspace{-0.07cm}L_\mu\overline{f}_1, 2c_1\hspace{-0.07cm}+\hspace{-0.07cm}L_\mu\overline{f}_3, 4c_1\hspace{-0.07cm}+\hspace{-0.07cm}L_\mu, 2\big\}/c_2>0$, and recalling that $\mathcal{E}^k := \|\rho_u'^k\|^2 + \|\rho_\theta^k\|^2 + \|\nabla\rho_u^k\|^2 + \|\rho_\delta'^k\|_{\Gamma_1}^2 + \|\rho_\delta^k\|_{\Gamma_1}^2$, see \eqref{energia}, follows
\begin{align}
& \mathcal{E}^m 
\leqslant 
c_3 \mathcal{E}^0 
+ c_3\tau \sum_{k=0}^{m-1}\mathcal{E}^k
+ c_3\tau \sum_{k=1}^m 
\big[
 \mathcal{F}_u^k 
+ \mathcal{\overline{F}}_\theta^k
+ \mathcal{F}_\delta^k
+ \mathcal{\overline{R}}_u^k 
+ \mathcal{\overline{R}}_{\theta}^k
+ \mathcal{\overline{R}}_{\delta}^k
\big].
\label{soma3}
\end{align}
In \eqref{soma3},  applying the discrete Gronwall lemma, see \cite{Artigo1960_Lees}, we obtain
\begin{align}
\max_{0\leqslant m\leqslant N}\mathcal{E}^m 
\leqslant 
c_4\Big[\mathcal{E}^0 
+ \tau \sum_{k=1}^N                 
\big(
 \mathcal{F}_u^k 
+ \mathcal{\overline{F}}_\theta^k
+ \mathcal{F}_\delta^k
+ \mathcal{\overline{R}}_u^k 
+ \mathcal{\overline{R}}_{\theta}^k
+ \mathcal{\overline{R}}_{\delta}^k
\big)\Big],
\label{Gronwaal}
\end{align}
with $c_4 = c_3e^{c_3T}$, and the terms $\mathcal{F}_u^k, \mathcal{\overline{F}}_\theta^k , \mathcal{F}_\delta^k$ and 
$ \mathcal{\overline{R}}_u^k, \mathcal{\overline{R}}_{\theta}^k, \mathcal{\overline{R}}_{\delta}^k$ are defined in \{\eqref{f_u}-\eqref{f_delta} and \eqref{f_theta1}\}, and  \{\eqref{resto_u}-\eqref{resto_delta}\} , respectively. To simplify, from now on, $c>0$ will denote various constants throughout the estimates.

\noindent{\it Final bounds}.
In \eqref{Gronwaal}, using the projection estimates (Lemma \ref{lem:proj}) for the terms $\mathcal{F}_u^k, \mathcal{\overline{F}}_\theta^k, \mathcal{F}_\delta^k$, we obtain:

\noindent\textbf{Estimate of $c_4 \tau\sum_{k=1}^N \mathcal{F}_u^k$}. 

\noindent\textbf{Term $\frac{c_4}{4\epsilon_1}\sum_{k=1}^N \tau \|\widetilde{\partial}\omega_u'^k\|^2$}. The finite difference $\widetilde{\partial}\omega_u'^k$ can be expressed as:
$
\widetilde{\partial}\omega_u'^k = \frac{\omega_u'^k - \omega_u'^{k-1}}{\tau}
= \frac{1}{\tau} \int_{t_{k-1}}^{t_k} \omega_u''(t) \, dt.
$
Now, applying the Cauchy–Schwarz inequality,
$
\|\widetilde{\partial}\omega_u'^k\|^2
\leqslant \frac{1}{\tau} \int_{t_{k-1}}^{t_k} \|\omega_u''(t)\|^2 dt.
$
So,
\begin{align}
\hspace{-0.08cm}\frac{c_4}{4\epsilon_1}\hspace{-0.10cm}\sum_{k=1}^N \hspace{-0.08cm}\tau \|\widetilde{\partial}\omega_u'^k\|^2
\hspace{-0.08cm}
\leqslant \hspace{-0.08cm}
\frac{c_4}{4\epsilon_1} \sum_{k=1}^N \hspace{-0.08cm}\int_{t_{k-1}}^{t_k} \hspace{-0.20cm}\|\omega_u''(t)\|^2 dt
\hspace{-0.08cm}
&=\hspace{-0.08cm}\frac{c_4}{4\epsilon_1}  \int_{0}^{t_N}\hspace{-0.20cm} \|\omega_u''(t)\|^2 dt
\nonumber
\\
\hspace{-0.08cm} &=\hspace{-0.08cm}\frac{c_4}{4\epsilon_1}\|\omega_u''\|^2_{L^2(L^2(\Omega))}\hspace{-0.08cm}\leqslant\hspace{-0.08cm}
c_4c\frac{1}{4\epsilon_1} h^{2(r+1)}\|u''\|^2_{L^2(H^{r+1}(\Omega))}.\label{est-f_u1}
\end{align}

\noindent\textbf{Term $c_4\frac{3c_f^2}{4\epsilon_1} \sum_{k=1}^N \tau \|\widetilde{\omega}_u^{\,k}\|^2$}.
\begin{align}
c_4\frac{3c_f^2}{4\epsilon_1} \sum_{k=1}^N \tau \|\widetilde{\omega}_u^{\,k}\|^2
\leqslant
c_4\frac{3c_f^2}{4\epsilon_1} \sum_{k=1}^N \tau \|\omega_u^{\,k}\|^2
&\leqslant
c_4\frac{3c_f^2}{4\epsilon_1}c\|\omega_u\|^2_{L^2(L^2(\Omega))}\leqslant
c_4c\frac{3c_f^2}{4\epsilon_1} h^{2(r+1)}\|u\|^2_{L^2(H^{r+1}(\Omega))}.\label{est-f_u2}
\end{align}

\noindent\textbf{Term $c_4\frac{\|\mathbf{a}\|_\infty^2}{\epsilon_{4}} \sum_{k=1}^N \tau \|\widetilde{\omega}_u'^{\,k}\|^2$}. Similarly to \eqref{est-f_u2}, 
\begin{align}
\hspace{-0.08cm}c_4\frac{\|\mathbf{a}\|_\infty^2}{\epsilon_{4}} \hspace{-0.08cm}\sum_{k=1}^N \tau \|\widetilde{\omega}_u'^{\,k}\|^2
&\hspace{-0.08cm}\leqslant\hspace{-0.08cm}
c_4\frac{\|\mathbf{a}\|_\infty^2}{\epsilon_{4}} \hspace{-0.08cm}\sum_{k=1}^N \tau \|\omega_u'^{\,k}\|^2
\hspace{-0.08cm}\leqslant\hspace{-0.08cm}
c_4c\frac{\|\mathbf{a}\|_\infty^2}{\epsilon_{4}}\|\omega_u'\|^2_{L^2(L^2(\Omega))}
\hspace{-0.08cm}\leqslant\hspace{-0.08cm}
c_4c\frac{\|\mathbf{a}\|_\infty^2}{\epsilon_{4}} h^{2(r+1)}\|u'\|^2_{L^2(H^{r+1}(\Omega))}.\label{est-f_u3}
\end{align}

\noindent\textbf{Term $c_4\Big(\frac{L^2_g\mu_1^{2}}{2\epsilon_2} \hspace{-0.07cm}+\hspace{-0.07cm} \frac{\mu_1^2}{4\epsilon_6}\Big) \sum_{k=1}^N \tau 
\|\widetilde{\omega}_u'^{\,k}\|_{\Gamma_1}^2$}. Completely analogous to \eqref{est-f_u3},
\begin{align}
\hspace{-0.08cm}c_4\Big(\hspace{-0.08cm}\frac{L^2_g\mu_1^{2}}{2\epsilon_2} \hspace{-0.09cm}+\hspace{-0.09cm}\frac{\mu_1^2}{4\epsilon_6}\Big) \hspace{-0.08cm}\sum_{k=1}^N\hspace{-0.08cm} \tau \|\widetilde{\omega}_u'^{\,k}\|_{\Gamma_1}^2
&\hspace{-0.09cm}\leqslant\hspace{-0.09cm}
c_4c\Big(\hspace{-0.08cm}\frac{L^2_g\mu_1^{2}}{2\epsilon_2} \hspace{-0.09cm}+\hspace{-0.09cm} \frac{\mu_1^2}{4\epsilon_6}\hspace{-0.08cm}\Big)\|\omega_u'\|^2_{L^2(L^2(\Gamma_1))}
 \hspace{-0.08cm}\leqslant\hspace{-0.08cm}
c_4c\Big(\frac{L^2_g\mu_1^{2}}{2\epsilon_2} \hspace{-0.09cm}+\hspace{-0.09cm} \frac{\mu_1^2}{4\epsilon_6}\Big) h^{2(r+\frac{1}{2})}\|u'\|^2_{L^2(H^{r+1}(\Omega))}.\label{est-f_u4}
\end{align}

We emphasize that the final convergence order is suboptimal due to this domain–boundary relation in \eqref{est-f_u4}.

\noindent\textbf{Estimate of $c_4\tau \sum_{k=1}^N \mathcal{\overline{F}}_\theta^k$}.

\noindent\textbf{Term $c_4\frac{\|\mathbf{a}\|_{\infty}^2}{4}\tau \sum_{k=1}^N \|\widetilde{\partial}\widetilde{\omega}^k_\theta\|^2$}. Analogously to what was done in \eqref{est-f_u1} and 
\eqref{est-f_u2}, it follows that
\begin{align}
 c_4\frac{\|\mathbf{a}\|_{\infty}^2}{4}\hspace{-0.07cm}\sum_{k=1}^N \tau \|\widetilde{\partial}\widetilde{\omega}^k_\theta\|^2
\hspace{-0.07cm}
\leqslant \hspace{-0.07cm}
c_4c\frac{\|\mathbf{a}\|_{\infty}^2}{4}\|\omega_\theta'\|^2_{L^2(L^2(\Omega))}
\hspace{-0.07cm}\leqslant\hspace{-0.07cm}
c_4c\frac{\|\mathbf{a}\|_{\infty}^2}{4}h^{2(r+1)}\|\theta'\|^2_{L^2( H^{r+1}(\Omega))}.
\label{est-f_theta1}
\end{align}

\noindent\textbf{Term $c_4\frac{1}{4\epsilon_3}\tau \sum_{k=1}^N \|\widetilde{\partial}\omega_\theta^k\|^2$}. Proceeding analogously to \eqref{est-f_u1} and \eqref{est-f_u2}, it follows that
\begin{align}
 c_4\frac{1}{4\epsilon_3}\hspace{-0.07cm}\sum_{k=1}^N \tau \|\widetilde{\partial}\omega_\theta^k\|^2
\hspace{-0.07cm}
\leqslant \hspace{-0.07cm}
c_4c\frac{1}{4\epsilon_3}\|\omega_\theta'\|^2_{L^2(L^2(\Omega))}
\leqslant
c_4c\frac{1}{4\epsilon_3} h^{2(r+1)}\|\theta'\|^2_{L^2(H^{r+1}(\Omega))}.
\label{est-f_theta2}
\end{align}

\noindent\textbf{Term $c_4\Big(\frac{3c_\beta^2}{4\epsilon_4}\hspace{-0.07cm}+\hspace{-0.07cm}1 \Big) \sum_{k=1}^N \tau \|\widetilde{\omega}_\theta^{\,k}\|^2$}. Similarly to \eqref{est-f_u2}, it follows that 
\begin{align}
c_4\Big(\frac{3c_\beta^2}{4\epsilon_4}\hspace{-0.07cm}+\hspace{-0.07cm}1 \Big) \sum_{k=1}^N \tau \|\widetilde{\omega}_\theta^{\,k}\|^2
&\leqslant
c_4c\Big(\frac{3c_\beta^2}{4\epsilon_4}\hspace{-0.07cm}+\hspace{-0.07cm}1 \Big)\|\omega_\theta\|^2_{L^2(L^2(\Omega))}
\leqslant
c_4c\Big(\frac{3c_\beta^2}{4\epsilon_4}\hspace{-0.07cm}+\hspace{-0.07cm}1 \Big)h^{2(r+1)}\|\theta\|^2_{L^2(H^{r+1}(\Omega))}.\label{est-f_theta2+}
\end{align}

\noindent\textbf{Term $c_4\frac{2\|\mathbf{a}\|^2_\infty}{\mu_0} \sum_{k=1}^N \tau \|\widetilde{\omega}_\theta^m\|^2$}.  
We observe that the summation does not depend on the index $k$, since $\widetilde{\omega}_\theta^m$ is constant with respect to $k$. By the projection estimate (Lemma \ref{lem:proj}), we have $\|\widetilde{\omega}_\theta^m\|\leqslant c h^{r+1} \|\theta(t_m)\|_{H^{r+1}(\Omega)}.$ Hence,
\begin{align}
& c_4\frac{2\|\mathbf{a}\|^2_\infty}{\mu_0}\sum_{k=1}^N \tau \|\widetilde{\omega}_\theta^m\|^2 
= c_4\frac{2\|\mathbf{a}\|^2_\infty}{\mu_0}\|\widetilde{\omega}_\theta^m\|^2 \sum_{k=1}^N \tau 
=c_4\frac{2\|\mathbf{a}\|^2_\infty}{\mu_0}T\|\widetilde{\omega}_\theta^m\|^2\nonumber
\\
&
\qquad\leqslant
Tc_4c\frac{2\|\mathbf{a}\|^2_\infty}{\mu_0} h^{2(r+1)} \|\theta(t_m)\|^2_{H^{r+1}(\Omega))}
\leqslant
Tc_4c\frac{2\|\mathbf{a}\|^2_\infty}{\mu_0} h^{2(r+1)} \|\theta\|^2_{C(0,T; H^{r+1}(\Omega))}.
\label{est-f_theta3}
\end{align}

\noindent\textbf{Term $c_4\frac{\|\mathbf{a}\|^2_\infty}{\mu_1} \sum_{k=1}^N \tau \|\widetilde{\omega}_\theta^0\|^2$}. 
Since the Crank–Nicolson scheme does not define $\widetilde{\omega}^{\,0}$ from a mean (because there is no $\omega^{-1}$), we adopt the natural convention $\widetilde{\omega}^{\,0} = \omega(0)$. Consequently,
$\widetilde{\omega}_\theta^0 = R_h^\theta\theta(0) - \theta(0) = \omega_\theta(0)=\omega_\theta^0$. Thus, following analogously to \eqref{est-f_theta3},
\begin{align} 
\hspace{-0.08cm}c_4\frac{\|\mathbf{a}\|^2_\infty}{\mu_1}
\hspace{-0.09cm}\sum_{k=1}^N \hspace{-0.09cm}\tau \|\widetilde{\omega}_\theta^0\|^2
&\hspace{-0.08cm}=\hspace{-0.08cm} c_4\frac{\|\mathbf{a}\|^2_\infty}{\mu_1}\hspace{-0.08cm} T \|\omega_\theta^0\|^2
\hspace{-0.08cm}\leqslant\hspace{-0.08cm}
Tc_4c\frac{\|\mathbf{a}\|^2_\infty}{\mu_1}h^{2(r+1)}\hspace{-0.03cm}\|\theta(0)\|_{H^{r+1}(\Omega)}
\hspace{-0.08cm}\leqslant\hspace{-0.08cm}
Tc_4c\frac{\|\mathbf{a}\|^2_\infty}{\mu_1} h^{2(r+1)}\hspace{-0.03cm}\|\theta\|^2_{\hspace{-0.04cm}C(H^{r+1}(\Omega))}.\label{est-f_theta4} 
\end{align}
From \eqref{theta-regularidade}, $\theta \in H^1(H^{r+1}(\Omega)\cap H^1_0(\Omega)) \subset H^1(H^{r+1}(\Omega))\hookrightarrow C([0,T]; H^{r+1}(\Omega)),$ justifying \eqref{est-f_theta3} and \eqref{est-f_theta4}.

\noindent\textbf{Estimate of $c_4 \tau\sum_{k=1}^N \mathcal{F}_\delta^k$}. 

\noindent\textbf{Term $c_4\frac{\mu_1^2\overline{f}_1}{4\epsilon_6} \sum_{k=1}^N \tau \|\widetilde{\partial}\omega_\delta'^k\|_{\Gamma_1}^2 $}. Analogously to \eqref{est-f_u1}, 
\begin{align}
 c_4\frac{\mu_1^2\overline{f}_1}{4\epsilon_6}\sum_{k=1}^N \tau \|\widetilde{\partial}\omega_{\delta}'^k\|_{\Gamma_1}^2
\leqslant
c_4\frac{\mu_1^2\overline{f}_1}{4\epsilon_6}\|\omega_{\delta}''\|^2_{L^2(L^2(\Gamma_1))}\leqslant
c_4c\frac{\mu_1^2\overline{f}_1}{4\epsilon_6} h^{2(r+1)}\|\delta''\|^2_{L^2(H^{r+1}(\Gamma_1))}.\label{est-f_delta1}
\end{align}

\noindent\textbf{Term $c_4 \Big(\frac{\mu_1^{2}}{4\epsilon_2}
+ \frac{\mu_1^2\overline{f}_2}{4\epsilon_6}\Big)\sum_{k=1}^N \tau \|\widetilde{\omega}_\delta'^{\,k}\|_{\Gamma_1}^2$}. Similar to \eqref{est-f_u3},
\begin{align}
c_4 \Big(\hspace{-0.08cm}\frac{\mu_1^{2}}{4\epsilon_2}
\hspace{-0.09cm} + \hspace{-0.09cm} \frac{\mu_1^2\overline{f}_2}{4\epsilon_6}\hspace{-0.08cm}\Big) \hspace{-0.08cm}\sum_{k=1}^N \hspace{-0.08cm}\tau \|\widetilde{\omega}_\delta'^{\,k}\|_{\Gamma_1}^2
\hspace{-0.08cm}\leqslant\hspace{-0.08cm}
c_4 c\Big(\hspace{-0.08cm}\frac{\mu_1^{2}}{4\epsilon_2}
\hspace{-0.07cm}+\hspace{-0.07cm} \frac{\mu_1^2\overline{f}_2}{4\epsilon_6}\hspace{-0.08cm}\Big)\|\omega_{\delta}'\|^2_{L^2(L^2(\Gamma_1))}
\hspace{-0.08cm} \leqslant\hspace{-0.08cm}
c_4 c\Big(\hspace{-0.08cm}\frac{\mu_1^{2}}{4\epsilon_2}
\hspace{-0.07cm}+\hspace{-0.07cm} \frac{\mu_1^2\overline{f}_2}{4\epsilon_6}\hspace{-0.08cm}\Big) h^{2(r+1)}\|\delta'\|^2_{L^2(H^{r+1}(\Gamma_1))}.\label{est-f_delta2}
\end{align}

\noindent\textbf{Term $c_4 \frac{\mu_1^2\overline{f}_3}{4\epsilon_6}\sum_{k=1}^N \tau \|\widetilde{\omega}_\delta^{\,k}\|_{\Gamma_1}^2$}. In a manner analogous to \eqref{est-f_u2},
\begin{align}
c_4 \frac{\mu_1^2\overline{f}_3}{4\epsilon_6} \sum_{k=1}^N \tau \|\widetilde{\omega}_\delta^{\,k}\|_{\Gamma_1}^2
\leqslant
c_4c \frac{\mu_1^2\overline{f}_3}{4\epsilon_6}\|\omega_\delta\|^2_{L^2(L^2(\Gamma_1))}
\leqslant
c_4c \frac{\mu_1^2\overline{f}_3}{4\epsilon_6} h^{2(r+1)}\|\delta\|^2_{L^2(H^{r+1}(\Gamma_1))}.\label{est-final-f}
\end{align}

\noindent\textbf{Estimate of the terms relating to $\mathcal{\overline{R}}_u^k, \mathcal{\overline{R}}_{\theta}^k, \mathcal{\overline{R}}_{\delta}^k$}. 
Based on the Taylor expansions and on results similar to those found in \cite{Artigo2004_ErvinHeuer}, \cite{Artigo2025_Adriano}, and \cite{Livro2009_Kendall_Atkinson}, all the following estimates hold.

\noindent\textbf{Estimate of $c_4 \tau\sum_{k=1}^N \mathcal{\overline{R}}_u^k$}. 
\begin{align}
& \triangleright\; 
c_4 \frac{3c_f^2}{4\epsilon_1}\sum_{k=1}^N \tau \|\widetilde{u}^k - u(t_{k-\frac{1}{2}})\|^2
\leqslant 
c_4c\frac{3c_f^2}{4\epsilon_1}\tau^4\|u''\|^2_{L^2(L^2(\Omega))};\label{est-inicio-resto}
\\[-6pt]
& \triangleright\; 
c_4\frac{1}{4\epsilon_1}\sum_{k=1}^N \tau\|\widetilde{\partial}u'^{k} - u''(t_{k-\frac{1}{2}})\|^2 
\leqslant
c_4c\frac{1}{4\epsilon_1}\tau^4\|u^4\|^2_{L^2(L^2(\Omega))};
\\[-5pt]
& \triangleright\; 
c_4\Big(\frac{\mu_1^{2}}{4\epsilon_0}\hspace{-0.08cm}+\hspace{-0.08cm}\frac{\|\mathbf{a}\|_{\infty}^2}{4}\Big)\hspace{-0.08cm}\sum_{k=1}^N \tau\|\nabla(\widetilde{\partial}u^k \hspace{-0.08cm}-\hspace{-0.08cm}\widetilde{u}'^k)\|^2 
\hspace{-0.08cm}\leqslant\hspace{-0.08cm}
c_4c\Big(\frac{\mu_1^{2}}{4\epsilon_0}\hspace{-0.08cm}+\hspace{-0.08cm}\frac{\|\mathbf{a}\|_{\infty}^2}{4}\Big)\hspace{-0.08cm}\;
\tau^4\|\nabla u'''\|^2_{L^2(L^2(\Omega))}; 
\\[-5pt]
& \triangleright\; 
c_4\frac{\mu_1^2}{4\epsilon_1}\sum_{k=1}^N \tau\|\Delta (\widetilde{u}^k-u(t_{k-\frac{1}{2}}))\|^2 
\leqslant
c_4c\frac{\mu_1^2}{4\epsilon_1} \tau^4\|\Delta u''\|^2_{L^2(L^2(\Omega))};
\\[-5pt]
& \triangleright\; 
c_4\frac{\mu_1^2}{4\epsilon_2}\sum_{k=1}^N \tau\|\frac{\partial}{\partial \nu}(\widetilde{u}^k-u(t_{k-\frac{1}{2}}))\|^2_{\Gamma_1} 
\leqslant
c_4c\frac{\mu_1^2}{4\epsilon_2} \tau^4\|\frac{\partial}{\partial \nu} u''\|^2_{L^2(L^2(\Gamma_1))};
\\[-5pt]
& \triangleright\; 
c_4\frac{\|a\|^2_{\infty}}{4\epsilon_3}\sum_{k=1}^N \tau\|\nabla(\widetilde{u}'^k - u'(t_{k-\frac{1}{2}}))\|^2
\leqslant
c_4c\frac{\|a\|^2_{\infty}}{4\epsilon_3}\; \tau^4\|\nabla u'''\|^2_{L^2(L^2(\Omega))}; 
\\[-5pt]
& \triangleright\; 
c_4 \Big(\frac{\mu_1^2}{4\epsilon_6} 
\hspace{-0.08cm}+\hspace{-0.08cm} \frac{L^2_g\mu_1^{2}}{2\epsilon_2}\Big)\sum_{k=1}^N \tau\|\widetilde{u}'^k \hspace{-0.08cm}-\hspace{-0.08cm} u'(t_{k-\frac{1}{2}})\|_{\Gamma_1}^2
\leqslant
c_4c \Big(\frac{\mu_1^2}{4\epsilon_6} 
\hspace{-0.08cm}+\hspace{-0.08cm} \frac{L^2_g\mu_1^{2}}{2\epsilon_2}\Big) \tau^4\| u'''\|^2_{L^2(L^2(\Gamma_1))};
\end{align}

\noindent\textbf{Estimate of $c_4 \tau\sum_{k=1}^N \mathcal{\overline{R}}_{\theta}^k$}. 
\begin{align}
& \triangleright\; 
c_4 \frac{1}{4\epsilon_3}\sum_{k=1}^N \tau \|\widetilde{\partial}\theta^{k} - \theta'(t_{k-\frac{1}{2}})\|^2
\leqslant 
c_4c\frac{1}{4\epsilon_3}\tau^4\|\theta'''\|^2_{L^2(L^2(\Omega))};
\\[-5pt]
& \triangleright\; 
c_4\frac{3c_\beta^2}{4\epsilon_4}\sum_{k=1}^N \tau\|\widetilde{\theta}^k-\theta(t_{k-\frac{1}{2}})\|^2
\leqslant
c_4c\frac{3c_\beta^2}{4\epsilon_4}\tau^4\|\theta''\|^2_{L^2(L^2(\Omega))};
\\[-5pt]
& \triangleright\; 
c_4\Big(\hspace{-0.08cm}\frac{\|a\|^2_{\infty}}{4\epsilon_1}\hspace{-0.08cm}+\hspace{-0.08cm} \frac{\overline{c}_\beta^2}{4\epsilon_4}\hspace{-0.08cm}\Big)\hspace{-0.08cm}\sum_{k=1}^N \tau\|\nabla(\widetilde{\theta}^k \hspace{-0.08cm}-\hspace{-0.08cm}\theta(t_{k-\frac{1}{2}}))\|^2
\hspace{-0.08cm}\leqslant\hspace{-0.08cm}
c_4c\Big(\frac{\|a\|^2_{\infty}}{4\epsilon_1} \hspace{-0.08cm}+\hspace{-0.08cm} \frac{\overline{c}_\beta^2}{4\epsilon_4}\Big) \tau^4\|\nabla \theta''\|^2_{L^2(L^2(\Omega))}; 
\end{align}

\noindent\textbf{Estimate of $c_4 \tau\sum_{k=1}^N \mathcal{\overline{R}}_{\delta}^k$}.
\begin{align}
& \triangleright\; 
c_4 \frac{\overline{f}_3\mu_1^2}{4\epsilon_6}\sum_{k=1}^N \tau \|\widetilde{\delta}^k - \delta(t_{k-\frac{1}{2}})\|_{\Gamma_1}^2
\leqslant 
c_4c\frac{\overline{f}_3\mu_1^2}{4\epsilon_6}\tau^4\|\delta''\|^2_{L^2(L^2(\Gamma_1))};
\\[-5pt]
& \triangleright\; 
c_4\frac{\mu_1^2\overline{f}_3}{4\epsilon_5}\sum_{k=1}^N \tau\|\widetilde{\partial}\delta^k - \widetilde{\delta}'^k\|_{\Gamma_1}^2
\leqslant
c_4c\frac{\mu_1^2\overline{f}_3}{4\epsilon_5}\tau^4\|\delta'''\|^2_{L^2(L^2(\Gamma_1))};
\\[-5pt]
& \triangleright\; 
c_4\frac{\overline{f}_1\mu_1^2}{4\epsilon_6}\sum_{k=1}^N \|\widetilde{\partial}\delta'^{k} - \delta''(t_{k-\frac{1}{2}})\|_{\Gamma_1}^2
\leqslant         
c_4c\frac{\overline{f}_1\mu_1^2}{4\epsilon_6}\tau^4\|\delta^4\|^2_{L^2(L^2(\Omega))}; 
\\[-5pt]
& \triangleright\; 
c_4\Big(\frac{\overline{f}_2\mu_1^2}{4\epsilon_6}\hspace{-0.08cm}+\hspace{-0.08cm}\frac{\mu_1^2}{4\epsilon_2}\Big)\hspace{-0.08cm}\sum_{k=1}^N \tau\|\widetilde{\delta}'^k \hspace{-0.08cm}-\hspace{-0.08cm} \delta'(t_{k-\frac{1}{2}})\|_{\Gamma_1}^2
\hspace{-0.08cm}\leqslant\hspace{-0.08cm}
c_4c\Big(\frac{\overline{f}_2\mu_1^2}{4\epsilon_6}\hspace{-0.08cm}+\hspace{-0.08cm}\frac{\mu_1^2}{4\epsilon_2}\Big)\tau^4\|\delta'''\|^2_{L^2(L^2(\Gamma_1))}.\label{est-final-delta}
\end{align}
Now, returning to \eqref{Gronwaal}, using the estimates \eqref{est-f_u1}–\eqref{est-final-f} and \eqref{est-inicio-resto}–\eqref{est-final-delta}, and defining the constant $c_5 > 0$ such that:
\begin{align}
\hspace{-0.08cm}c_5&\hspace{-0.07cm}=\hspace{-0.07cm}\max\hspace{-0.07cm}\Big\{\hspace{-0.08cm}c_4, c_4c\frac{1}{4\epsilon_1}, c_4c\frac{3c_f^2}{4\epsilon_1}, c_4c\frac{\|\mathbf{a}\|_\infty^2}{\epsilon_{4}}, c_4c\hspace{-0.05cm}\Big(\hspace{-0.08cm}\frac{L^2_g\mu_1^{2}}{2\epsilon_2} \hspace{-0.08cm}+\hspace{-0.08cm} \frac{\mu_1^2}{4\epsilon_6}\hspace{-0.08cm}\Big), c_4c\frac{\|\mathbf{a}\|_{\infty}^2}{4}, c_4c\frac{1}{4\epsilon_3},
c_4c\hspace{-0.05cm}\Big(\hspace{-0.08cm}\frac{3c_\beta^2}{4\epsilon_4}\hspace{-0.08cm}+\hspace{-0.08cm}1 \hspace{-0.08cm}\Big),
 T\hspace{-0.04cm}c_4c\frac{2\|\mathbf{a}\|^2_\infty}{\mu_0}, c_4c\hspace{-0.05cm}\frac{1}{4\epsilon_1},
 \nonumber
\\
&
\;
Tc_4c\hspace{-0.06cm}\frac{\|\mathbf{a}\|^2_\infty}{\mu_1},
c_4c\frac{\mu_1^2\overline{f}_1}{4\epsilon_6},
   c_4 c\Big(\hspace{-0.08cm}\frac{\mu_1^{2}}{4\epsilon_2}
\hspace{-0.08cm}+\hspace{-0.08cm} \frac{\mu_1^2\overline{f}_2}{4\epsilon_6}\hspace{-0.08cm}\Big), c_4c \frac{\mu_1^2\overline{f}_3}{4\epsilon_6},  
c_4c\frac{3c_f^2}{4\epsilon_1}, c_4c\Big(\hspace{-0.08cm}\frac{\mu_1^{2}}{4\epsilon_0}\hspace{-0.08cm}+\hspace{-0.08cm}\frac{\|\mathbf{a}\|_{\infty}^2}{4}\hspace{-0.08cm}\Big), c_4c\frac{\mu_1^2}{4\epsilon_1}, c_4c\frac{\mu_1^2}{4\epsilon_2}, c_4c\frac{\|a\|^2_{\infty}}{4\epsilon_3},     \nonumber
\\
&
\;
c_4c\frac{1}{4\epsilon_3},c_4c\frac{3c_\beta^2}{4\epsilon_4},
c_4c \Big(\hspace{-0.08cm}\frac{\mu_1^2}{4\epsilon_6} 
\hspace{-0.08cm}+\hspace{-0.08cm} \frac{L^2_g\mu_1^{2}}{2\epsilon_2}\hspace{-0.08cm}\Big),  c_4c\Big(\hspace{-0.08cm}\frac{\|a\|^2_{\infty}}{4\epsilon_1} \hspace{-0.09cm}+\hspace{-0.09cm} \frac{\overline{c}_\beta^2}{4\epsilon_4}\hspace{-0.08cm}\Big), 
c_4c\frac{\overline{f}_3\mu_1^2}{4\epsilon_6}, c_4c\frac{\mu_1^2\overline{f}_3}{4\epsilon_5}, c_4c\frac{\overline{f}_1\mu_1^2}{4\epsilon_6},  c_4c\Big(\hspace{-0.08cm}\frac{\overline{f}_2\mu_1^2}{4\epsilon_6}\hspace{-0.09cm}+\hspace{-0.09cm}\frac{\mu_1^2}{4\epsilon_2}\hspace{-0.08cm}\Big)
\hspace{-0.08cm}\Big\},\label{constante-geral1}
\end{align}
thus,
\begin{align}
&\max_{0\leqslant m\leqslant N}\hspace{-0.09cm}\mathcal{E}^m 
\leqslant 
c_5\Big[\mathcal{E}^0         
\hspace{-0.09cm}+\hspace{-0.09cm} h^{2(r+\frac{1}{2})}\Big(\hspace{-0.02cm}\|u''\|^2_{L^2(H^{r+1}(\Omega))} 
\hspace{-0.09cm}+\hspace{-0.09cm} \|u\|^2_{L^2(H^{r+1}(\Omega))}
\hspace{-0.09cm}+\hspace{-0.09cm} 2\|u'\|^2_{L^2(H^{r+1}(\Omega))}
\hspace{-0.09cm}+\hspace{-0.09cm} 2\|\theta'\|^2_{L^2(H^{r+1}(\Omega))} 
\nonumber
\\
&
\quad
\hspace{-0.09cm}+\hspace{-0.09cm} \|\theta\|^2_{L^2(H^{r+1}(\Omega))}
\hspace{-0.09cm}+\hspace{-0.09cm} 2\|\theta\|^2_{C(0,T; H^{r+1}(\Omega))}
\hspace{-0.09cm}+\hspace{-0.09cm} \|\delta''\|^2_{L^2(H^{r+1}(\Gamma_1))}
\hspace{-0.09cm}+\hspace{-0.09cm} \|\delta'\|^2_{L^2(H^{r+1}(\Gamma_1))}
\hspace{-0.09cm}+\hspace{-0.09cm} \|\delta\|^2_{L^2(H^{r+1}(\Gamma_1))}
\hspace{-0.02cm}\Big)
\nonumber
\\
&
\quad
\hspace{-0.09cm}+\hspace{-0.09cm} \tau^4\Big(\hspace{-0.09cm}\|u''\|^2_{L^2(L^2(\Omega))}
\hspace{-0.09cm}+\hspace{-0.09cm} \|u^4\|^2_{L^2(L^2(\Omega))}
\hspace{-0.09cm}+\hspace{-0.09cm} 2\|\nabla u'''\|^2_{L^2(L^2(\Omega))}
\hspace{-0.09cm}+\hspace{-0.09cm} \|\Delta u''\|^2_{L^2(L^2(\Omega))}
\hspace{-0.09cm}+\hspace{-0.09cm} \|\frac{\partial}{\partial \nu} u''\|^2_{L^2(L^2(\Gamma_1))}
\hspace{-0.09cm}+\hspace{-0.09cm} \| u'''\|^2_{L^2(L^2(\Gamma_1))}\nonumber
\\
&
\quad
\hspace{-0.09cm}+\hspace{-0.09cm} \|\theta'''\|^2_{L^2(L^2(\Omega))}
\hspace{-0.09cm}+\hspace{-0.09cm} \|\theta''\|^2_{L^2(L^2(\Omega))}
\hspace{-0.09cm}+\hspace{-0.09cm} \|\nabla \theta''\|^2_{L^2(L^2(\Omega))}
\hspace{-0.09cm}+\hspace{-0.09cm} \|\delta''\|^2_{L^2(L^2(\Gamma_1))}
\hspace{-0.09cm}+\hspace{-0.09cm} 2\|\delta'''\|^2_{L^2(L^2(\Gamma_1))}
\hspace{-0.09cm}+\hspace{-0.09cm} \|\delta^4\|^2_{L^2(L^2(\Omega))}
\hspace{-0.09cm}\Big)
\Big],\label{Gronwaal1}
\end{align}
since $0<h<1$, it follows that $h^{2(r+1)}<h^{2(r+\frac{1}{2})}$, that is, $h^{2(r+\frac{1}{2})}$ dominates.

\noindent\textit{Transition to the estimate in linear norms}. From inequality  \eqref{Gronwaal1} for the term $\mathcal{E}^m$ (see \eqref{energia}), we have
\begin{align}
\max_{0\leqslant m\leqslant N}\mathcal{E}^m \leqslant c_5\Big( \mathcal{E}^0 + h^{2(r+\frac{1}{2})}\mathcal{N}_1 + \tau^4\mathcal{N}_2 \Big),\label{energia-01}
\end{align}
where $\mathcal{N}_1$ and $\mathcal{N}_2$ are the combinations of norms of the exact solutions and their derivatives appearing on the right-hand side of \eqref{Gronwaal1}, and the constant $c_5>0$ is the one defined in \eqref{constante-geral1}. Furthermore, from the initial condition and the properties of the Ritz projection, we know that $\mathcal{E}^0 \leqslant c_0 h^{2(r+\frac{1}{2})}$, with $c_0>0$ independent of $h$ and $\tau$; see this hypothesis in Theorem \ref{teo:convergencia}. To obtain an estimate for the sum of the linear norms, we apply the Cauchy–Schwarz inequality in the space $\mathbb{R}^5$. We define the vectors
$
\mathbf{x_1} = (1,1,1,1,1)$ and
$\mathbf{x_2} = \big(\|\rho_u'^m\|,\; \|\rho_\theta^m\|,\; \|\nabla\rho_u^m\|,\; \|\rho_\delta'^m\|_{\Gamma_1},\; \|\rho_\delta^m\|_{\Gamma_1}\big).
$
The Cauchy-Schwarz inequality $\langle \mathbf{x}_1, \mathbf{x}_2 \rangle \leqslant \|\mathbf{x}_1\|_2 \|\mathbf{x}_2\|_2$ (where $ \|\cdot\|_2 $ denotes the Euclidean norm induced by the inner product $ \langle \cdot, \cdot \rangle $) yields
\begin{align}
\|\rho_u'^m\| + \|\rho_\theta^m\| + \|\nabla\rho_u^m\| + \|\rho_\delta'^m\|_{\Gamma_1} + \|\rho_\delta^m\|_{\Gamma_1}
\leqslant \sqrt{5}\,\sqrt{\mathcal{E}^m}.\label{energia-02}
\end{align}
Therefore, for all $m$, from \eqref{energia-02} in \eqref{energia-01},
\begin{align*}
\|\rho_u'^m\| + \|\rho_\theta^m\| + \|\nabla\rho_u^m\| + \|\rho_\delta'^m\|_{\Gamma_1} + \|\rho_\delta^m\|_{\Gamma_1}
\leqslant \sqrt{5}\,\sqrt{c_5\big( c_0 h^{2(r+\frac{1}{2})} + h^{2(r+\frac{1}{2})}\mathcal{N}_1 + \tau^4\mathcal{N}_2 \big)}.
\end{align*}
As for $b_1,b_2,b_3 \geqslant0$, $\sqrt{b_1+b_2 +b_3} \leqslant \sqrt{b_1}+\sqrt{b_2}+\sqrt{b_3}$ and $\sqrt{h^{2(r+\frac{1}{2})}} = h^{r+\frac{1}{2}}$, we obtain
\begin{align*}
\max_{0\leqslant m\leqslant N}\Big( \|\rho_u'^m\| \hspace{-0.07cm}+\hspace{-0.07cm} \|\rho_\theta^m\| \hspace{-0.07cm}+\hspace{-0.07cm} \|\nabla\rho_u^m\| \hspace{-0.07cm}+\hspace{-0.07cm} \|\rho_\delta'^m\|_{\Gamma_1} \hspace{-0.07cm}+\hspace{-0.07cm} \|\rho_\delta^m\|_{\Gamma_1} \Big)
\leqslant \sqrt{5c_5}\,\Big( \sqrt{c_0\hspace{-0.07cm}+\hspace{-0.07cm}\mathcal{N}_1}\; h^{r+\frac{1}{2}} \hspace{-0.07cm}+\hspace{-0.07cm} \sqrt{\mathcal{N}_2}\;\tau^2 \Big).\hspace{-0.07cm}
\end{align*}
Finally, we define the explicit constant
$C_{\star} := \sqrt{5c_5}\; \max\Big\{ \sqrt{c_0+\mathcal{N}_1},\; \sqrt{\mathcal{N}_2} \Big\},$
which depends only on the constants ($c_5$, $c_0$) and the norms of the exact solutions ($\mathcal{N}_1$, $\mathcal{N}_2$), but is independent of $h$ and $\tau$. Consequently,
\begin{align}
\max_{0\leqslant m\leqslant N}\Big( \|\rho_u'^m\| + \|\rho_\theta^m\| + \|\nabla\rho_u^m\| + \|\rho_\delta'^m\|_{\Gamma_1} + \|\rho_\delta^m\|_{\Gamma_1} \Big) \leqslant C_{\star}\big( h^{r+\frac{1}{2}} + \tau^2 \big).\label{estimativa-final1} 
\end{align}
\begin{theorem}[Fully Discrete Convergence]
\label{teo:convergencia}
Assume that the conditions of Theorem \ref{T2.1} and the additional hypotheses in \eqref{hipoteses} are satisfied. Furthermore, assume that the exact solutions $u$, $\theta$, and $\delta$ of the system \eqref{modelo1} satisfy the following regularity conditions, for some $1 \leqslant r \leqslant j$:
\begin{align}
u &\in L^\infty(0,T; H^{r+1}(\Omega)\cap V) \cap H^2(0,T; H^{r+1}(\Omega)\cap V)\cap H^3(0,T; V)\cap H^4(0,T; L^2(\Omega));\label{u-regularidade} \\[4pt]
\theta &\in H^1(0,T; H^{r+1}(\Omega)\cap H^1_0(\Omega)) \cap H^2(0,T; H^1_0(\Omega))\cap H^3(0,T; L^2(\Omega));\label{theta-regularidade} \\[4pt]
\delta &\in H^2(0,T; H^{r+1}(\Gamma_1)) \cap H^4(0,T; L^2(\Gamma_1)).\label{delta-regularidade}
\end{align}
Let $\{U_h^k,P_h^k,\Theta_h^k,D_h^k,Q_h^k\}$ be the solution to the fully discrete problem \eqref{PTD} com $m\in\{1, \ldots, N\}$. Assume that the initial approximations are given by the Ritz projections of the initial data, that is,
$\|\rho_u'^0\| \hspace{-0.07cm}+\hspace{-0.07cm} \|\rho_\theta^0\| \hspace{-0.07cm}+\hspace{-0.07cm} \|\nabla\rho_u^0\| \hspace{-0.07cm}+\hspace{-0.07cm} \|\rho_\delta'^0\|_{\Gamma_1} \hspace{-0.07cm}+\hspace{-0.07cm} \|\rho_\delta^0\|_{\Gamma_1} \leqslant c_0 h^{r+\frac{1}{2}}$ for some $c_0>0$ and $\tau = \mathcal{O}(h^{n/4})$ (where $n=2,3$ is the dimension of $\Omega$). 
Then there exist constants $h_0>0$, $\tau_0>0$, and $C>0$, independent of $h$ and $\tau$, such that, for all $0<h\leqslant h_0$ and $0<\tau\leqslant\tau_0$, the error estimate holds:
\begin{align}
&\max_{0\leqslant k\leqslant N}\Big(\hspace{-0.07cm} \|P_h^k \hspace{-0.07cm}-\hspace{-0.07cm} u'(t_k)\| 
\hspace{-0.07cm}+\hspace{-0.07cm} \|U_h^k \hspace{-0.07cm}-\hspace{-0.07cm} u(t_k)\| 
\hspace{-0.07cm}+\hspace{-0.07cm} \|\Theta_h^k \hspace{-0.07cm}-\hspace{-0.07cm} \theta(t_k)\|
\hspace{-0.07cm}+\hspace{-0.07cm} \|Q_h^k \hspace{-0.07cm}-\hspace{-0.07cm} \delta'(t_k)\|_{\Gamma_1} 
\hspace{-0.07cm}+\hspace{-0.07cm} \|D_h^k \hspace{-0.07cm}-\hspace{-0.07cm} \delta(t_k)\|_{\Gamma_1} \hspace{-0.07cm}\Big) \leqslant C\big( h^{r+\frac{1}{2}} \hspace{-0.07cm}+\hspace{-0.07cm} \tau^2 \big). \label{eq:teorema1}
\end{align}
\end{theorem}
\noindent $\mathbf{ Proof.}$
For each $k$, $k\in\{0,1, \ldots, N\}$, using error decomposition and the triangle inequality, it follows
\begin{align}
&\|P_h^k - u'(t_k)\| 
+ \|U_h^k - u(t_k)\| 
+ \|\Theta_h^k - \theta(t_k)\|
+ \|Q_h^k - \delta'(t_k)\|_{\Gamma_1} 
+ \|D_h^k - \delta(t_k)\|_{\Gamma_1}\nonumber
\\
&
\quad\leqslant 
\|\rho_u'^k\| \hspace{-0.07cm}+\hspace{-0.07cm} \|\rho_\theta^k\| 
\hspace{-0.07cm}+\hspace{-0.07cm} \|\rho_u^k\| 
\hspace{-0.07cm}+\hspace{-0.07cm} \|\rho_\delta'^k\|_{\Gamma_1} 
\hspace{-0.07cm}+\hspace{-0.07cm} \|\rho_\delta^k\|_{\Gamma_1}
\hspace{-0.07cm}+\hspace{-0.07cm}\|\omega_u'^k\| 
\hspace{-0.07cm}+\hspace{-0.07cm} \|\omega_\theta^k\| 
\hspace{-0.07cm}+\hspace{-0.07cm} \|\omega_u^k\| 
\hspace{-0.07cm}+\hspace{-0.07cm} \|\omega_\delta'^k\|_{\Gamma_1} 
\hspace{-0.07cm}+\hspace{-0.07cm} \|\omega_\delta^k\|_{\Gamma_1}\nonumber
\\
&
\quad\leqslant C_0\big(
\|\rho_u'^k\| \hspace{-0.07cm}+\hspace{-0.07cm} \|\rho_\theta^k\| 
\hspace{-0.07cm}+\hspace{-0.07cm} \|\nabla \rho_u^k\| 
\hspace{-0.07cm}+\hspace{-0.07cm} \|\rho_\delta'^k\|_{\Gamma_1} 
\hspace{-0.07cm}+\hspace{-0.07cm} \|\rho_\delta^k\|_{\Gamma_1}\big)
\hspace{-0.07cm}+\hspace{-0.07cm}\|\omega_u'^k\| 
\hspace{-0.07cm}+\hspace{-0.07cm} \|\omega_\theta^k\| 
\hspace{-0.07cm}+\hspace{-0.07cm} \|\omega_u^k\| 
\hspace{-0.07cm}+\hspace{-0.07cm} \|\omega_\delta'^k\|_{\Gamma_1} \hspace{-0.07cm}+\hspace{-0.07cm} \|\omega_\delta^k\|_{\Gamma_1}
\nonumber
\\
&
\quad
\leqslant C_0 C_{\star}\big( h^{r+\frac{1}{2}} + \tau^2 \big) + C_1h^{r+1}
\leqslant (C_0 C_{\star}+C_1) h^{r+\frac{1}{2}} + C_0 C_{\star}\tau^2,\label{estimativa-final2}
\end{align}
where $0<h<1$, implying $h^{r+1}<h^{r+\frac{1}{2}}$ (see \eqref{h});
$C_0>0$ is the constant from the applied Poincaré inequality; $C_{\star}>0$ comes from the estimate result in \eqref{estimativa-final1}, since by \eqref{u-regularidade}-\eqref{delta-regularidade} we have the regularities required throughout the calculations; and $C_1>0$ comes from the projection estimates (see Lemma \ref{lem:proj}).
All these constants are independent of  $h$ and $\tau$. Hence, there exist $h_0>0$ (see \eqref{h}) and $\tau_0$ (see \eqref{tau}) such that, taking $C:=\max \{C_0 C_{\star}+C_1, C_0 C_{\star}\}$ in \eqref{estimativa-final2}, and finally taking the maximum over all $k \in \{0, 1, \ldots, N\} $ in \eqref{estimativa-final2}, the desired result in \eqref{eq:teorema1} follows, with $C$ independent of $h$ and $\tau$. \hfill$\square$
\section{Numerical results}\label{numerical_results}
In this section, we illustrate the convergence order between the numerical and exact solutions established by the error estimate of Theorem \ref{teo:convergencia}.
The numerical examples were implemented in the open-source programming language Julia.
The complete source code, together with the instructions for reproducing the numerical results of this section, is available in the reproducibility repository \cite{ThermoelasticAcoustic}.

The model under study is an extension of Model II presented in \citet{Artigo2025_Adriano}, obtained by adding the thermal component to the system. 
While the previous model considered only the wave equation with acoustic boundary conditions, the present formulation incorporates the heat equation for $\theta$, the bidirectional coupling terms $(\mathbf{a}\cdot \nabla)\theta$ and $(\mathbf{a}\cdot \nabla)u'$, and the nonlocal diffusion coefficient $\beta(\int_\Omega \theta \, dx)$. 
These additional terms introduce new nonlinearities and couplings that demand a more sophisticated numerical treatment.
To assess the impact of this new equation on the results, we use the same input data common to both models, complemented by the input data associated with the thermal part.
More precisely, we consider the following manufactured solutions for the nonhomogeneous version of the model under study \eqref{modelo1}:
\begin{align}\label{sol_manuf}
u(x,y,t) = (x^b-x)(y^b-1)(t^2+4),\quad
\theta(x,y,t) = (x^b-x)(y^b-y)(t^2+4),
\end{align}
with $\delta(x,t)$ constructed from the identity $\delta'(x,t) = -u_y(x,0,t) + g(x,u'(x,0,t))$, for $x\in]0,1[$ and $t\in[0,T]$, where $T$ denotes the final time.
Here, $\Omega = \,]0,1[\,\times\,]0,1[\,$, $\Gamma_1$ denotes the lower edge of the square, where the acoustic equation acts, and $\Gamma_0$ denotes the remainder of the boundary, so that $\Gamma = \Gamma_0\cup\Gamma_1$ represents the complete boundary of the domain.
The remaining input data are as follows:
\begin{align*}
&
\delta_0(x) = \sin(\pi x),\quad
\mu(t) = 1 + e^{-t},\quad
\beta(s) = 1 + e^{-s^2},\quad
g(x,s) = (1 + e^{-x^2})(2s+\sin(s)),\quad
f(s) = s|s|^3,
\\
&
f_1 = f_2 = f_3 = 1,\quad
\mathbf{a} = (1,1),\quad
T = 1.
\end{align*}

With these data, and taking $f_1,\,f_2,\,f_3$ as constants for simplicity, the nonhomogeneous version of the model \eqref{modelo1} on the unit square, already incorporating the source terms required for \eqref{sol_manuf} to satisfy the system, is given by
\begin{align*}
& u^{\prime\prime}(x,y,t) 
- \mu(t)\Delta u(x,y,t) 
+ f\big(u(x,y,t)\big)
+ (\mathbf{a}\cdot\nabla)\theta(x,y,t) 
= \mathbb{f}_1(x,y,t),
\\
& \theta^{\prime}(x,y,t) 
- \beta\Big(\int_\Omega\theta(t)d\Omega\Big)\Delta\theta(x,y,t)
+ (\mathbf{a}\cdot\nabla)u^\prime(x,y,t)
= \mathbb{f}_2(x,y,t),
\\[5pt]
& f_1 \delta^{\prime\prime}(x,t)
+ f_2 \delta^\prime(x,t)
+ f_3 \delta(x,t)
= \mathbb{f}_3(x,t),
\\
& -\frac{\partial u}{\partial y}(x,0,t)
= \delta^\prime(x,t)
- g\big(x,u^\prime(x,0,t)\big),
\\
& 
u = 0\; \text{in }\Gamma_0,\; 
\theta = 0\;\text{in }\Gamma,\;
u(0) = u_0, \; 
u^\prime(0)=u_1, \; 
\theta_0(0)=\theta_0,\; 
\delta(0) = \delta_0,\;
\delta^\prime(x,0) = -\frac{\partial u_0}{\partial y}(x,0) + g\big(x,u_1(x,0)\big).
\end{align*}

Theorem \ref{teo:convergencia} establishes the error estimate
\begin{align}\label{result_estimative}
  \Vert P_h^k-u^\prime(t_k)\Vert
+ \Vert U_h^k-u(t_k)\Vert
+ \Vert \Theta_h^k-\theta(t_k)\Vert
+ \Vert Q_h^k-\delta^\prime(t_k)\Vert_{\Gamma_1}
+ \Vert D_h^k-\delta(t_k)\Vert_{\Gamma_1}
\leqslant
c(h^{r+\frac{1}{2}} + \tau^2),
\end{align}
valid for all $k\in\{0,1,\dots,N\}$, where $N$ is the integer part of $T/\tau$, $r$ is the polynomial degree of the basis functions generating the approximation subspace and $\tau$ is the time step of the temporal discretization.
This estimate corresponds to a convergence order of $\mathcal{O}(h^{r+\frac{1}{2}} + \tau^2)$, which is suboptimal in space and optimal in time.

This loss of order in space originates from the coupling between the domain and boundary equations, as well as from the solution regularity.
As discussed in \cite{Artigo2015_Kashiwabara,Artigo2025_Adriano}, requiring higher regularity of the solution makes it possible to overcome this loss and achieve the optimal convergence order in the $L^2$ norm, namely $\mathcal{O}(h^{r+1})$.
In this work, as in \cite{Artigo2025_Adriano}, we numerically illustrate this order gain by means of the parameter $b$ in \eqref{sol_manuf}, which controls the regularity of the manufactured solution and is chosen according to the polynomial degree of the basis used, as detailed in Table \ref{tab:regularities}, where
\begin{align}\label{def:Hm_Woo}
H^m:= \left\{ v \in L^2(\Omega) : D^\alpha v \in L^2(\Omega)\;\forall |\alpha| \leqslant m \right\};\quad
W_\infty^m := \left\{ v \in L^\infty(\Omega) : D^\alpha v \in L^\infty(\Omega)\;\forall |\alpha| \leqslant m \right\}.
\end{align}

Table \ref{tab:convergence_study__spatial} presents the spatial convergence study resulting from these combinations, obtained on uniform meshes of quadrilateral finite elements with $N_x$ partitions per side of the unit square; this value is reported in the column $N_x$ of the table.
To isolate the spatial error, we fix the temporal discretizations with a sufficiently small time step, namely $\tau = 2^{-13}$.
To highlight the error behavior of each variable separately, which is more informative than the sum, since it can reveal differences in convergence order among them, we take the temporal maximum of each norm composing the estimate \eqref{result_estimative}, defining
\begin{align}\label{def:LooL2_discrete}
\mathrm{error}(u') 
= \max_{k\in\{0,1,\dots,N\}}
\Vert P_h^k-u^\prime(t_k)\Vert,
\end{align}
with $\mathrm{error}(u)$, $\mathrm{error}(\theta)$, $\mathrm{error}(\delta')$, and $\mathrm{error}(\delta)$ defined analogously.

As shown in Table \ref{tab:convergence_study__spatial}, the minimum convergence order established by Theorem \ref{teo:convergencia} is satisfied asymptotically for all combinations of $b$ and polynomial degree considered.
The most challenging case occurs for the cubic basis: on the coarsest meshes, the observed rates for variables $\delta'$ and $\delta$ fall short of the theoretical order.
In the optimal case ($b = 4.4$), for example, $\mathrm{error}(\delta')$ and $\mathrm{error}(\delta)$ exhibit rates of $2.45$ and $2.52$, respectively, between $N_x = 4$ and $8$, against the expected order $\mathcal{O}(h^4)$.
However, as the mesh is refined, the rates quickly converge to the theoretical order.

Turning to the temporal convergence order study, Table \ref{tab:convergence_study__temporal} shows the order $\mathcal{O}(\tau^2)$ for each variable.
In this case, we use the linear Lagrange basis, a spatial mesh with $512$ partitions in each direction, and a manufactured-solution regularity parameter $b = 2.4$.

In comparison with the numerical results presented in \cite{Artigo2025_Adriano}, it is worth noting the increase in computational cost arising from the thermal coupling.
In the previous model, without this component, the nonlinear algebraic system solved at each time step has the form
\begin{align*}
sparse\_matrix*X + nonlinear\_vector(X) = vector.
\end{align*}
For the current model, which incorporates the thermal part, the system takes on a block structure, namely,
\begin{align*}
\begin{bmatrix}
sparse\_matrix_{1,1} & sparse\_matrix_{1,2}\\
sparse\_matrix_{2,1} & sparse\_matrix_{2,2}
\end{bmatrix}*
\begin{bmatrix}
X_1\\
X_2
\end{bmatrix}
+\begin{bmatrix}
nonlinear\_vector_1(X_1)\\
nonlinear\_vector_2(X_2)
\end{bmatrix}
=\begin{bmatrix}
vector_1\\
vector_2
\end{bmatrix}.
\end{align*}
Moreover, since Newton's method is employed to solve the nonlinear algebraic systems, the term $nonlinear\_vector_2(X_2)$, associated with the component $\beta(\int_\Omega\theta d\Omega)\Delta\theta$, causes the Jacobian matrix to take the form $A + \mathbb{u}\mathbb{v}^t$, where $A$ is a sparse matrix and $\mathbb{u}$, $\mathbb{v}$ are column vectors defining the rank-one matrix $\mathbb{u}\mathbb{v}^t$.
In this case, the Sherman--Morrison formula is applied to solve the resulting linear system at each Newton iteration, thereby avoiding the explicit assembly of the rank-one matrix $\mathbb{u}\mathbb{v}^t$, a dense block that would otherwise be added to the sparse matrix $A$, at the cost of solving two linear systems with the same sparse matrix $A$.
Details on the algebraic systems and on the computation of the Jacobian matrices for both models can be found in the package documentation at
\cite{ThermoelasticAcoustic} and \cite{WaveAcoustics}.

Motivated by the increase in computational cost reported above, we numerically investigated linearized variants of the scheme, obtained by replacing the nonlinear terms with backward approximations (extrapolations built from previous time steps), so that the unknown variable no longer appears nonlinearly in the resulting system.
When applying this strategy to the model without the thermal part of \cite{Artigo2025_Adriano}, we found that linearizing the boundary term $g(x,u^\prime(x,0,t))$ compromises the unconditional convergence of the method, requiring $\tau$ smaller than $h$ for convergence.
By preserving the nonlinearity only in this term and linearizing the remaining terms, including the coupling between the wave and heat equations in the wave equation, which allows the system to be decoupled, we obtained a computationally cheaper scheme, without loss of convergence order or any additional restriction between $h$ and $\tau$.
We emphasize that we did not carry out a theoretical error-estimate analysis for these alternative schemes.
This is a purely numerical investigation, the implementation details and the corresponding convergence tests are available in the package documentation at
\cite{ThermoelasticAcoustic} and \cite{WaveAcoustics}.

\begin{table}[width=.9\linewidth,cols=4,pos=h]
\caption{Chosen values of $b$ for \eqref{sol_manuf}, corresponding regularity, and expected order of convergence in the $L^2$ norm, with $H^m$ and $W_\infty^m$ defined as in \eqref{def:Hm_Woo}.}
\label{tab:regularities}
\begin{tabular*}{\tblwidth}{@{} LLLL@{} }
\toprule
$b$-value & $b$-options           & Regularity                                                           & $L^2$ error
\\ \midrule
$1.76$    & $\frac{3}{2}<b<2$     & $u,\theta\in H^2$, $u,\theta\notin H^3$, $u,\theta\notin W_\infty^2$ & $h^{1.5}$\, (Lagrange linear)
\\
$2.4$     & $2<b\leqslant \frac{5}{2}$ & $u,\theta\in H^2$, $u,\theta\notin H^3$, $u,\theta\in W_\infty^2$    & $h^2$ \quad (Lagrange linear)
\\
$2.58$    & $\frac{5}{2}<b<3$     & $u,\theta\in H^3$, $u,\theta\notin H^4$, $u,\theta\notin W_\infty^3$ & $h^{2.5}$ \;(Lagrange quadratic)
\\
$3.4$     & $3<b\leqslant \frac{7}{2}$ & $u,\theta\in H^3$, $u,\theta\notin H^4$, $u,\theta\in W_\infty^3$    & $h^3$ \quad(Lagrange quadratic)
\\
$3.51$    & $\frac{7}{2}<b<4$     & $u,\theta\in H^4$, $u,\theta\notin H^5$, $u,\theta\notin W_\infty^4$ & $h^{3.5}$ \;(Lagrange cubic)
\\
$4.4$     & $4<b\leqslant \frac{9}{2}$ & $u,\theta\in H^4$, $u,\theta\notin H^5$, $u,\theta\in W_\infty^4$    & $h^4$ \;\; (Lagrange cubic)
\\ \bottomrule
\end{tabular*}
\end{table}
\begin{table}[width=.9\linewidth,cols=11,pos=h]
\caption{
Spatial convergence order study with time step $\tau = 2^{-13}$ for each combination of $b$-value and polynomial basis specified in Table \ref{tab:regularities}. 
Error columns follow the discrete norm defined in \eqref{def:LooL2_discrete}.
}
\label{tab:convergence_study__spatial}
\begin{tabular*}{\tblwidth}{@{} LLLLLLLLLLL@{} }
\toprule
$N_x$ & error($u'$)       & rate   & error($u$)       & rate   & error($\theta$)  & rate   & error($\delta'$)       & rate   & error($\delta$)       & rate  \\
\midrule
\multicolumn{11}{l}{Lagrange linear, sub-optimal case, $b$-value $1.76$:  $\mathcal{O}(h^{1.5})$}\\
\midrule
$2^2$ & $1.07\text{e-}2$ & $-$    & $3.08\text{e-}2$ & $-$    & $9.00\text{e-}3$ & $-$    & $6.68\text{e-}2$ & $-$    & $4.64\text{e-}2$ & $-$   \\
$2^3$ & $2.75\text{e-}3$ & $1.97$ & $7.90\text{e-}3$ & $1.96$ & $2.30\text{e-}3$ & $1.97$ & $1.51\text{e-}2$ & $2.14$ & $1.07\text{e-}2$ & $2.12$\\
$2^4$ & $7.08\text{e-}4$ & $1.96$ & $2.04\text{e-}3$ & $1.95$ & $5.87\text{e-}4$ & $1.97$ & $3.56\text{e-}3$ & $2.09$ & $2.56\text{e-}3$ & $2.06$\\
$2^5$ & $1.89\text{e-}4$ & $1.91$ & $5.38\text{e-}4$ & $1.92$ & $1.50\text{e-}4$ & $1.97$ & $8.56\text{e-}4$ & $2.06$ & $6.28\text{e-}4$ & $2.03$\\
\midrule
\multicolumn{11}{l}{Lagrange linear, optimal case, $b$-value $2.4$:  $\mathcal{O}(h^{2})$}\\
\midrule
$2^2$ & $1.92\text{e-}2$ & $-$    & $4.51\text{e-}2$ & $-$    & $1.86\text{e-}2$ & $-$    & $6.68\text{e-}2$ & $-$    & $4.95\text{e-}2$ & $-$   \\
$2^3$ & $4.95\text{e-}3$ & $1.96$ & $1.11\text{e-}2$ & $2.02$ & $4.62\text{e-}3$ & $2.01$ & $1.50\text{e-}2$ & $2.15$ & $1.15\text{e-}2$ & $2.11$\\
$2^4$ & $1.25\text{e-}3$ & $1.99$ & $2.77\text{e-}3$ & $2.01$ & $1.15\text{e-}3$ & $2.00$ & $3.61\text{e-}3$ & $2.06$ & $2.79\text{e-}3$ & $2.04$\\
$2^5$ & $3.13\text{e-}4$ & $2.00$ & $6.92\text{e-}4$ & $2.00$ & $2.88\text{e-}4$ & $2.00$ & $8.88\text{e-}4$ & $2.02$ & $6.92\text{e-}4$ & $2.01$\\
\midrule
\multicolumn{11}{l}{Lagrange quadratic, sub-optimal case, $b$-value $2.58$:  $\mathcal{O}(h^{2.5})$}\\
\midrule
$2^2$ & $5.91\text{e-}4$ & $-$    & $1.45\text{e-}3$ & $-$    & $6.13\text{e-}4$ & $-$    & $4.14\text{e-}3$ & $-$    & $3.75\text{e-}3$ & $-$   \\
$2^3$ & $8.34\text{e-}5$ & $2.83$ & $1.93\text{e-}4$ & $2.91$ & $8.14\text{e-}5$ & $2.91$ & $6.15\text{e-}4$ & $2.75$ & $5.23\text{e-}4$ & $2.84$\\
$2^4$ & $1.06\text{e-}5$ & $2.97$ & $2.53\text{e-}5$ & $2.93$ & $1.07\text{e-}5$ & $2.93$ & $8.37\text{e-}5$ & $2.88$ & $6.84\text{e-}5$ & $2.93$\\
$2^5$ & $1.46\text{e-}6$ & $2.86$ & $3.31\text{e-}6$ & $2.94$ & $1.39\text{e-}6$ & $2.94$ & $1.10\text{e-}5$ & $2.93$ & $8.77\text{e-}6$ & $2.96$\\
\midrule
\multicolumn{11}{l}{Lagrange quadratic, optimal case, $b$-value $3.4$:  $\mathcal{O}(h^{3})$}\\
\midrule
$2^2$ & $1.49\text{e-}3$ & $-$    & $3.34\text{e-}3$ & $-$    & $1.63\text{e-}3$ & $-$    & $5.06\text{e-}3$ & $-$    & $3.70\text{e-}3$ & $-$   \\
$2^3$ & $1.82\text{e-}4$ & $3.03$ & $4.18\text{e-}4$ & $3.00$ & $2.04\text{e-}4$ & $3.00$ & $6.40\text{e-}4$ & $2.98$ & $4.84\text{e-}4$ & $2.93$\\
$2^4$ & $2.32\text{e-}5$ & $2.97$ & $5.23\text{e-}5$ & $3.00$ & $2.55\text{e-}5$ & $3.00$ & $8.32\text{e-}5$ & $2.94$ & $6.16\text{e-}5$ & $2.98$\\
$2^5$ & $2.90\text{e-}6$ & $3.00$ & $6.53\text{e-}6$ & $3.00$ & $3.19\text{e-}6$ & $3.00$ & $1.06\text{e-}5$ & $2.97$ & $7.74\text{e-}6$ & $2.99$\\
\midrule
\multicolumn{11}{l}{Lagrange cubic, sub-optimal case, $b$-value $3.51$:  $\mathcal{O}(h^{3.5})$}\\
\midrule
$2^2$ & $3.05\text{e-}5$ & $-$    & $7.14\text{e-}5$ & $-$    & $3.53\text{e-}5$ & $-$    & $3.97\text{e-}4$ & $-$    & $1.35\text{e-}4$ & $-$   \\
$2^3$ & $2.06\text{e-}6$ & $3.89$ & $4.89\text{e-}6$ & $3.87$ & $2.42\text{e-}6$ & $3.87$ & $4.03\text{e-}5$ & $3.30$ & $1.13\text{e-}5$ & $3.59$\\
$2^4$ & $1.37\text{e-}7$ & $3.91$ & $3.29\text{e-}7$ & $3.89$ & $1.63\text{e-}7$ & $3.89$ & $2.54\text{e-}6$ & $3.99$ & $7.13\text{e-}7$ & $3.98$\\
$2^5$ & $1.01\text{e-}8$ & $3.77$ & $2.19\text{e-}8$ & $3.91$ & $1.08\text{e-}8$ & $3.91$ & $1.49\text{e-}7$ & $4.09$ & $4.30\text{e-}8$ & $4.05$\\
\midrule
\multicolumn{11}{l}{Lagrange cubic, optimal case, $b$-value $4.4$:  $\mathcal{O}(h^{4})$}\\
\midrule
$2^2$ & $1.04\text{e-}4$ & $-$    & $2.21\text{e-}4$ & $-$    & $1.19\text{e-}4$ & $-$    & $5.48\text{e-}4$ & $-$    & $1.38\text{e-}4$ & $-$   \\
$2^3$ & $6.59\text{e-}6$ & $3.97$ & $1.40\text{e-}5$ & $3.98$ & $7.52\text{e-}6$ & $3.99$ & $1.00\text{e-}4$ & $2.45$ & $2.39\text{e-}5$ & $2.52$\\
$2^4$ & $3.79\text{e-}7$ & $4.12$ & $8.76\text{e-}7$ & $4.00$ & $4.70\text{e-}7$ & $4.00$ & $7.04\text{e-}6$ & $3.83$ & $1.69\text{e-}6$ & $3.83$\\
$2^5$ & $2.39\text{e-}8$ & $3.99$ & $5.48\text{e-}8$ & $4.00$ & $2.94\text{e-}8$ & $4.00$ & $4.16\text{e-}7$ & $4.08$ & $9.99\text{e-}8$ & $4.08$\\
\bottomrule
\end{tabular*}
\end{table}
\begin{table}[width=.9\linewidth,cols=11,pos=h]
\caption{
Temporal convergence order study with $N_x = 2^{9}$, $b$-value $2.4$ (see Table~\ref{tab:regularities}), and linear Lagrange basis.
}
\label{tab:convergence_study__temporal}
\begin{tabular*}{\tblwidth}{@{} LLLLLLLLLLL@{} }
\toprule
$\tau$   & error($u'$)       & rate   & error($u$)       & rate   & error($\theta$)  & rate   & error($\delta'$)       & rate   & error($\delta$)       & rate  \\
\midrule
$2^{-2}$ & $1.09\text{e-}2$ & $-$    & $3.47\text{e-}3$ & $-$    & $1.45\text{e-}3$ & $-$    & $7.12\text{e-}3$ & $-$    & $5.41\text{e-}3$ & $-$   \\
$2^{-3}$ & $2.72\text{e-}3$ & $2.00$ & $8.95\text{e-}4$ & $1.96$ & $3.00\text{e-}4$ & $2.27$ & $1.79\text{e-}3$ & $1.99$ & $1.36\text{e-}3$ & $1.99$\\
$2^{-4}$ & $6.76\text{e-}4$ & $2.01$ & $2.27\text{e-}4$ & $1.98$ & $6.17\text{e-}5$ & $2.28$ & $4.50\text{e-}4$ & $2.00$ & $3.41\text{e-}4$ & $2.00$\\
$2^{-5}$ & $1.68\text{e-}4$ & $2.01$ & $5.80\text{e-}5$ & $1.97$ & $1.60\text{e-}5$ & $1.95$ & $1.13\text{e-}4$ & $2.00$ & $8.54\text{e-}5$ & $2.00$\\
\bottomrule
\end{tabular*}
\end{table}

\section{Conclusion}
In this work, we present the first numerical analysis in the $L^2$-norm for the nonlinear thermoelastic system with acoustic boundary conditions introduced by Braz e Silva, Clark, and Frota (2017) in \cite{Artigo2017_SilvaClarkFrota}.
Using the Crank-Nicolson scheme in time and  finite element methods with basis of degree $r \geqslant 1$ in space, we establish the error estimate
$
\mathcal{O}(h^{r+\frac{1}{2}} + \tau^2)
$
for all variables of the system. The loss of half a power in $u$, arising from the domain-boundary coupling, propagates to $\theta$ and $\delta$ through the thermoelastic and boundary couplings. 
The nonlinearity  $f(u) = \lambda |u|^d u$ was handled via a bootstrap argument, ensuring the uniform boundedness of $\|U_h^n\|_{L^\infty}$ under the condition $\tau = \mathcal{O}(h^{n/4})$. 

Numerical experiments using finite element spaces with Lagrange basis of degrees $r = 1, 2, 3$ corroborate the estimates and illustrate the recovery of the optimal order when the exact solution possesses additional $L^\infty$ regularity.
The experiments were implemented in the open-source programming language Julia, with source code available in the reproducibility repository \cite{ThermoelasticAcoustic}, allowing full and independent verification of the numerical results presented here.

In summary, this work not only extends $L^2$-norm error estimate results to a more general class of thermoelastic systems but also provides a solid foundation for future investigations. In addition to extending the estimation and numerical results of \cite{Artigo2018_Adriano, Artigo2021_Adriano, Artigo2022_Adriano, Artigo2025_Adriano}, \cite{Artigo1976_Baker}, \cite{Artigo2015_Kashiwabara}, and \cite{Artigo2020_HippKovacs_L2}.

\printcredits 

\section*{Acknowledgments}
The following authors obtained funding: Adriano A. Alc\^{a}ntara was supported by a post-doc grant from CNPq-Brazil (N. 151903/2024-7), 
Bruno A. Carmo was supported by a post-doc grant from CNPq-Brazil (N. 151197/2025-3),   Carlos M. Guzmán is partially supported by FAPERJ, process (N. E-26/SEI-260003/003317/2022), and CNPq, grant (N. 312032/2022-6). Furthermore, Juan B. Límaco is partially supported by CNPq.

\section*{Data availability}
The complete source code is provided in a reproducibility repository.

\section*{Declaration of competing interests}
The authors declare that they have no conflict of interest.

\section*{Declaration of generative AI use}
During the preparation of this work and implementation of the Julia packages \cite{ThermoelasticAcoustic,WaveAcoustics} the author B. A. Carmo used Claude (Antropic, free browser-based chat) in order to: (i) support translation, grammatical review, and refinement of the scientific writing; (ii) review selected code sections to identify and fix bugs; (iii) suggest more idiomatic or efficient Julia implementations; (iv) refine docstrings and documentation text. No IDE autocomplete or local AI agent was used. After using this tool/service, the author reviewed and edited the content as needed and take full responsibility for the content of the published article and packages.

During the preparation of this manuscript, the author Adriano A. Alcantara used ChatGPT (OpenAI, free browser-based version) to assist with English language editing, grammar correction, translation, and refinement of the scientific writing. After using this tool, the author reviewed and edited the content as needed and takes full responsibility for the content of the published article.

\bibliographystyle{cas-model2-names}

\bibliography{cas-refs}
\end{document}